\documentclass{amsart}
\usepackage{amssymb,pb-diagram}

\title{Affine Weyl groups in $K$-theory and representation theory}

\author{Cristian Lenart and Alexander Postnikov}

\address{Department of Mathematics and Statistics, State University of New York, Albany, NY 12222}
\email{lenart@csc.albany.edu}

\address{Department of Mathematics, M.I.T., Cambridge, MA 02139}
\email{apost@math.mit.edu}


\thanks{Cristian Lenart was supported by SUNY Albany Faculty Research Award 1032354}
\thanks{Alexander Postnikov was supported by National Science Foundation 
grant DMS-0201494 and by Alfred P.\ Sloan Foundation research fellowship}


\date{original version: September~12, 2003; updated: May~27, 2004, June~1, 2005}

\numberwithin{equation}{section}

\theoremstyle{plain}
\newtheorem{theorem}{Theorem}[section]
\newtheorem{proposition}[theorem]{Proposition}
\newtheorem{lemma}[theorem]{Lemma}
\newtheorem{corollary}[theorem]{Corollary}

\theoremstyle{definition}
\newtheorem{definition}[theorem]{Definition}
\newtheorem{conjecture}[theorem]{Conjecture}
\newtheorem{example}[theorem]{Example}

\theoremstyle{remark}
\newtheorem{remark}[theorem]{Remark}

\def\R{\mathbb{R}}
\def\Z{\mathbb{Z}}
\def\C{\mathbb{C}}
\def\Q{\mathbb{Q}}

\def\O{\mathcal{O}}
\def\SL{\mathit{SL}}
\def\L{\mathcal{L}}
\def\I{\mathcal{I}}
\def\weight{\mathrm{weight}}
\def\Waff{W_{\mathrm{aff}}}
\def\ds{\displaystyle}
\def\h{\mathfrak{h}}
\def\hR{\mathfrak{h}^*_\mathbb{R}}
\def\hH{\mathcal{\hat H}}

\begin{document}
\bibliographystyle{plain}

\begin{abstract}
We give an explicit combinatorial Chevalley-type formula for the equivariant
$K$-theory of generalized flag varieties $G/P$ which is a direct generalization
of the classical Chevalley formula.  Our formula implies a simple combinatorial
model for the characters of the irreducible representations of $G$ and, more
generally, for the Demazure characters. This model can be viewed as a discrete
counterpart of the Littelmann path model, and has several advantages. Our
construction is given in terms of a certain $R$-matrix, that is, a collection
of operators satisfying the Yang-Baxter equation.  It reduces to combinatorics
of decompositions in the affine Weyl group and enumeration of saturated chains
in the Bruhat order on the (nonaffine) Weyl group.  Our model easily implies
several symmetries of the coefficients in the Chevalley-type formula. We also
derive a simple formula for multiplying an arbitrary Schubert class by a
divisor class, as well as a dual Chevalley-type formula.  The paper contains
other applications and examples.  
\end{abstract}

\maketitle

\vspace{-0.9cm}

\tableofcontents

\vspace{-0.9cm}

\section{Introduction}
\label{sec:intro}

The Chevalley formula~\cite{Chev} from Schubert calculus expresses 
the products of the classes of Schubert varieties with the classes 
of certain line bundles in the cohomology ring of the
generalized flag variety $G/B$, where $G$ is a complex semisimple Lie group and
$B$ is a Borel subgroup. 
This formula implies a rule for products of divisor classes
with arbitrary Schubert classes, known as Monk's rule in type~$A$.
Fulton and Lascoux~\cite{FL} extended the Chevalley formula to
the equivariant Grothendieck ring $K_T(SL_n/B)$ of the classical flag variety,
using combinatorics of Young tableaux. Other 
Chevalley-type and Monk-type formulas in $K(SL_n/B)$ were given in \cite{Len}. Pittie and Ram~\cite{PR} extended the
Chevalley formula to the equivariant Grothendieck ring $K_T(G/B)$ using
LS-paths, which are special cases of Littelmann paths.  However, the Pittie-Ram
formula is often hard to use for explicit calculations.  It works for dominant
weights only and involves some nontrivial recursive procedures.  In this
article, we present a simple nonrecursive combinatorial Chevalley-type formula
for products in the equivariant Grothendieck ring $K_T(G/P)$, where $P$ is a
parabolic subgroup in $G$.  Our formula implies a nonnegative combinatorial
model for the characters of the irreducible representations of $G$ and for 
the Demazure
characters.  This model is more efficient computationally than other known
models for characters, such as the Littelmann path model.  Our formula easily
explains two symmetries of Chevalley coefficients in the equivariant
$K$-theory, clarifies their connection with a Monk-type formula in
this ring, and
implies positivity (or negativity) of these coefficients.  One of these
symmetries was earlier derived by Brion~\cite{Bri} using a nontrivial geometric
argument.  Our formula is based on a collection of operators that satisfy the
Yang-Baxter equation.  Its proof is completely elementary.  It does not rely on
any geometric arguments, and it just uses combinatorics of the affine Weyl group
and some algebraic manipulations with $R$-matrices and Demazure operators.

Littelmann paths give a model for the characters of the irreducible
representations $V_\lambda$ of $G$.  Littelmann~\cite{Li1,Li2} showed that the
characters can be described by counting certain continuous paths in $\hR$.
These paths are constructed recursively starting with an initial one, by using
certain operators acting on them, which are known as root operators. By making
specific choices for the initial path, one can obtain special cases which are
described combinatorially. One such class of paths, corresponding to a straight
line initial path, is known as the class of Lakshmibai-Seshadri paths
(LS-paths). These paths were introduced before Littelmann's work, in the
context of standard monomial theory~\cite{LS1}.  They have a nonrecursive
characterization in terms of the Bruhat order on the quotient $W/W_\lambda$ of
the corresponding Weyl group $W$ modulo the stabilizer $W_\lambda$ of
$\lambda$.  Recently, Gaussent and Littelmann~\cite{GaLi}, motivated by the
study of Mirkovi\'c-Vilonen cycles, defined another combinatorial model for the
irreducible characters of a complex semisimple Lie group.  This model is based
on LS-galleries, which are certain sequences of faces of alcoves for the
corresponding affine Weyl group.

A geometric application of LS-paths was given by Pittie and Ram~\cite{PR}, who
used them to derive a Chevalley-type multiplication formula in the
$T$-equivariant $K$-theory of the generalized flag variety $G/B$.  Let
$K_T(G/B)$ be the Grothendieck ring of $T$-equivariant coherent sheaves on
$G/B$.  According to Kostant and Kumar~\cite{KK}, the ring $K_T(G/B)$ is a free
module over the representation ring $R(T)$ of the maximal torus, with basis
given by the classes $[\O_{X_w}]$, $w\in W$, of structure sheaves of Schubert
varieties.  Pittie and Ram showed that the basis expansion of the product of
$[\O_{X_w}]$ with the class $[\L_\lambda]$ of a negative line bundle (corresponding to the character of $B$ determined by the antidominant weight $-\lambda$) can be expressed as a nonnegative sum over certain special LS-paths.
The fact that the product in the Pittie-Ram formula expands as a nonnegative
linear combination was also explained geometrically by Brion~\cite{Bri} and
Mathieu~\cite{Mat}.  The coefficients in the Pittie-Ram
formula were identified as certain characters by Lakshmibai and
Littelmann~\cite{LL} using geometry.  Littelmann and Seshadri~\cite{LS3} showed
that the Pittie-Ram formula is a consequence of standard monomial
theory~\cite{LLM,LS1,Li3}, and, furthermore, that it is almost equivalent to
standard monomial theory.

In this paper, we present an alternative simple Chevalley-type 
formula\footnote{Notational remark:  We call a rule for
$[\L_\lambda]\cdot [\O_{X_w}]$ a {\it Chevalley-type formula} and
use the term {\it Monk-type formula}
for a rule for products $[\O_{X_{w_\circ s_i}}]\cdot [\O_{X_w}]$ of divisor classes 
$[\O_{X_{w_\circ s_i}}]$ with arbitrary classes $[\O_{X_w}]$. 
The term {\it Pieri-type formula} refers to multiplication with the 
special Schubert classes pulled back from a Grassmannian.}
for the product of $[\O_{X_w}]$ and $[\L_\lambda]$ in the equivariant 
Grothendieck ring
$K_T(G/P)$.  The formula is based on enumerating certain {\it saturated
chains\/} in the Bruhat order on the corresponding Weyl group $W$. This
enumeration is determined by an {\it alcove path}, which is a sequence of
adjacent alcoves for the affine Weyl group $\Waff$ of the Langland's dual group
$G^\vee$. Alcove paths correspond to decompositions of elements in the affine
Weyl group into products of generators.  Our Chevalley-type formula is
conveniently formulated in terms of a certain {\it R-matrix}, that is, in terms
of a collection of operators satisfying the {\it Yang-Baxter equation}.  We
express the operator $E^\lambda$ of multiplication by the class of a line
bundle as a composition $R^{[\lambda]}$ of elements of the
$R$-matrix given by a certain alcove path.  In order to prove
the formula, we simply verify that the operators $R^{[\lambda]}$ satisfy the
same commutation relations with the elementary Demazure operators $T_i$ as the
operators $E^\lambda$.

Our equivariant $K$-theory Chevalley formula has the following nice features, including several ones developed or to be developed in subsequent publications.
\begin{itemize}
\item  
The formula works for line bundles corresponding to arbitrary 
weights.  The Pittie-Ram formula works for dominant weights only.
Note that several applications require working with nondominant weights.
\item
The formula is equally simple for all weights (regular and nonregular,
dominant and nondominant) and is a direct generalization of the classical 
Chevalley formula.
The Pittie-Ram formula and standard monomial theory 
require Deodhar's lift operators
$W/W_\lambda\to W$ from cosets modulo $W_\lambda$,
which are defined by a nontrivial recursive procedure
\cite{Deo2}. 
In our construction, no lift operators are needed, since we are working in $W$.
\item 
Our formula easily implies a Monk-type formula for products of the 
classes $[\O_{X_w}]$ with divisor classes.  Indeed, the special classes 
are expressed in terms of the classes of line bundles $\L_{-\omega_i}$, where $\omega_i$ denotes a fundamental weight.
It is more difficult to apply the Pittie-Ram formula for this computation,
because the latter formula makes sense for dominant weights only.
Our formula implies a Pieri-type formula in $K(SL_n/B)$ \cite{lasptf} 
for multiplying arbitrary Schubert classes with certain special Schubert 
classes pulled back from a Grassmannian that are indexed by cycles.  No other such formulas based
on other models are known.
\item 
Our formula exposes
 explicitly why products with the classes of line bundles $\L_\lambda$ result in nonnegative coefficients if $\lambda$ is a dominant weight, and in coefficients
with alternating signs if $\lambda$ is antidominant.  The present model facilitates the study of certain symmetries of the Chevalley 
coefficients in equivariant $K$-theory, which is not easily carried 
out based on other methods.
\item
Our formula immediately implies the dual Chevalley-type formula for
products of $[\L_\lambda]$ with elements of the dual basis 
to $\{[\O_{X_w}]\mid w\in W\}$.
\item
The independence of our formula from the choice involved in our model (i.e., the choice of an alcove path) 
follows from the fact that the $R$-matrices used in the construction 
satisfy the Yang-Baxter equation.
No such simple explanation is available for other models.  
\item
The proof of our formula is completely algebraic/combinatorial. 
\item 
For dominant weights $\lambda$, our formula implies a simple 
combinatorial model for the characters of the irreducible 
representations $V_\lambda$ and for the Demazure characters 
$ch(V_{\lambda,w})$. In fact, in \cite{LP} we develop our model entirely within the representation theory of complex symmetrizable Kac-Moody algebras; in this context, we derive an explicit Littlewood-Richardson
rule for decomposing tensor products of irreducible representations and describe the corresponding crystal graph structures.
\item Our model enables a far-reaching generalization of the combinatorics of
Young tableaux.  We will discuss these issues in forthcoming 
publication(s).  For example, we will discuss a generalization of 
Sch\"utzenberger's {\em evacuation procedure} for tableaux, 
which describes the action of the fundamental involution on canonical 
bases for irreducible $U_q({\mathfrak g})$-modules.
\end{itemize}

As a preview of our main result, let us present here a formula for the
product $[\L_\lambda]\cdot [\O_{X_w}]$ of classes in the usual (nonequivariant)
Grothendieck  ring\footnote{The ring $K(G/B)$ is not related to Russian
security services.} $K(G/B)$.   
Let $\mathcal{A}$ be the affine Coxeter arrangement for the Langland's 
dual group $G^\vee$.  The regions of $\mathcal{A}$, called alcoves, correspond
to the elements of the affine Weyl group $\Waff$.
Fix a weight $\lambda$.  Let $\pi(t)$ be a
continuous path in $\hR$ that connects a point $\pi(0)$ inside the fundamental
alcove with the point $\pi(1)=\pi(0)-\lambda$.  Assume that $\pi(t)$ 
does not pass through pairwise intersections of hyperplanes in $\mathcal{A}$. 
As $t$ changes from 0 to 1, the
path $\pi(t)$ crosses the hyperplanes $H_1,\dots,H_l\in \mathcal{A}$. 
Let $\beta_i$ be the root perpendicular to $H_i$ with the opposite orientation 
to the path $\pi(t)$.
We call a sequence of roots $(\beta_1,\dots,\beta_l)$ 
obtained in such a way a {\it $\lambda$-chain.} 
In fact, $\lambda$-chains are in a bijective correspondence 
with decompositions of a certain
element $v_{-\lambda}$ of the affine Weyl group into products $v_{-\lambda} =
s_{i_1}\cdots s_{i_l}$ of the generators of $\Waff$.

For positive roots $\alpha\in \Phi^+$, let us define the
{\it Bruhat operators\/} $B_\alpha$ that act on the 
Grothendieck ring $K(G/B)$ by
$$
B_\alpha : [\O_{X_w}]\longmapsto\left\{\begin{array}{cl}
[\O_{X_{ws_\alpha}}] & \textrm{if }  \ell(ws_\alpha) = \ell(w)-1, \\[.05in]
0 & \textrm{otherwise.}
\end{array}\right.
$$
Also let $B_{-\alpha} = - B_\alpha$.
These operators are specializations of the quantum Bruhat 
operators from \cite{BFP}.  
The operators $1+B_\alpha$ satisfy the Yang-Baxter equation.

\begin{theorem} 
{\rm ($K$-theory Chevalley formula) }
Let $\lambda$ be any weight (dominant or nondominant, regular or nonregular).
Let $(\beta_1,\dots,\beta_l)$ be a $\lambda$-chain.  Then, for any $w\in W$, 
we have
$$
[\L_\lambda]\cdot [\O_{X_w}] = (1+B_{\beta_l})\cdots (1+B_{\beta_1}) ([\O_{X_w}])
$$
in the Grothendieck ring $K(G/B)$.
\label{th:intro}
\end{theorem}

The number of times a root $\alpha$ appears
in the $\lambda$-chain $(\beta_1,\dots,\beta_l)$ minus
the number of times $-\alpha$ appears in the $\lambda$-chain
equals $(\lambda,\alpha^\vee)$.
Thus the linear part of the expansion of $(1+B_{\beta_l})\cdots
(1+B_{\beta_1})$ is precisely $\sum_{\alpha>0} (\lambda,\alpha^\vee) 
\,B_\alpha$.
This linear part produces the classical Chevalley formula for 
products of classes in the cohomology ring $H^*(G/B)$.

We say that a $\lambda$-chain is {\it reduced\/} if it has minimal 
possible length.  Reduced $\lambda$-chains correspond to reduced decompositions
in the affine Weyl group.
If $\lambda$ is a dominant weight, then all roots 
in a reduced $\lambda$-chain are positive.  In this case,
Theorem~\ref{th:intro} involves only positive terms.  If $\lambda$ is an
antidominant weight, then all roots 
in a reduced $\lambda$-chain are negative.
In this case, the sign of the coefficient of 
$[\O_{X_w}]$ in $[\L_\lambda]\cdot
[\O_{X_u}]$ equals $(-1)^{\ell(u)-\ell(w)}$,
and Theorem~\ref{th:intro} gives a subtraction-free expression
for this coefficient.

Let $s_1,\dots,s_r$ be the system of simple reflections in
the Weyl group (compatible with our choice of Borel subgroup),
let $\omega_1,\dots,\omega_r$ be the corresponding set of fundamental
weights, and let $w_\circ$ be the longest element in $W$.
The {\it special\/} classes $[\O_{X_{w_\circ s_i}}]\in K(G/B)$ 
for codimension one Schubert varieties can be expressed as
$ [\O_{X_{w_\circ s_i}}]= 1 - [\L_{-\omega_i}]$. 
Note that $(\beta_1,\dots,\beta_l)$ is a $\lambda$-chain if and only if
$(-\beta_l,\dots,-\beta_1)$ is a $(-\lambda)$-chain.

\begin{corollary}  {\rm ($K$-theory Monk formula) }
Let us fix a simple reflection $s_i$. 
Let $(\beta_1,\dots,\beta_l)$ be a $(-\omega_i)$-chain.
Then, for any $w\in W$,  we have
$$
[\O_{X_{w_\circ s_i}}]\cdot [\O_{X_w}] = 
(1 - (1- B_{\beta_1})\cdots (1-B_{\beta_l}))([\O_{X_w}])
$$
in the Grothendieck ring $K(G/B)$.
\label{cor:K-hyperplane-sections}
\end{corollary}

The special classes $[\O_{X_{w_\circ s_i}}]$ generate the 
Grothendieck ring $K(G/B)$.
Thus Corollary~\ref{cor:K-hyperplane-sections} gives a complete
characterization of the multiplicative structure of the 
Grothendieck ring.

Our construction was developed independently of the LS-galleries of Gaussent 
and Littelmann~\cite{GaLi}.  Learning about the latter prompted us to subsequently
reformulate the model for characters of $V_\lambda$ that follows from our
formula by using admissible foldings of galleries. 
For regular  weights, our admissible foldings are similar
(but not equivalent!) to LS-galleries.  However, for nonregular weights,
these two models diverge.  Our model is simpler and more efficient 
computationally than the models based on LS-paths and LS-galleries.
It eliminates several choices that appear in the definitions of
LS-galleries and LS-paths.  Also it is harder to work with
sequences of lower dimensional faces of alcoves
(LS-galleries) than with reduced decompositions in the affine Weyl group (our
model). 
Note that
we cannot discard the case of nonregular weights as something of less
importance than regular weights.  The fundamental weights, which are highly
nonregular, are, in a sense, the most important weights for our purposes.
Indeed, these weights appear in Monk-type product formulas.
Also note that LS-galleries were not applied to Demazure characters and to 
the $K$-theoretic Chevalley formula; no generalization of them to Kac-Moody algebras is known either.

\medskip
The general outline of the paper is as follows.  In Section~\ref{sec:notation},
we review basic notions related to roots systems and fix our notation.  In
Section~\ref{sec:K-theory}, we present some background on the Grothendieck ring
$K_T(G/B)$.  In Section~\ref{sec:demazure}, we discuss the relationship between
the Grothendieck ring and the Demazure characters.  In
Section~\ref{sec:affine-weyl-group}, we remind a few facts about affine Weyl
groups.  In particular, we show that decompositions of affine Weyl group
elements correspond to sequences of adjacent alcoves, which we call alcove
paths.  In Section~\ref{sec:K-chevalley-formula}, we state our combinatorial
formula for products in equivariant $K$-theory, that is, our $K_T$-Chevalley
formula.  As a corollary of the $K_T$-Chevalley formula, we obtain a combinatorial model 
for the characters of the irreducible representations $V_\lambda$ 
and for the Demazure characters.  
In Section~\ref{sec:parabolic}, we extend the $K_T$-Chevalley formula
to equivariant $K$-theory of $G/P$.  In
Section~\ref{sec:duality}, we present several applications 
of our $K_T$-Chevalley formula.  We derive the $K_T$-Monk formula 
for the product of an arbitrary class $[\O_{X_w}]$ 
with a divisor class $[\O_{X_{w_\circ s_i}}]$, as well as the dual $K_T$-Chevalley formula.  Then we study two symmetries 
of the coefficients in the $K_T$-Chevalley formula.    
In the following sections, we develop tools needed to
reformulate our rule in a compact operator notation and to prove this rule.  In
Section~\ref{sec:yang-baxter}, we discuss the Yang-Baxter equation.  In
Section~\ref{sec:bruhat-operators}, we construct a certain $R$-matrix and show
that it satisfies the Yang-Baxter equation.  In
Section~\ref{sec:commutation-relations}, we derive commutation relations
between the elements of the $R$-matrix and the Demazure operators $T_i$.  These
commutation relations are the core of the proof of our formula.  In
Section~\ref{sec:path-operators}, we define compositions $R^{[\lambda]}$ of
elements of the $R$-matrix.  We use tail-flips of alcove paths to prove that
the operators $R^{[\lambda]}$ satisfy the same commutation relations with $T_i$
as the operators $E^\lambda$.   
In Section~\ref{sec:k-chevalley-operator-notation}, we reformulate and prove
our main result---the $K_T$-Chevalley formula---using the $R$-matrix notation.
We show that $R^{[\lambda]}$ coincides with the operator $E^\lambda$ 
of multiplication by the class of a line bundle in the Grothendieck ring $K_T(G/B)$. 
In Section~\ref{sec:central-points}, we use central points of alcoves to
prove the equivalence of the two formulations of our main result.
In Sections~\ref{sec:Type-A} and~\ref{sec:other-types}, we give several examples
for types $A$, $B$, $C$, and $G_2$.  In Section~\ref{sec:Quantum}, we
conjecture a natural generalization of our $K$-theory Monk formula to
quantum $K$-theory.
In Appendix~\ref{sec:galleries}, we reformulate our model for characters 
using admissible foldings of galleries and compare our 
model with LS-galleries and LS-paths.

\medskip
{\sc Acknowledgments:}
We are indebted to Shrawan Kumar for several geometric explanations and useful
suggestions.  We are grateful to V.~Lakshmibai for interesting discussions and
thoughtful comments.  We thank Allen Knutson, Yuan-Pin Lee, and Andrei
Zelevinsky for helpful remarks.

\section{Notation}
\label{sec:notation}

Let $G$ be a  connected, simply connected, simple complex Lie group.
Fix a Borel subgroup $B$ and a maximal torus $T$ such that 
$G\supset B\supset T$.
Let $\h$ be the corresponding Cartan subalgebra
of the Lie algebra $\mathfrak{g}$ of $G$. 
Let $r$ be the rank of the Cartan subalgebra $\h$.
Let $\Phi\subset \h^*$ be the 
corresponding irreducible {\it root system}. 
Let $\hR\subset \h^*$ be the real span of the roots.
Let $\Phi^+\subset \Phi$ be the set 
of positive roots corresponding to our choice of $B$. 
Then $\Phi$ is the disjoint union of $\Phi^+$ and $\Phi^- = -\Phi^+$.
Let $\alpha_1,\dots,\alpha_r\in\Phi^+$ be the corresponding 
{\it simple roots}.  They form a basis of $\hR$.
Let $(\lambda,\mu)$ denote the nondegenerate scalar product on $\hR$ induced by
the Killing form.  
Given a root $\alpha$, the corresponding {\it coroot\/} is $\alpha^\vee := 2\alpha/(\alpha,\alpha)$.  The collection of coroots
$\Phi^\vee:=\{\alpha^\vee \mid \alpha\in\Phi\}$ forms the 
{\it dual root system.}

The {\it Weyl group\/} $W\subset\mathrm{Aut}(\hR)$ 
of the Lie group $G$ is generated by the reflections 
$s_{\alpha}: \hR \to \hR$, for $\alpha\in\Phi$,
given by 
$$
s_\alpha:  \lambda \mapsto \lambda - (\lambda,\alpha^\vee)\,\alpha.
$$  
In fact, the Weyl group $W$ is generated by the
{\it simple reflections\/} $s_1,\dots,s_r$ corresponding 
to the simple roots $s_i := s_{\alpha_i}$, subject to the 
{\it Coxeter relations:}
$$(
s_i)^2=1
\quad\textrm{and}\quad 
(s_i s_j)^{m_{ij}}=1\quad\textrm{for any }i,j\in\{1,\dots,r\},
$$
where $m_{ij}$ is half of the order of the dihedral subgroup generated 
by $s_i$ and $s_j$.
An expression of a Weyl group element $w$ as a product 
of generators $w=s_{i_1}\cdots s_{i_l}$
which has minimal length is called a {\it reduced decomposition\/}
for $w$; its length $\ell(w)=l$ is called the {\it length\/} of $w$.
The Weyl group contains a unique {\it longest element\/} $w_\circ$
with maximal length $\ell(w_\circ)=|\Phi^+|$.
For $u,w\in W$, we say that $u$ {\it covers\/} $w$, and write $u\gtrdot w$,
if $w=u s_{\beta}$, for some $\beta\in\Phi^+$, and $\ell(u)=\ell(w)+1$.
The transitive closure ``$>$'' of the relation ``$\gtrdot$'' is called
the {\it Bruhat order\/} on $W$.

The {\it weight lattice\/} $\Lambda$ is given by
\begin{equation}
\Lambda:=\{\lambda\in \hR \mid (\lambda,\alpha^\vee)\in\Z
\textrm{ for any } \alpha\in\Phi\}.
\label{eq:weight-lattice}
\end{equation}
The weight lattice $\Lambda$ is generated by the 
{\it fundamental weights\/}
$\omega_1,\dots,\omega_r$, which are defined as the elements of the dual basis to the 
basis of simple coroots, i.e., $(\omega_i,\alpha_j^\vee)=\delta_{ij}$.
The set $\Lambda^+$ of {\it dominant weights\/} is given by
$$
\Lambda^+:=\{\lambda\in\Lambda \mid (\lambda,\alpha^\vee)\geq 0
\textrm{ for any } \alpha\in\Phi^+\}.
$$

Let $\rho:=\omega_1+\cdots+\omega_r=\frac{1}{2}\sum_{\beta\in\Phi^+}\beta$.
The {\it height\/} of a coroot $\alpha^\vee\in\Phi^\vee$ is 
$(\rho,\alpha^\vee) =  c_1+\cdots+c_r$ if 
$\alpha^\vee=c_1\alpha_1^\vee+\cdots + c_r \alpha_r^\vee$.
Since we assumed that $\Phi$ is irreducible, there is 
a unique {\it highest coroot\/} $\theta^\vee\in\Phi^\vee$ that has 
maximal height.  (In other words, $\theta^\vee$ is the highest root
of the dual root system $\Phi^\vee$.  It should not be confused with 
the coroot of the highest root of $\Phi$.)
We will also use the {\it Coxeter number}, that 
can be defined as $h:=(\rho,\theta^\vee)+1$.

\section{Equivariant $K$-theory of Generalized Flag Varieties}
\label{sec:K-theory}

In this section, we remind a few facts about the Grothendieck ring 
$K_T(G/B)$.
For more details on the Grothendieck ring, we refer to Kostant 
and Kumar~\cite{KK},  see also Pittie and Ram~\cite{PR}.  
\medskip

The {\it generalized flag variety\/} $G/B$ is a smooth projective
variety.  It decomposes into a disjoint union of {\it Schubert cells\/} 
$X_w^\circ := BwB/B$ indexed by elements $w\in W$ of the Weyl group.
The closures of Schubert cells $X_w := \overline{X_w^\circ}$
are called {\it Schubert varieties.} 
We have $u>w$ in the Bruhat order (defined as above) 
if and only if $X_u\supset X_w$. 
Let $\O_{X_w}$ be the structure sheaf of the Schubert variety $X_w$.

Let $\Z[\Lambda]$ be the group algebra of the weight lattice $\Lambda$.  It has
a $\Z$-basis of formal exponents $\{e^\lambda \mid \lambda\in\Lambda\}$ with
multiplication $e^\lambda\cdot e^\mu := e^{\lambda+\mu}$, i.e.,
$\Z[\Lambda]=\Z[e^{\pm\omega_1},\cdots,e^{\pm\omega_r}]$ is the algebra of
Laurent polynomials in  $r$ variables.  The group of characters  $X=X(T)$ of
the maximal torus $T$ is isomorphic to the weight lattice $\Lambda$.  Its group
algebra $\Z[X]=R(T)$ is the representation ring of $T$.  The rings
$\Z[\Lambda]$ and $\Z[X]$ are isomorphic.  (However we will distinguish these
two rings.) Let us denote by $x^\lambda$ the element of $\Z[X]$ corresponding
to the character determined by $\lambda$, as well as to $e^\lambda\in\Z[\Lambda]$.  Thus
$\Z[X]=\Z[x^{\pm\omega_1},\cdots,x^{\pm\omega_r}]$.  Let $\mathcal{L}_\lambda$ be the line bundle over $G/B$ associated with the
weight $\lambda$, that is, $\mathcal{L}_\lambda:= G\times_B \C_{-\lambda}$, where $B$ acts on $G$ by right multiplication, and the
$B$-action on $\C_{-\lambda}=\C$ corresponds to the character determined by $-\lambda$. 
(This character of $T$ extends to $B$ by defining it to be 
identically one on the commutator subgroup $[B,B]$.) 

Denote by $K_T(G/B)$ the {\it Grothendieck ring\/} of coherent $T$-equivariant
sheaves on $G/B$.  According to Kostant and Kumar~\cite{KK}, the Grothendieck
ring $K_T(G/B)$ is a free $\Z[X]$-module, and the classes $[\O_{X_w}]\in K_T(G/B)$
of the structure sheaves of Schubert varieties form its $\Z[X]$-basis. 
The classes $[\mathcal{L}_\lambda]$ of the line bundles
$\L_\lambda$ also span $K_T(G/B)$ as a $\Z[X]$-module. 

We now discuss the presentation of the Grothendieck ring $K_T(G/B)$ 
as a quotient of $\Z[X]\otimes\Z[\Lambda]$.
The Weyl group $W$ acts on the group algebra $\Z[\Lambda]$
by $w(e^\lambda) := e^{w(\lambda)}$.
Let $\Z[\Lambda]^W$ be the subalgebra of $W$-invariant elements.
The tensor product $\Z[X]\otimes\Z[\Lambda]$
is the algebra of Laurent polynomials in $2r$ variables 
$x^{\omega_1},\dots,x^{\omega_r}, 
e^{\omega_1},\dots,e^{\omega_r}$ with integer coefficients. 
Let $i:\Z[\Lambda]\to \Z[X]$ be the natural isomorphism given by
$i(e^\lambda):=x^\lambda$.
Let $\mathcal{I}$ be the ideal in $\Z[X]\otimes\Z[\Lambda]$
generated by the following elements:
$$
\mathcal{I}:=
\left<i(f)\otimes 1 - 1\otimes f \mid 
f\in \Z[\Lambda]^W \right>.
$$
The Grothendieck ring $K_T(G/B)$ is canonically isomorphic 
to the quotient ring
\begin{equation}
K_T(G/B)\simeq (\Z[X]\otimes\Z[\Lambda])/ \mathcal{I}.
\label{eq:KGB=ZL/I}
\end{equation}
The isomorphism is given by the $\Z[X]$-linear map 
$[\L_\lambda]\mapsto e^{-\lambda}$, for 
$\lambda\in\Lambda$. From now on, we will identify the two rings. Recall from \cite{KK} the ${\mathbb Z}$-linear involution $*\,:\,K_T(G/B)\rightarrow K_T(G/B)$ given by $x^\mu \otimes e^\lambda\mapsto x^{-\mu} \otimes e^{-\lambda}$; in other words, this map takes a vector bundle to its dual. Let $[\O_w]:=*[\O_{X_w}]$. Throughout most of this paper, we will work with the classes $[\O_w]$ instead of $[\O_{X_w}]$; it is straightforward to rephrase all results in terms of $[\O_{X_w}]$.

It is possible to express all classes $[\O_w]$ 
as Laurent polynomials in $\Z[X]\otimes\Z[\Lambda]$ by choosing 
a representative of the class $[\O_1]$ 
and by applying Demazure operators, as described below. 
The action of the Weyl group on $\Z[\Lambda]$ defined above is extended $\Z[X]$-linearly to 
$\Z[X]\otimes\Z[\Lambda]$. 
For $i=1,\dots,r$, the elementary {\it Demazure operator\/} 
$T_i:\Z[X]\otimes\Z[\Lambda]\to \Z[X]\otimes\Z[\Lambda]$ is 
the $\Z[X]$-linear operator given by
\begin{equation}
T_i(f) :=  \frac{f -  e^{-\alpha_i} s_i(f)}{1- e^{-\alpha_i}}\,.
\label{eq:isobaric-div-diff}
\end{equation}
Note that the numerator is always divisible 
by the denominator\footnote{Check this for $f=e^\lambda$.},
so the right-hand side is a valid expression in the algebra
$\Z[X]\otimes\Z[\Lambda]$.  One can verify directly from the 
definition that the operators $T_i$ satisfy the following relations:
\begin{eqnarray}
&&T_i^2 = T_i \,,
\label{eq:T-reduce}
\\[.05in]
&&(T_i\, T_j)^{m_{ij}} = 1\,,
\label{eq:T-coxeter}
\\[.05in]
&&T_i(fg) = f\cdot T_i(g),\quad\textrm{if } s_i(f) = f\,.
\label{eq:Ti(fg)}
\end{eqnarray}
Equations~(\ref{eq:T-reduce}) and~(\ref{eq:T-coxeter}) imply that
the operators $T_i$ give an action of the corresponding Hecke algebra $\mathcal{H}_q$
specialized at $q=0$, e.g., see~\cite{Hum}.
Equation~(\ref{eq:Ti(fg)}) implies that the operators $T_i$ 
preserve the ideal $\mathcal{I}$.
Thus the elementary Demazure operators $T_i$ induce operators acting 
on the Grothendieck ring 
$K_T(G/B)\simeq (\Z[X]\otimes\Z[\Lambda])/\mathcal{I}$, 
which will be denoted by the same symbols.

For a reduced decomposition $w=s_{i_1}\cdots s_{i_l}\in W$, 
the {\it Demazure operator\/} $T_w$ is defined as the following 
composition of elementary Demazure operators:
\begin{equation}
T_w := T_{i_1}\cdots T_{i_l}.
\label{eq:demazure-w}
\end{equation}
The Coxeter relations~(\ref{eq:T-coxeter})
imply that the operator $T_w$ depends only on $w$,  not on
the choice of a reduced decomposition.
Equation~(\ref{eq:T-reduce}) implies that an arbitrary 
product $T_{j_1}\cdots T_{j_m}$ reduces to $T_w$ for some $w\in W$.
Kostant and Kumar~\cite{KK} showed that, for any $w\in W$,
\begin{equation}
[\O_w] = T_{w^{-1}}([\O_1]).
\label{eq:KK}
\end{equation}

For type $A$, the elementary Demazure operators $T_i$ are also called
{\it isobaric divided difference operators.}  The polynomial
representatives of the classes $[\O_w]$ obtained by applying 
these operators to a certain polynomial representative of $[\O_1]$ are the {\it double Grothendieck polynomials} of Lascoux and Sch\"utzenberger \cite{LS2}.

The product $e^\lambda\cdot [\O_u]$ in the Grothendieck ring 
$K_T(G/B)$ can be written as a finite sum
\begin{equation}
e^\lambda\cdot [\O_u] = \sum_{w\in W,\,\mu\in\Lambda} 
c_{u,w}^{\lambda,\mu} \ x^\mu\,[\O_w]\,,
\label{eq:LlOu}
\end{equation}
where $c_{u,w}^{\lambda,\mu}$ are some integer coefficients. Equivalently, we can write
\begin{equation}\label{eq:kch1}
[\L_\lambda]\cdot [\O_{X_u}] = \sum_{w\in W,\,\mu\in\Lambda} 
c_{u,w}^{\lambda,\mu} \ x^{-\mu}\,[\O_{X_w}]\,.
\end{equation}
We will call the coefficients $c_{u,w}^{\lambda,\mu}$ 
{\it $K_T$-Chevalley coefficients,} because they extend
the coefficients in the usual Chevalley formula, as shown
below in this section. 
In this paper, we present an explicit combinatorial formula for  
$c_{u,w}^{\lambda,\mu}$, see 
Theorems~\ref{th:main-combinatorial-formulation} and~\ref{thm:K-Chevalley}.
We will see that $c_{u,w}^{\lambda,\mu}= 0$ unless $w\leq u$ in the 
Bruhat order, and that $c_{u,u}^{\lambda,\mu} = \delta_{\lambda,\mu}$.
If $\lambda$ is a dominant weight, then we will see that
all coefficients $c_{u,w}^{\lambda,\mu}$ are nonnegative.
In this case, Pittie and Ram~\cite{PR} showed that  
$c_{u,w}^{\lambda,\mu}$ count certain LS-paths,
cf.\ also Lakshmibai-Littelmann~\cite{LL} and
Littelmann-Seshadri~\cite{LS3}.

For a weight $\lambda$,	
let $E^\lambda:f\mapsto e^\lambda f$ be the operator of multiplication
by the exponent $e^\lambda$ in the ring $\Z[X]\otimes\Z[\Lambda]$.  
The induced operator on $K_T(G/B)$, which will be denoted by the same symbol $E^\lambda$, acts as the operator 
of multiplication by the class $[\L_{-\lambda}]$ of a line bundle.
It follows from the definitions that $E^\lambda$ and $T_i$ satisfy
the following commutation relation:
\begin{equation}
E^\lambda \,T_i = T_i \,E^{s_i(\lambda)} +
\frac{E^\lambda - E^{s_i(\lambda)}}{1 - E^{-\alpha_i}}. 
\label{eq:YT=TY}
\end{equation}
The quotient in this expression expands as the Laurent polynomial
$$
\frac {E^\lambda-E^{s_i(\lambda)}}{1-E^{-\alpha_i}}
= \sum_{0\leq k<(\lambda,\alpha_i^\vee)} E^{\lambda- k\alpha_i}
- \sum_{(\lambda,\alpha_i^\vee)\leq k< 0} E^{\lambda- k\alpha_i}.
$$
Also, we have
\begin{equation}
E^\lambda ([\O_1])= x^\lambda\,[\O_1].
\label{eq:YO1}
\end{equation}

Let $\hH$ be the ring generated by the operators
$T_1,\dots,T_r$ and $E^\lambda$, $\lambda\in\Lambda$. 
Then $\hH$ is described by relations~(\ref{eq:T-reduce}), 
(\ref{eq:T-coxeter}), and~(\ref{eq:YT=TY}), i.e.,
$\hH$ is a certain degeneration of the affine Hecke algebra.
This follows from the fact that the elements $T_{w^{-1}} E^\mu$,
$w\in W$, $\mu\in\Lambda$, form a $\Z$-basis of $\hH$.  
Indeed, according to the relations, 
the elements $T_{w^{-1}} E^\mu$ span $\hH$.  On the other hand, 
these elements are  linearly independent, because
$T_{w^{-1}} E^\mu([\O_1])= x^\mu [\O_{w}]$.

Using the commutation relation in~(\ref{eq:YT=TY})
repeatedly, we obtain, for any $u\in W$ and $\lambda\in \Lambda$,
the following identity in the ring $\mathcal{\hat H}$: 
\begin{equation}
E^\lambda \,T_{u^{-1}} = \sum_{w\in W,\,\mu\in\Lambda}
c_{u,w}^{\lambda,\mu}\ T_{w^{-1}} \,E^{\mu},
\label{eq:YT=TY-long}
\end{equation}
for some integer coefficients $c_{u,w}^{\lambda,\mu}$.
Applying both sides of this expression to the class $[\O_1]$
and using~(\ref{eq:KK}) and~(\ref{eq:YO1}), we deduce that
the coefficients $c_{u,w}^{\lambda,\mu}$ in~(\ref{eq:YT=TY-long}) are 
equal to the $K_T$-Chevalley coefficients in~(\ref{eq:LlOu}).

The commutation relation~(\ref{eq:YT=TY}) gives a recursive procedure for
calculating the product $e^\lambda\cdot[\O_u]$ in $K_T(G/B)$.
In this paper, we present a simple nonrecursive rule for this product.  
The proof of our rule is based on the following trivial observation,
which is implied by the above discussion.

\begin{lemma}
Let $A$ be an algebra that contains $\Z[X]$, and let 
$\tilde K = K_T(G/B)\otimes_{\Z[X]} A$.
The action of the Demazure operators $T_i$ extends $A$-linearly to 
$\tilde K$. 
Suppose that $R^\lambda$, $\lambda\in\Lambda$, 
is a family of $A$-linear operators acting on 
the space $\tilde K$ such that relations~{\rm(\ref{eq:YT=TY})}
and~{\rm(\ref{eq:YO1})} hold with $E^\lambda$ replaced by $R^\lambda$.
Then the operator $R^\lambda$ preserves $K_T(G/B)\subset \tilde K$ and
coincides with $E^\lambda$ for all $\lambda$.
\label{lem:simple-observation}
\end{lemma}

\begin{proof} 
The conditions imply that relation~(\ref{eq:YT=TY-long}) holds
with $E^\lambda$ replaced by $R^\lambda$.  Applying this expression 
to $[\O_1]$, we deduce that 
$R^\lambda([\O_u]) = E^\lambda ([\O_u])$, for any $u\in W$.
\end{proof}

Let us also mention another basis of $K_T(G/B)$ studied by
Kostant and Kumar~\cite{KK}, see also the recent paper~\cite{GR}
by Griffeth and Ram.
One can easily check that there is an algebra involution $\psi$ of the ring $\mathcal{\hat H}$ given by $\psi:T_i\mapsto 1-T_i$, 
$i=1,\dots,r$, and $\psi:E^\lambda\mapsto E^{-\lambda}$.
In other words, the operators $\varepsilon_i = 1-T_i$, for $i=1,\dots,r$,  
satisfy relations~(\ref{eq:T-reduce}), (\ref{eq:T-coxeter}), 
and~(\ref{eq:YT=TY}) with $T_i$ replaced by $\varepsilon_i$ and $E^\lambda$
replaced by $E^{-\lambda}$.  Thus one can correctly define the elements 
$\varepsilon_w := 
\varepsilon_{i_1}\cdots \varepsilon_{i_l}\in\mathcal{\hat H}$, 
for a reduced decomposition $w=s_{i_1}\cdots s_{i_l}\in W$.
For $w\in W$, let $[\I_w]$ be the element of $K_T(G/B)$ given by
\begin{equation}
[\I_w] = \varepsilon_{w^{-1}}([\O_1]).
\label{eq:I-dual-basis}
\end{equation}
It turns out that the elements $[\I_w]$, $w\in W$,
form a $\Z[X]$-basis of $K_T(G/B)$, as well.  Moreover, the bases 
$\{[\I_w] \mid w\in W\}$ and $\{[\O_w] \mid w\in W\}$
are related to each other as follows:
$$
[\I_w] = \sum_{u\leq w} (-1)^{\ell(u)} [\O_u]\qquad\textrm{and}\qquad
[\O_w] = \sum_{u\leq w} (-1)^{\ell(u)} [\I_u].
$$
These two relations are easy to check by induction on the length of $w$.  

The element $[\I_w]$ can be described geometrically. Up to sign, it is the class $*[\I_{X_w}]$, where $\I_{X_w}$ is the sheaf given by
the exact sequence $0\to \I_{X_w}\to\O_{X_w}\to\O_{\partial X_w}\to 0$, and 
$\partial X_w = \bigcup_{u<w} X_u$ is the boundary of the Schubert variety
$X_w$ (cf. \cite[Theorem 2.1 (ii)]{Mat}, \cite[Equation (4)]{LS3}, and \cite[Section 2]{GR}). Brion and Lakshmibai~\cite{BL} showed that the classes $[\I_{X_w}]$
form the dual basis to $\{[\O_{X_w}]\mid w\in W\}$ with respect to the natural 
intersection pairing in $K$-theory.

Applying the above involution $\psi$ to both sides of~(\ref{eq:YT=TY-long}),
we obtain
$$
E^{-\lambda} \,\varepsilon_{u^{-1}} = \sum_{w\in W,\,\mu\in\Lambda}
c_{u,w}^{\lambda,\mu}\ \varepsilon_{w^{-1}} \,E^{-\mu}.
$$
Then applying  both sides of this relation to $[\O_1]$,
we immediately deduce the following {\it dual form\/} 
of~{\rm(\ref{eq:LlOu})}
\begin{equation}
e^{-\lambda}\cdot [\I_u] = \sum_{w\in W,\,\mu\in\Lambda}
c_{u,w}^{\lambda,\mu} \,x^{-\mu}\, [\I_w],
\label{eq:L-I}
\end{equation}
where $c_{u,w}^{\lambda,\mu}$ are the same $K_T$-Chevalley coefficients
as those in~{\rm(\ref{eq:LlOu})} and~(\ref{eq:YT=TY-long}).

Note that relations~(\ref{eq:T-reduce}), (\ref{eq:T-coxeter}), 
and~(\ref{eq:YT=TY})
in the algebra $\mathcal{\hat H}$ are equivalent to the relations
obtained from them by reversing the order of all terms.
This symmetry of the relations implies that the expression
\begin{equation}
T_{u}\, E^\lambda = \sum_{w\in W,\,\mu\in\Lambda}
c_{u,w}^{\lambda,\mu}\ E^{\mu}\, T_{w} 
\label{eq:DE=ED-dual}
\end{equation}
has the same $K_T$-Chevalley coefficients $c_{u,w}^{\lambda,\mu}$.

The (nonequivariant) Grothendieck ring $K(G/B)$ of coherent sheaves on $G/B$ 
can be obtained by the specialization $x^\mu \mapsto 1$, for all $\mu$, i.e., 
by ignoring all exponents $x^\mu$ in equivariant $K$-theory. By a slight abuse of notation, we will use the same symbols $[\O_{X_w}]$  
and $[\L_\lambda]$ for the obvious classes in $K(G/B)$
as in equivariant $K$-theory.
 The classes $[\O_{X_w}]$, $w\in W$, form a $\Z$-basis of $K(G/B)$.

Let us also recall the way in which Schubert calculus in cohomology 
can be recovered from $K$-theory.
Let $H^*(G/B):=H^*(G/B,\Q)$ be the cohomology ring of $G/B$ with 
rational coefficients.  It has a linear basis of classes of Schubert varieties
$[X_w]$, $w\in W$, called {\it Schubert classes.} 
The cohomology ring is $2\Z$-graded by $\deg([X_w]) = 2(\ell(w_\circ)-\ell(w))$.
Let $\h^*_\Q\subset \h^*$ be the $\Q$-span of the weight lattice $\Lambda$,
and let $\mathit{Sym}(\h^*_\Q)$ be its symmetric algebra,
i.e., the ring of  polynomials on $\h_\Q$.
The classical {\it Borel theorem\/} says that the cohomology ring $H^*(G/B)$ 
is isomorphic to the following quotient of the symmetric algebra:
$$
H^*(G/B)\simeq \mathit{Sym}(\h^*_\Q)/\mathcal{J},
$$
where $\mathcal{J}:=\left<f\in \mathit{Sym}(\h^*_\Q)^W \mid f(0)=0\right>$
is the ideal generated by $W$-invariant polynomials without constant term.
The isomorphism identifies the Chern class $[\lambda]\in H^2(G/B)$ 
of the line bundle $\L_\lambda$ with the coset of $\lambda$ 
modulo $\mathcal{J}$.
The product of $[\lambda]$ and a Schubert class $[X_u]$ 
in the cohomology ring is given by the following classical formula due to Chevalley~\cite{Chev}:
\begin{equation}
[\lambda]\cdot [X_u]  = \sum_{\alpha\in\Phi^+,\,
\ell(u s_\alpha) = \ell(u) - 1} (\lambda,\alpha^\vee)\,[X_{us_\alpha}].
\label{eq:chevalley}
\end{equation}
The {\it Chern character\/} is the ring isomorphism
$\mathit{ChCh}:K(G/B)\otimes\Q\to H^*(G/B)$  that sends the class 
$[\L_\lambda]\in K(G/B)$ 
of the line bundle $\L_\lambda$ to 
$\exp{[\lambda]} := 1 + [\lambda] + [\lambda]^2/2! + \cdots\in H^*(G/B)$.
Then
$$
\mathit{ChCh}([\O_w]) = [X_w] + \textrm{higher degree terms}.
$$

This shows that the Chevalley formula~(\ref{eq:chevalley}) for the product 
$[\lambda]\cdot [X_u]$ in $H^*(G/B)$ is obtained  
from the expression $[\L_\lambda]\cdot [\O_{X_u}]-[\O_{X_u}]$ in $K_T(G/B)$
by expanding it using~(\ref{eq:kch1}),
ignoring the exponents $x^{-\mu}$, applying the Chern character map,
and then extracting terms of degree $\deg([X_u]) + 2$. 
In other words, for $\lambda\in\Lambda$, $u\in W$, $\alpha\in\Phi^+$
such that $\ell(us_\alpha)= \ell(u)-1$, the coefficient
in the Chevalley formula equals
\begin{equation}
 (\lambda,\alpha^\vee)=
\sum_{\mu\in\Lambda}c_{u,us_\alpha}^{\lambda,\mu}.
\label{eq:K-theory->Chevalley}
\end{equation}
A rule for computing the coefficients $c_{u,w}^{\lambda,\mu}$ can be 
thought of as a generalization of the Chevalley formula
to $T$-equivariant $K$-theory.

\begin{remark} In fact, Pittie and Ram~\cite{PR} worked in 
a more general setup than the 
Grothendieck ring $K_T(G/B)$.
Their construction implies that the same $K_T$-Chevalley coefficients 
$c_{u,w}^{\lambda,\mu}$ as in~(\ref{eq:kch1})
give the product of the classes of $\L_\lambda$ and 
$\O_{X_u}$ in the $K$-theory of a $G/B$-bundle over a smooth base.  
Thus, the results of the present paper apply to this more general case as well.
\end{remark}

\section{Demazure Characters}
\label{sec:demazure}

Lakshmibai-Littelmann~\cite{LL} and Littelmann-Seshadri~\cite{LS3}
indicated that
the product $[\L_\lambda]\cdot [\O_{X_u}]$ in the Grothendieck ring 
$K_T(G/B)$ is related to representation theory.  
This relation is also implicit in the Pittie-Ram formula~\cite{PR}.
Kumar~\cite{Kum} pointed out that the Demazure characters can be expressed
in terms of the $K_T$-Chevalley coefficients, as shown below.
\medskip

For a dominant weight $\lambda\in \Lambda^+$,
let $V_\lambda$ denote the finite dimensional irreducible representation of 
the Lie group $G$ with highest weight $\lambda$. 
For $\lambda\in\Lambda^+$ and $w\in W$, 
the {\it Demazure module\/} $V_{\lambda,w}$ is the $B$-module that 
is dual to the space of global 
sections of the line bundle $\L_\lambda$ on the Schubert variety $X_w$:
\begin{equation}\label{demaz}
V_{\lambda,w} := H^0(X_w,\L_\lambda)^*.  
\end{equation}
For the longest Weyl group element $w=w_\circ$,
the space $V_{\lambda,w_\circ} = H^0(G/B,\L_\lambda)^*$ has the structure of
a $G$-module.  The classical {\it Borel-Weil theorem\/} says that 
$V_{\lambda,w_\circ}$ is isomorphic to the irreducible $G$-module $V_\lambda$.
The formal characters of these modules, called {\it Demazure characters}, 
are given by $ch(V_{\lambda,w})=\sum_{\mu\in\Lambda} 
m_{\lambda,w}(\mu)\,e^\mu\in\Z[\Lambda]$,
where $m_{\lambda,w}(\mu)$ is the multiplicity of the weight $\mu$ in 
$V_{\lambda,w}$.
They generalize the characters of the irreducible representations
$ch(V_\lambda)=ch(V_{\lambda,w_\circ})$.
The {\it Demazure character formula}~\cite{Dem} says that 
the character $ch(V_{\lambda,w})$ is given by
\begin{equation}
ch(V_{\lambda,w}) = T_w(e^\lambda),
\label{eq:demazure-character-formula}
\end{equation}
where $T_w$ is the Demazure operator~(\ref{eq:demazure-w}).

\begin{lemma}
For any $\lambda\in\Lambda^+$ and $u\in W$, the Demazure character
$ch(V_{\lambda,u})$ can be expressed
in terms of the $K_T$-Chevalley coefficients $c_{u,w}^{\lambda,\mu}$ 
in~{\rm(\ref{eq:LlOu})} as follows:
$$
ch(V_{\lambda,u}) = \sum_{w\in W,\, \mu\in\Lambda} c_{u,w}^{\lambda,\mu}
\, e^\mu.
$$
In particular, the character of the irreducible representation $V_\lambda$
of $G$ is equal to
$$
ch(V_{\lambda}) = \sum_{w\in W,\, \mu\in\Lambda} c_{w_\circ,w}^{\lambda,\mu}
\, e^\mu.
$$
\label{lem:demazure-character}
\end{lemma}

\begin{proof}
Applying both sides of identity~(\ref{eq:DE=ED-dual})
to $[\O_{w_\circ}]=1$ and using $T_w(1) = 1$, we obtain
$$
T_u(e^\lambda) = \sum_{w\in W,\,\mu\in\Lambda}
c_{u,w}^{\lambda,\mu}\ e^{\mu},
$$
which, together with the Demazure character 
formula~(\ref{eq:demazure-character-formula}), proves the lemma.
\end{proof}

Let us also give a geometric argument that proves Lemma~\ref{lem:demazure-character}.  It is implicit in~\cite{LL} and~\cite{LS3}
and was reported to us by Kumar~\cite{Kum}.
Let $\chi:K_T(G/B)\to \Z[\Lambda]$ be the {\it Euler characteristic map\/}
given by
$$
\chi:[\mathcal{V}] \longmapsto
\sum_{i\geq 0} (-1)^i\, ch(H^i(G/B,\mathcal{V})^*),
$$
for a coherent sheaf $\mathcal{V}$ on $G/B$.
For a dominant weight $\lambda$, the Euler characteristic
$\chi([\L_\lambda]\cdot [\O_{X_u}])$ is equal to the Demazure character
$ch(V_{\lambda,u})$.  Indeed, this follows from (\ref{demaz}), the fact that
$$
H^i(G/B,\L_\lambda\otimes \O_{X_u})=H^i(X_u,\L_\lambda)\,,
$$
and the vanishing of the cohomologies $H^i(X_u,\L_\lambda)$, for $i\geq 1$. 
In particular, we have $\chi([\O_{X_w}])=1$, for any $w\in W$.
Thus $\chi(x^{-\mu}[\O_{X_w}]) = e^\mu$.
Applying the Euler characteristic map $\chi$ to both sides 
of~{\rm(\ref{eq:kch1})}, we obtain Lemma~\ref{lem:demazure-character}.

\section{Affine Weyl Groups}
\label{sec:affine-weyl-group}

In this section, we remind a few basic facts about 
affine Weyl groups and alcoves, cf.  
Humphreys~\cite[Chaper~4]{Hum} for more details.
Then we define $\lambda$-chains that will be used in the rest of the paper.
\medskip

Let $\Waff$ be  the {\it affine Weyl group\/} for the 
Langland's dual group $G^\vee$.
The affine Weyl group $\Waff$ is generated by the affine reflections 
$s_{\alpha,k}: \hR \to \hR$, for $\alpha\in\Phi$ and $k\in\Z$, 
that reflect the space $\hR$ with respect to the affine hyperplanes
\begin{equation}
H_{\alpha,k} := \{\lambda\in \hR \mid (\lambda,\alpha^\vee)=k\}.
\label{eq:H-alpha-k}
\end{equation}
Explicitly, the affine reflection $s_{\alpha,k}$ is given by 
$$
s_{\alpha,k}: \lambda \mapsto 
s_{\alpha}(\lambda) + k\,\alpha =
\lambda - ((\lambda,\alpha^\vee)-k)\,\alpha.
$$
The hyperplanes $H_{\alpha,k}$ divide the real vector space $\hR$ into open
regions, called {\it alcoves.} 
Each alcove $A$ is given by inequalities of the form
$$
A:=\{\lambda\in \hR \mid m_\alpha<(\lambda,\alpha^\vee)<m_\alpha+1
\textrm{ for all } \alpha\in\Phi^+\},
$$
where $m_\alpha=m_\alpha(A)$, $\alpha\in\Phi^+$, are some integers.

A proof of the following important property of the affine Weyl group 
can be found, e.g., in~\cite[Chapter~4]{Hum}.

\begin{lemma}
The affine Weyl group $\Waff$ acts simply transitively 
on the collection of all alcoves.
\label{lem:simply-transitively}
\end{lemma}

The {\it fundamental alcove\/} $A_\circ$ is given by 
$$
A_\circ :=\{\lambda\in \hR \mid 0<(\lambda,\alpha^\vee)<1 \textrm{ for all }
\alpha\in\Phi^+\}.
$$
Lemma~\ref{lem:simply-transitively} implies that, for any alcove $A$, 
there exists a unique element $v_A$ of the affine Weyl group $\Waff$
such that $v_A(A_\circ) = A$.  Hence the map $A\mapsto v_A$ is a one-to-one
correspondence between alcoves and elements of the affine Weyl group.

Recall that $\theta^\vee\in\Phi^\vee$ is the highest coroot. 
Let $\theta\in\Phi^+$ be the corresponding root,
and let $\alpha_0:=-\theta$.
The fundamental alcove $A_\circ$ is, in fact, the simplex given by
\begin{equation}
A_\circ =\{\lambda\in \hR \mid 0<(\lambda,\alpha_i^\vee) \textrm{ for }
i=1,\dots,r, \textrm{ and }(\lambda,\theta^\vee)<1\},
\label{eq:fund-alcove}
\end{equation}
Lemma~\ref{lem:simply-transitively} also implies that the affine Weyl group 
is generated by the set of reflections $s_0,s_1,\dots,s_r$ with respect 
to the walls of the fundamental alcove $A_\circ$, where 
$s_0 := s_{\alpha_0,-1}$ and $s_1,\dots,s_r\in W$ 
are the simple reflections $s_i=s_{\alpha_i,0}$.
As before, a decomposition $v=s_{i_1}\cdots s_{i_l}\in \Waff$ 
is called {\it reduced\/} if it has minimal length;
its length $\ell(v)=l$ is called the length of $v$.

Like the Weyl group, the affine Weyl group $\Waff$ is a 
Coxeter group, i.e., it is described by the relations 
\begin{equation}
(s_i)^2=1
\quad\textrm{and}\quad 
(s_i s_j)^{m_{ij}}=1,\quad\textrm{for any }i,j\in\{0,\dots,r\},
\label{eq:coxeter-relations-affine}
\end{equation}
where $m_{ij}$ is half of the order of the dihedral subgroup generated 
by $s_i$ and $s_j$.

We say that two alcoves $A$ and $B$ are {\it adjacent} 
if $B$ is obtained by an affine reflection of $A$ with respect to one of its 
walls.  In other words, two alcoves are adjacent if they are
distinct and have a common wall.  
For a pair of adjacent alcoves, let us write 
$A\stackrel{\beta}\longrightarrow B$ if the common wall of $A$ and $B$ 
is of the form $H_{\beta,k}$ and the root $\beta\in\Phi$ points 
in the direction from $A$ to $B$.  
By the definition, all alcoves that are adjacent to the fundamental 
alcove $A_\circ$ are obtained from $A_\circ$ by the reflections 
$s_0,\cdots,s_r$, and $A_\circ\stackrel{-\alpha_i}\longrightarrow s_i(A_\circ)$.

\begin{definition}
An {\it alcove path\/} is a sequence of alcoves
$(A_0,A_1,\dots,A_l)$ such that $A_{j-1}$ and $A_j$ are adjacent, for
$j=1,\dots,l$.
Let us say that an alcove path is {\it reduced\/} if it has minimal 
length among all alcove paths from $A_0$ to $A_l$.
\end{definition}

Let $v\mapsto \bar v$ be the homomorphism $\Waff\to W$ defined by ignoring the affine translation.
In other words,   $\bar s_{\alpha,k} = s_\alpha\in W$.

The following lemma, which is essentially well-known,
summarizes some properties of decompositions in affine Weyl groups,
cf.~\cite{Hum}.

\begin{lemma}  
Let $v$ be any element of $\Waff$,
and let $A=v(A_\circ)$ be the corresponding alcove.
Then the decompositions $v=s_{i_1}\cdots s_{i_l}$ of $v$ (reduced or not) as a product of generators in $\Waff$ are in one-to-one correspondence
with alcove paths 
$A_0\stackrel{-\beta_1}\longrightarrow A_1
\stackrel{-\beta_2}\longrightarrow 
\cdots \stackrel{-\beta_l}\longrightarrow A_l$
from the fundamental alcove $A_0=A_\circ$ 
to $A_l=A$.  This correspondence is explicitly given by
$A_j = s_{i_1}\cdots s_{i_j} (A_\circ)$, for $j=0,\dots,l$;
and the roots $\beta_1,\dots,\beta_l$ are given by
$$
\beta_1 = \alpha_{i_1}, \  \beta_2=\bar s_{i_1}(\alpha_{i_2}), \ 
\beta_3= \bar s_{i_1} \bar s_{i_2} (\alpha_{i_3}), \dots, \ 
\beta_l= \bar s_{i_1} \cdots \bar s_{i_{l-1}} (\alpha_{i_l}).
$$
Let $r_j\in\Waff$ denote the affine reflection with respect to the 
common wall of the alcoves $A_{j-1}$ and $A_j$, for $j=1,\dots,l$.
Then the affine reflections $r_1,\dots,r_l$ are given by
$$
r_1 = s_{i_1}, \ r_2 = s_{i_1} s_{i_2} s_{i_1}, \ 
r_3 =  s_{i_1} s_{i_2} s_{i_3} s_{i_2} s_{i_1}, \ \dots, \
r_l =  s_{i_1} \cdots s_{i_r} \cdots s_{i_1}. 
$$
We have $\bar{r}_i=s_{\beta_i}$ and $v=s_{i_1}\cdots s_{i_l} = r_l \cdots r_1$.
Moreover, the following claims are equivalent:
\begin{enumerate}
\item[(a)] $v=s_{i_1}\cdots s_{i_l}$ is a reduced decomposition;
\item[(b)] 
$(A_0,A_1,\dots,A_l)$ is a reduced alcove path;
\item[(c)] all affine reflections $r_1,\dots,r_l$ are distinct;
\item[(d)] $\beta_i\ne-\beta_j$, for any $i$ and $j$.
\end{enumerate}
Finally, for any $\alpha\in\Phi^+$, we have 
$m_\alpha(A) = \#\{j \mid \beta_j = -\alpha\} -
\#\{j \mid \beta_j = \alpha\}$.
\label{lem:chains=decompostions}
\end{lemma}  

\begin{proof}
Let $v=s_{i_1} \cdots s_{i_l}$ be a decomposition and 
$A_j = s_{i_1}\cdots s_{i_j} (A_\circ)$, for $j=0,\dots,l$.
Then $A_0=A_\circ$ and $A_l = v(A_\circ) = A$.  Applying 
$s_{i_1}\cdots s_{i_{j-1}}$ to the adjacent pair
$A_\circ\stackrel{-\alpha_{i_j}}\longrightarrow s_{i_j}(A_\circ)$,
we deduce that the pair 
$A_{j-1}\stackrel{-\beta_j}\longrightarrow A_j$
is adjacent as well,  
where $\beta_j = \bar s_{i_1} \cdots \bar s_{i_{j-1}}(\alpha_{i_j})$.
Thus $(A_0,\dots,A_l)$ is an alcove path from $A_\circ$ to $A$.
The reflection $s_{i_j}$ switches
the alcoves $A_\circ$ and $s_{i_j}(A_\circ)$.
Thus the reflection $r_j=s_{i_1}\cdots s_{i_j}\cdots s_{i_1}$ 
is the reflection with respect to the common wall of 
$A_{j-1}$ and $A_j$.

On the other hand, let  $(A_0,\dots,A_l)$ 
be any alcove path from $A_\circ$ to $A$,  and let $r_j$ be the reflection
with respect to the common wall of $A_{j-1}$ and $A_j$, for $j=1,\dots,l$.
Then $A_j = r_j\cdots r_1(A_\circ)$.  Applying 
$(r_{j-1}\cdots r_{1})^{-1}= r_1\cdots r_{j-1}$
to the adjacent pair $(A_{j-1}, A_j)$, 
we obtain the adjacent pair 
$(A_\circ, s(A_\circ))$,
where $s= r_1\cdots r_{j-1}r_j r_{j-1}\cdots r_1$.
Thus $s$ should be a reflection with respect to one of the walls 
of $A_\circ$. 
Thus there are $i_1,\dots,i_l\in\{0,\dots,r\}$ such that
$r_1\cdots r_{j-1}r_j r_{j-1}\cdots r_1 = s_{i_j}$, for $j=1,\dots,l$.
The affine Weyl group element $s_{i_1}\cdots s_{i_l} = r_l \cdots r_1$ 
maps $A_\circ$ to $A$, and is equal to $v$.

(a) $\Leftrightarrow$ (b).  \
This is clear, because a decomposition and the corresponding alcove path
have the same length.  

(b) $\Leftrightarrow$ (c). \ The fact that all affine reflections
$r_1,\dots,r_l$ are distinct for a reduced decomposition 
is given in~\cite[Lemma~4.5]{Hum}.
On the other hand, the length $l$ of any alcove path should
be at least the number of hyperplanes of the form $H_{\alpha,k}$ 
that separate $A_0$ and $A_l$.  If all affine reflections 
$r_1,\dots,r_l$ are distinct, then the path never crosses the same 
hyperplane twice, and, thus, its length equals the number of hyperplanes 
that separate $A_0$ and $A_l$.

(c)  $\Leftrightarrow$ (d). \ 
If $\beta_i = - \beta_j=\alpha$, then the alcove 
path crosses two parallel hyperplanes $H_{\alpha,k}$ and $H_{\alpha,l}$ 
in opposite directions.  It follows that the path crosses one of these 
hyperplanes twice, and, thus, the affine reflections $r_1,\dots,r_l$ are 
not distinct.
On the other hand, if $r_1,\dots,r_l$ are not distinct, then the path
crosses the same hyperplane more than once.  It follows that the path should 
cross this hyperplane in opposite directions.  Thus $\beta_i=-\beta_j$ 
for some $i$ and $j$.

The last claim follows from the fact that, each time
the alcove path crosses a hyperplane of the form 
$H_{\alpha,k}$, $\alpha\in\Phi^+$,
in positive (respectively negative) direction, the number $m_\alpha$
increases (respectively decreases) by 1, and all other $m_\beta$'s do not change.
\end{proof}

The affine translations by weights preserve the set
of affine hyperplanes $H_{\alpha,k}$, 
cf.~(\ref{eq:weight-lattice}) and~(\ref{eq:H-alpha-k}).
It follows that these affine translations map alcoves to alcoves.
Let $A_\lambda=A_\circ + \lambda$
be the alcove obtained by the affine translation of the fundamental alcove
$A_\circ$ by a weight $\lambda\in\Lambda$.  Let $v_\lambda=v_{A_\lambda}$
be the corresponding element of $\Waff$, i.e,.
$v_{\lambda}$ is defined by $v_\lambda(A_\circ) = A_\lambda$.
Note that the element $v_\lambda$ may not be an affine translation
itself.

\begin{definition}
Let $\lambda$ be a weight, and let 
$v_{-\lambda} = s_{i_1}\cdots s_{i_l}$ be any decomposition, reduced or not, 
of $v_{-\lambda}$ as a product of generators of $\Waff$.
Let us say that the {\it $\lambda$-chain of roots\/} associated with 
this decomposition is the sequence $(\beta_1,\dots,\beta_l)$ 
of the roots in $\Phi$ given by
$$
\beta_1 = \alpha_{i_1}, \  \beta_2=\bar s_{i_1}(\alpha_{i_2}), \ 
\beta_3 = \bar s_{i_1} \bar s_{i_2} (\alpha_{i_3}), \dots, \ 
\beta_l= \bar s_{i_1} \cdots \bar s_{i_{l-1}} (\alpha_{i_l})\,.
$$
Sometimes we will abbreviate ``$\lambda$-chain of roots'' as,
simply, ``$\lambda$-chain.''
Let us also say that the {\it $\lambda$-chain of reflections\/} associated 
with the above decomposition for $v_{-\lambda}$ is the sequence 
$(r_1,\dots,r_l)$ of the affine reflections in $\Waff$ given by
$$
r_1 = s_{i_1}, \ r_2 = s_{i_1} s_{i_2} s_{i_1}, \ 
r_3 =  s_{i_1} s_{i_2} s_{i_3} s_{i_2} s_{i_1}, \ \dots, \
r_l =  s_{i_1} \cdots s_{i_r} \cdots s_{i_1}. 
$$
In particular, $\bar r_i = s_{\beta_i}$.

According to Lemma~\ref{lem:chains=decompostions}, 
we can equivalently define a $\lambda$-chain as 
a sequence of roots $(\beta_1,\dots,\beta_l)$ 
such that there exists an alcove path
$A_0\stackrel{-\beta_1}\longrightarrow \cdots
\stackrel{-\beta_l}\longrightarrow A_{l}$
from $A_0=A_\circ$ to $A_l=A_{-\lambda}$
with edges labeled by the roots $-\beta_1,\dots,-\beta_l$.
The $j$-th element of the corresponding $\lambda$-chain of reflections 
$(r_1,\dots,r_l)$ is the affine reflection $r_j$ with respect to the common 
walls of the alcoves $A_{j-1}$ and $A_j$, for $j=1,\dots,l$.

Finally, we say that a $\lambda$-chain is {\it reduced\/} if it is associated
with a reduced decomposition for $v_{-\lambda}$.
\label{def:lambda-chain}
\end{definition}

\begin{remark}\label{remtransl}
If $A\stackrel{\beta}\longrightarrow B$
is a pair of adjacent alcoves, then 
$(A+\lambda)\stackrel{\beta}\longrightarrow (B+\lambda)$,
for any affine translation of the alcoves by the weight $\lambda$.
Thus, a translation of an alcove path by a weight $\lambda$ is an alcove
path labeled by the same sequence of roots.
For a $\lambda$-chain of roots $(\beta_1,\dots,\beta_l)$, let us translate
the corresponding alcove path $A_\circ
\stackrel{-\beta_1}\longrightarrow\cdots \stackrel{-\beta_l}\longrightarrow
A_{-\lambda}$ by 
the weight $\lambda$, and then reverse its direction.
We obtain the alcove path $A_\circ\stackrel{\beta_l}\longrightarrow \cdots
\stackrel{\beta_1}\longrightarrow A_{\lambda}$ associated with the 
$(-\lambda)$-chain $(-\beta_l,\dots,-\beta_1)$.
\end{remark}

\section{The $K_T$-Chevalley Formula}
\label{sec:K-chevalley-formula}

In this section, we formulate our main result and give 
its several specializations and applications to characters.

\begin{theorem}  {\rm ($K_T$-Chevalley formula) }
Fix any weight $\lambda$.  Let $(r_1,\dots,r_l)$ and 
$(\beta_1,\dots,\beta_l)$ be the $\lambda$-chain of reflections 
and the $\lambda$-chain of roots associated with a decomposition 
$v_{-\lambda} = s_{i_1}\cdots s_{i_l}\in\Waff$, which may or may 
not be reduced.
Let $u,w\in W$, and $\mu\in\Lambda$.
Then the $K_T$-Chevalley coefficient $c_{u,w}^{\lambda,\mu}$,
i.e., the coefficient of $x^\mu\,[\O_w]$ in the expansion of the product 
$e^\lambda\cdot [\O_u]$, can be expressed as follows:
\begin{equation}
c_{u,w}^{\lambda,\mu}=\sum_J(-1)^{n(J)}\,,
\label{eq:cuwlambdamu}
\end{equation}
where the summation is over
all subsets $J=\{j_1<\cdots<j_s\}$ of $\{1,\dots,l\}$ satisfying the following conditions: 
\begin{enumerate}
\item[(a)]
$u \gtrdot u\, \bar r_{j_1} \gtrdot 
u \,\bar r_{j_1} \bar r_{j_2} \gtrdot \cdots \gtrdot 
u \,\bar r_{j_1} \bar r_{j_2} \cdots \bar r_{j_s} = w$
is a saturated decreasing chain from $u$ to $w$ in the Bruhat order on 
the Weyl group $W$;
\item[(b)]
$-\mu = u \,r_{j_1} \cdots r_{j_s}(-\lambda)$,
\end{enumerate}
and $n(J)$ is the number of negative roots in 
$\{\beta_{j_1},\dots,\beta_{j_s}\}$.
\label{th:main-combinatorial-formulation}
\end{theorem}

In Section~\ref{sec:k-chevalley-operator-notation},
we reformulate this theorem in a compact form
and then prove it, using a certain $R$-matrix.
In Sections~\ref{sec:Type-A}
and~\ref{sec:other-types}, we give several examples
that illustrate this theorem.

Lemma~\ref{lem:chains=decompostions} implies the following statement.
\begin{lemma}
Let $(\beta_1,\dots,\beta_l)$ be a reduced $\lambda$-chain of roots.
Let $\alpha\in\Phi$ be a root such that $(\lambda,\alpha^\vee)\geq 0$. 
Then $\#\{i \mid \beta_i=\alpha\}
=(\lambda,\alpha^\vee)$ and $\#\{i \mid \beta_i=-\alpha\}=0$.

In particular, if $\lambda$ is a dominant weight,
then all roots $\beta_1,\dots,\beta_l$ are positive.
Also, if $\lambda$ is an antidominant weight, that is, $-\lambda\in\Lambda^+$,
then all roots $\beta_1,\dots,\beta_l$ are negative.
\label{lem:reduced-dominant-antidominant}
\end{lemma}

In the special cases corresponding to dominant and antidominant weights
$\lambda$, Theorem~\ref{th:main-combinatorial-formulation} can be reformulated
in a more explicit way. In these cases, for reduced $\lambda$-chains,
Theorem~\ref{th:main-combinatorial-formulation} gives
a manifestly positive formula, which is not the case in general.

\begin{corollary}
Consider the setup in Theorem~{\rm \ref{th:main-combinatorial-formulation}}.
Assume that $v_{-\lambda}=s_{i_1}\cdots s_{i_l}$ is a reduced decomposition
in $\Waff$. 

If $\lambda$ is a dominant weight, then $c_{u,w}^{\lambda,\mu}$
equals the number of subsets $J\subseteq\{1,\dots,l\}$ that satisfy conditions
{\rm(a)} and~{\rm(b)} in Theorem~{\rm \ref{th:main-combinatorial-formulation}}.

If $\lambda$ is an antidominant weight, then $(-1)^{\ell(u)-\ell(w)}\,c_{u,w}^{\lambda,\mu}$
equals the number of subsets $J\subseteq\{1,\dots,l\}$ that satisfy conditions
{\rm(a)} and~{\rm(b)} in Theorem~{\rm \ref{th:main-combinatorial-formulation}}.
\label{cor:dominant-antidominant}
\end{corollary}

\begin{proof}
For a dominant weight $\lambda$, all roots $\beta_1,\dots,\beta_l$ 
are positive; thus $n(J)=0$.
For an antidominant weight $\lambda$, all roots $\beta_1,\dots,\beta_l$ 
are negative; thus $n(J)=|J|=\ell(u)-\ell(w)$.
\end{proof}

Theorem~\ref{th:main-combinatorial-formulation} specializes to
following rule for products in the (nonequivariant) Grothendieck ring 
$K(G/B)$.

\begin{corollary}
The coefficient $c_{u,w}^\lambda$ of $[\O_w]$ in the 
product $e^\lambda\cdot [\O_u]$ of classes in $K(G/B)$ 
has the same combinatorial 
description as in Theorem~{\rm \ref{th:main-combinatorial-formulation}}, 
except that condition~{\rm(b)} on the weights involved is dropped.
\end{corollary}

\begin{proof}
We have $c_{u,w}^\lambda=\sum_{\mu\in\Lambda} c_{u,w}^{\lambda,\mu}$.
\end{proof}

Theorem~\ref{th:main-combinatorial-formulation} implies 
the following combinatorial model for the Demazure 
characters $ch(V_{\lambda,u})$ and, in particular, 
for the characters $ch(V_\lambda)$
of the irreducible representations $V_\lambda$ of the Lie group $G$.

\begin{corollary}
Let $\lambda$ be a dominant weight, let $u\in W$,
and let $(r_1,\ldots,r_l)$ be a reduced $\lambda$-chain of reflections.
Then the Demazure character 
$ch(V_{\lambda,u})$ 
is equal to the sum
$$
ch(V_{\lambda,u}) =
\sum_J e^{-u\, r_{j_1}\cdots r_{j_s}(-\lambda)}
$$
over all subsets $J=\{j_1<\cdots<j_s\}\subset\{1,\dots,l\}$ such that
$$
u \gtrdot u\, \bar r_{j_1} \gtrdot 
u \,\bar r_{j_1} \bar r_{j_2} \gtrdot \cdots \gtrdot 
u \,\bar r_{j_1} \bar r_{j_2} \cdots \bar r_{j_s} 
$$
is a saturated decreasing chain in the Bruhat order on 
the Weyl group $W$.
\label{cor:model-for-demazure-character}
\end{corollary}

\begin{proof}
Apply Corollary~\ref{cor:dominant-antidominant} and
Lemma~\ref{lem:demazure-character}.  
\end{proof}

We can slightly simplify the formula for the characters 
$ch(V_\lambda) = ch(V_{\lambda,w_\circ})$ 
of the irreducible representations of $G$, as follows. 

\begin{corollary}  Consider the setup in 
Corollary~{\rm \ref{cor:model-for-demazure-character}}. 
We have
$$
ch(V_\lambda) =
\sum_J e^{-r_{j_1}\cdots r_{j_s}(-\lambda)}\,,
$$
where the summation is over all subsets $J=\{j_1<\cdots<j_s\}\subset\{1,\dots,l\}$ such that
$$
1 \lessdot \bar r_{j_1} \lessdot \bar r_{j_1} \bar r_{j_2} \lessdot 
\cdots \lessdot \bar r_{j_1} \bar r_{j_2} \cdots \bar r_{j_s} 
$$
is a saturated increasing chain in the Bruhat order on 
the Weyl group $W$.
\label{cor:model-for-character}
\end{corollary}

\begin{proof}
Multiplying elements in a decreasing chain by $w_\circ$ on the left
results in an increasing chain in Bruhat order.  On the other hand, we can remove 
$w_\circ$ from
the exponent because the character $ch(V_\lambda)$ is $W$-invariant.
\end{proof}

In the rest of this section, we show how to construct 
$\lambda$-chains of reflections $(r_1,\dots,r_l)$ and 
$\lambda$-chains of roots $(\beta_1,\dots,\beta_l)$.
Clearly, there are many possible choices.

Let us fix an arbitrary weight $\lambda$.  Let $\pi:[0,1]\to \h_\R^*$ 
be a sufficiently generic continuous path such that 
$\pi(0)\in A_\circ$ and $\pi(1)\in A_{-\lambda}$.
Here ``sufficiently generic'' means that the path $\pi$ does not cross any
face of an alcove of codimension 2 or higher. 
For example, the path $\pi:t\mapsto -t\, \lambda  + \gamma$, 
where $\gamma$ is a generic point in $A_\circ$, will suffice.
Suppose that the path $\pi$ passes through the sequence of alcoves
$A_\circ,\dots,A_{-\lambda}$ as $t$ varies from 0 to 1.
This sequence is an alcove path.
Let $H_1,\dots,H_l$ be the affine hyperplanes of the form
$H_{\alpha,k}$ that the path $\pi$ crosses as $t$ varies from 0 to 1.
According to Lemma~\ref{lem:chains=decompostions}, the sequence 
$(r_1,\dots,r_l)$ of affine reflections with respect to 
$H_1,\dots,H_l$ is a $\lambda$-chain of reflections.

In order to make our formula completely combinatorial, we 
present one particular choice for a $\lambda$-chain of reflections
and the corresponding $\lambda$-chain of roots.
The construction  depends on the choice of a total order 
$\alpha_1<\cdots < \alpha_r$ on the simple roots in $\Phi$.
Suppose that $\pi=\pi_{\varepsilon}:[0,1]\to \h_\R^*$ 
is the path given by 
$$
\pi_{\varepsilon}:t\mapsto -t\,\lambda + \varepsilon\, \omega_1 +
\varepsilon^2 \omega_2 +\cdots +\varepsilon^r \omega_r,
$$
where $\varepsilon$ is a sufficiently small positive real number.
Let $\mathcal{R}=\mathcal{R}_\lambda\subset\Waff$ be the set of affine 
reflections with respect to affine hyperplanes $H_{\alpha,k}$ that 
separate the alcoves $A_\circ$ and $A_{-\lambda}$.  This set is given by
$$
\mathcal{R}= 
\mathcal{R}_\lambda= 
\bigcup_{\alpha\in\Phi^+}
\left\{
\begin{array}{cl} 
\{s_{\alpha,k} \mid 0\geq k > -(\lambda,\alpha^\vee)\}&
\textrm{if } (\lambda,\alpha^\vee)>0\,, \\[.05in]
\{s_{\alpha,k} \mid 0< k \leq -(\lambda,\alpha^\vee)\}&
\textrm{if } (\lambda,\alpha^\vee)<0\,, \\[.05in]
\emptyset &
\textrm{if }(\lambda,\alpha^\vee)=0\,.
\end{array}
\right.
$$
For any $s_{\alpha,k}\in\mathcal{R}$, $\alpha\in\Phi^+$,
the path $\pi_\varepsilon$ crosses the affine hyperplane $H_{\alpha,k}$ 
at the point $t=t_{\alpha,k}=(\lambda,\alpha^\vee)^{-1}
(-k+\sum_{i=1}^r(\omega_i,\alpha^\vee)\,\varepsilon^i)$.
Note that $(\lambda,\alpha^\vee)\ne 0$, for $s_{\alpha,k}\in\mathcal{R}$.
Let $h:\mathcal{R}\to\R^{r+1}$ be the map given by 
\begin{equation}
h:s_{\alpha,k}\mapsto 
(\lambda,\alpha^\vee)^{-1}\, (-k,(\omega_1,\alpha^\vee),\dots,
(\omega_r,\alpha^\vee)),
\label{eq:map-h}
\end{equation}
for any $s_{\alpha,k}\in\mathcal{R}$ with $\alpha\in\Phi^+$.
Then, for sufficiently small $\varepsilon>0$, we have 
$t_{\alpha,k}<t_{\alpha',k'}$ 
if and only if $h(s_{\alpha,k})$
is less than $h(s_{\alpha',k'})$ in the lexicographic order on $\R^{r+1}$.
We claim that the map $h$ is injective.  Indeed, if $h(s_{\alpha,k})=
h(s_{\alpha',k'})$, then $\alpha=\alpha'$.  Otherwise, 
the root system $\Phi^\vee$ would contain two proportional positive coroots 
$\alpha^\vee\ne(\alpha')^\vee$, which is not possible.
Also, the fact that $\alpha=\alpha'$ implies that $k=k'$.

Let $b:\{\textrm{affine reflections}\}\to\Phi$ be the map given by
$$
b: s_{\alpha,k}\longmapsto \left\{
\begin{array}{cl}
\alpha&\textrm{if } k\leq 0 \textrm{ and }\alpha\in\Phi^+,\\[.05in]
-\alpha&\textrm{if } k>0
\textrm{ and }\alpha\in\Phi^+.
\end{array}
\right.
$$

We obtain the following result by using Lemma~\ref{lem:chains=decompostions}.

\begin{proposition}
Let $\mathcal{R}=\{r_1<r_2<\cdots < r_l\}$ be the total order on the set 
$\mathcal{R}$ such that
$h(r_1)<h(r_2)<\cdots<h(r_l)$ in the lexicographic order on $\R^{r+1}$.
Then $(r_1,\dots,r_l)$ is the $\lambda$-chain of reflections
and $(\beta_1,\dots,\beta_l)=(b(r_1),\dots,b(r_l))$ is the  
$\lambda$-chain of roots associated with a certain reduced decomposition 
of $v_{-\lambda}$. 
\label{prop:construction-for-lambda-chain}
\end{proposition}

Example~\ref{example:G-2} illustrates this proposition.

\section{Generalization to $G/P$}
\label{sec:parabolic}

Let $P$ be a parabolic subgroup in $G$ such that $P\supset B$.
In this section, we show that the $K_T$-Chevalley formula can
be easily extended to equivariant $K$-theory of 
the {\it generalized partial flag variety\/} $G/P$.
\medskip

Let $\Delta_P$ be the subset of the simple roots associated with
the parabolic subgroup $P$.
Let $\Phi_P\subset \Phi$ be the set of roots that
can be written as sums of roots in $\Delta_P$, and let 
$\Phi_P^+ = \Phi_P\cap \Phi^+$.
Then $\Phi_P$ is a root system itself, with the Weyl group $W_P\subset W$
generated by the simple reflections $s_i$, for $\alpha_i\in \Delta_P$.  Each
coset $\bar w = wW_P$ in $W/W_P$ has a unique representative of maximal
length.   Let us denote the set of maximal coset representatives by
$W^P\subset W$, and let us identify it with $W/W_P$.  The Bruhat order on $W$
induces the Bruhat order on $W^P\simeq W/W_P$.  According to
Deodhar~\cite{Deo1}, the covering relations in $W^P$ are of the form $u\gtrdot
w$, where $w=us_\beta$, for some $\beta\in\Phi^+\setminus\Phi_P^+$, and
$\ell(u)=\ell(w)+1$.  In particular, every covering relation in $W^P$ is a
covering relation in the Bruhat order on $W$.

The generalized partial flag variety $G/P$ decomposes
into Schubert cells $X_{\bar w}^\circ=B\bar wP/P$ indexed by $\bar w\in W/W_P$.
Their closures $X_{\bar w} := \overline{X_{\bar w}^\circ}$
are called Schubert varieties.
Let $\O_{X_{\bar w}}^P$, $\bar w\in W/W_P$, be the structure sheaf
of the Schubert variety $X_{\bar w}$.
If $\lambda$ is a weight satisfying $(\lambda,\beta)=0$, 
for all $\beta$ in $\Delta_P$ (or, equivalently,
$W_P\subseteq W_\lambda$, where $W_\lambda$ is the stabilizer of $\lambda$),
then $-\lambda$ determines a character of $P$, and thus a line bundle
$\mathcal{L}_\lambda^P := G\times_P \C_{-\lambda}$ on $G/P$. 
Let $[\O_{X_{\bar w}}^P]$ and $[\L_\lambda^P]$ be the corresponding classes
in $K_T(G/P)$.
The classes $[\O_{X_{\bar w}}^P]$ form a $\Z[X]$-basis of $K_T(G/P)$,
and the classes $[\L_\lambda^P]$ span $K_T(G/P)$ over $\Z[X]$. Let $[\O_{\bar w}^P]:=*[\O_{X_{\bar w}}^P]$, where the involution $*$ on $K_T(G/P)$ is defined like the one on $K_T(G/B)$.

The equivariant $K$-theory of $G/P$ 
can be recovered from $K_T(G/B)$, as stated in \cite{KK}. 
We have the canonical projection $\pi_P\::\:G/B\rightarrow G/P$. This
determines an injective $\Z[X]$-linear homomorphism
$\pi_P^*\::\:K_T(G/P)\rightarrow K_T(G/B)$. Moreover, the image of this map,
with which $K_T(G/P)$ can be identified, consists precisely of the
$W_P$-invariants in $K_T(G/B)$.  It is straightforward to
show that 
\begin{equation}\label{eq:proj}
\pi_P^*([\O_{\bar w}^P])=[\O_{w}]\,,\quad
\textrm{and}\quad\pi_P^*([{\mathcal L}_\lambda^P])=[{\mathcal L}_\lambda]\,,
\end{equation}
where $w\in W^P$ is the maximal coset representative of $\bar w\in W/W_P$,
and the weight $\lambda$ is such that $W_P\subseteq W_\lambda$. By abuse of notation, we will denote the class $[{\mathcal L}_{-\lambda}^P]$ in $K_T(G/P)$ by $e^\lambda$, as well.

Let us define the integer coefficients $c_{\bar u,\bar w}^{\lambda,\mu}$,
for $\bar u,\bar w\in W/W_P$ and $\lambda,\mu\in \Lambda$,
with $W_P\subseteq W_\lambda$, by the following expansion
of the product in $K_T(G/P)$:
\begin{equation}\label{parchev}
e^\lambda\cdot [\O_{\bar u}^P] = \sum_{\bar w\in W/W_P,\,\mu\in\Lambda} 
c_{\bar u, \bar w}^{\lambda,\mu} \ x^\mu\,[\O_{\bar w}^P]\,.
\end{equation}

Our combinatorial Chevalley-type formula for $K_T(G/B)$ can be 
generalized to $K_T(G/P)$, as follows. 

\begin{corollary} 
Let $u,w\in W^P$ be the maximal coset representatives 
of $\bar u,\bar w\in W/W_P$, and let $\lambda,\mu\in \Lambda$ such that
$W_P\subseteq W_\lambda$.
Then we have 
$c_{\bar u,\bar w}^{\lambda,\mu}= c_{u,w}^{\lambda,\mu}$, 
where $c_{u,w}^{\lambda,\mu}$ is the $K_T$-Chevalley coefficient
for $K_T(G/B)$, which have the combinatorial description given in
Theorem~{\rm \ref{th:main-combinatorial-formulation}}. Moreover, if we work
with reduced $\lambda$-chains, then all the elements of the corresponding
saturated chains in the Bruhat order lie in $W^P$.  
\end{corollary}

\begin{proof} The first part of the proof is immediate by applying the map
$\pi_P^*$ to both sides of~(\ref{parchev}),
and by using~(\ref{eq:proj}). The second statement follows from the fact that,
given the choice of $\lambda$, we have $(\lambda,\beta^\vee)=0$, for all $\beta$
in $\Phi_P$. Indeed, by Lemma~\ref{lem:chains=decompostions}, a reduced
$\lambda$-chain of roots does not contain any roots in $\Phi_P$. 
Therefore, the conclusion follows from the above description of the 
Bruhat order on $W^P$.  
\end{proof}

\section{Applications: $K_T$-Monk Formula and 
Duality Formulas}
\label{sec:duality}

In this section, we present several applications of our $K_T$-Chevalley formula.
First, we give a rule for products $[\O_{w_\circ s_i}]\cdot [\O_u]$,
which we call the {\it $K_T$-Monk formula.}
We also give the {\it dual $K_T$-Chevalley formula\/} for products 
$e^\lambda\cdot [\I_u]$.
Then we derive two {\it duality formulas\/} 
for the $K_T$-Chevalley coefficients.
The first one has been already stated for $K(G/B)$, in a slightly imprecise
way, by Brion in~\cite[Theorem 4]{Bri}, and proved using some fairly involved
geometric arguments.  We present a concise combinatorial proof, based on our
$K_T$-Chevalley formula.   The two dualities came from the two involutions
$w\mapsto w w_\circ$ and $w\mapsto w_\circ w$ on $W$.  Our $K_T$-Chevalley
formula is symmetric with respect to these involutions, because they map
increasing chains in the Bruhat order to decreasing chains.
\medskip

Let us call the classes $[\O_{w_\circ s_i}]\in K_T(G/B)$ the {\it special classes}; they correspond to the structure sheaves of codimension one Schubert varieties $X_{w_\circ s_i}$,

\begin{lemma}
{\rm(a) \cite{Bri} } 
For a simple reflection $s_i$, we have 
$$
[\O_{w_\circ s_i}] = 1 - x^{w_\circ(\omega_i)} e^{-\omega_i}
$$
in the Grothendieck ring $K_T(G/B)$.

{\rm(b)} The special classes $[\O_{w_\circ s_i}]$, $i=1,\dots,r$, generate
the Grothendieck ring $K_T(G/B)$ as an algebra over $\Z[X]$.
\label{lem:O=L}
\end{lemma}

Brion proved that $[\O_{X_{w_\circ s_i}}] = 1 - [\L_{-\omega_i}]$
in $K(G/B)$ using a simple geometric argument 
based on the exact sheaf sequence
$0\to \L_{-\omega_i} \to \O_{G/B}\to \O_{X_{w_\circ s_i}} \to 0$.
Brion also mentioned that this argument extends 
to $T$-equivariant $K$-theory.

\begin{proof}
(a)
Let us apply  Theorem~\ref{th:main-combinatorial-formulation}, for 
$u=w_\circ$ and $\lambda=-\omega_i$. 
Every saturated chain in the Bruhat order decreasing from $w_\circ$
should start with a simple reflection.
For a reduced $(-\omega_i)$-chain of reflections $(r_1,\dots,r_l)$,
exactly one of the reflections $\bar r_1, \dots,\bar r_l$ is simple.
Namely, $\bar r_l = s_i$ and, moreover, $r_l = s_{\alpha_i,1}$.
Thus the expansion of the product $e^{-\omega_i}\cdot [\O_{w_\circ}]$
consists of the two terms corresponding to the subsets
$J=\emptyset$ and $J=\{l\}$.  This expansion is
$e^{-\omega_i}\cdot [\O_{w_\circ}] = 
x^{-w_\circ(\omega_i)} [\O_{w_\circ}] - 
x^{-w_\circ(\omega_i)} [\O_{w_\circ s_i}]. 
$
Since $[\O_{w_\circ}]=1$, we obtain the required identity.

(b)  Let us identify $K_T(G/B)$ with the quotient in~(\ref{eq:KGB=ZL/I}).  
There is a {\it finite\/} set $D$ of exponents $e^\mu$ that spans $K_T(G/B)$ 
as a $\Z[X]$-module.  Indeed, we can take all exponents in 
some representatives for the classes $[\O_w]$ in $\Z[X]\otimes\Z[\Lambda]$.
For a weight $\lambda\in\Lambda$, the exponent $e^\lambda$ is an invertible 
element in $K_T(G/B)$; and, thus, the set 
$e^\lambda D = \{e^{\lambda+\mu}\mid e^\mu\in D\}$ also spans $K_T(G/B)$.  
For a sufficiently large antidominant weight $\lambda$,
all exponents in the set $e^\lambda D$ correspond to antidominant weights.
On the other hand, according to~(a), we have
$e^{-\omega_i} = x^{-w_\circ(\omega_i)}(1-[\O_{w_\circ s_i}])$;
thus, all classes $e^\mu = [\L_{-\mu}]$, for antidominant weights $\mu$, 
can be expressed in terms of the special classes $[\O_{w_\circ s_i}]$.
This implies the statement.
\end{proof}

The second part of Corollary~\ref{cor:dominant-antidominant},  
for $\lambda=-\omega_i$, and Lemma~\ref{lem:O=L}(a) 
imply the combinatorial rule below for products of the special classes
with the basis elements in $K_T(G/B)$. Note that, if $\omega_i$ is a minuscule weight (i.e., $(\omega_i,\alpha^\vee)=0\textrm{ or }1$ for any $\alpha\in\Phi^+$), then all reflections $r_j$ in a reduced $(-\omega_i)$-chain of reflections have the form $r_j=s_{\beta_j,1}$, and therefore they all fix $\omega_i$.

\begin{corollary} {\rm ($K_T$-Monk formula) }
Fix a simple reflection $s_i$,
and let $(r_1,\dots,r_l)$ be a reduced $(-\omega_i)$-chain of reflections.
Then, for any $u\in W$, we have 
$$
[\O_{w_\circ s_i}]\cdot [\O_u] = 
(1-x^{w_\circ(\omega_i) - u(\omega_i)}) \, [\O_u] +
\sum_J
(-1)^{|J|-1} \,x^{\nu(J)} \,[\O_{w(J)}],
$$
where the sum is over nonempty subsets $J=\{j_1,\dots,j_s\}$
in $\{1,\dots,l\}$ such that
$u \gtrdot u\, \bar r_{j_1} \gtrdot 
u \,\bar r_{j_1} \bar r_{j_2} \gtrdot \cdots \gtrdot 
u \,\bar r_{j_1} \bar r_{j_2} \cdots \bar r_{j_s} = w$ 
is a saturated decreasing chain in the Bruhat order from $u$ to
$w=w(J)$, and $\nu(J) = 
w_\circ(\omega_i) - u\, r_{j_1}\cdots r_{j_s}(\omega_i)$.
\label{cor:OiO=O}
If $\omega_i$ is minuscule, then the above formula has the simpler form:
$$
[\O_{w_\circ s_i}]\cdot [\O_u]= [\O_u]+ x^{w_\circ(\omega_i) - u(\omega_i)} \left(
\sum_J
(-1)^{|J|-1} \,[\O_{w(J)}]\right),
$$
where the notation is as above, but we drop the condition $J\ne\emptyset$.
\end{corollary}

Since the special classes $[\O_{w_\circ s_i}]$ generate the Grothendieck ring 
$K_T(G/B)$, Corollary~\ref{cor:OiO=O} completely characterizes
the multiplicative structure of this ring.

\begin{remark}
In the equivariant case,
the expansion of $[\O_{w_\circ s_i}]\cdot [\O_u]$ contains
the term $[\O_u]$ with a nonzero coefficient.
This term vanishes in the nonequivariant case of $K(G/B)$.
A similar phenomenon happens in the Monk-type formula for 
equivariant cohomology,
which can be derived from Corollary~\ref{cor:OiO=O}.
\end{remark}

\medskip

Recall that the classes $[\I_w]$, $w\in W$, given 
by~(\ref{eq:I-dual-basis}) form the dual basis to $\{[\O_w]\mid w\in W\}$
with respect to the natural pairing in $K$-theory.
Define the {\it dual $K_T$-Chevalley coefficients\/} $d_{u,w}^{\lambda,\mu}$,
for $u,w\in W$, $\lambda,\mu\in \Lambda$, by the expansion
$$
e^\lambda\cdot [\I_u] = \sum_{w\in W,\, \mu\in \Lambda} 
d_{u,w}^{\lambda,\mu}\, x^\mu\,[\I_w].
$$

\begin{corollary}  {\rm (dual $K_T$-Chevalley formula) }
The dual $K_T$-Chevalley coefficients are related to the $K_T$-Chevalley 
coefficients as $d_{u,w}^{\lambda,\mu} = c_{u,w}^{-\lambda,-\mu}$.
Thus Theorem~{\rm \ref{th:main-combinatorial-formulation}} provides
a combinatorial description for the coefficients
$d_{u,w}^{\lambda,\mu}$.
\end{corollary}

\begin{proof}
Follows from~(\ref{eq:L-I}).
\end{proof}

\begin{remark}
In a recent paper\footnote{\cite{GR} appeared in {\tt arXiv} after the present
paper was finished.}, Griffeth and Ram~\cite{GR} provided more details
of the proof of the Pittie-Ram formula and gave a dual $K_T$-Chevalley
formula, for dominant weights $\lambda$, using LS-paths.  
They also derived Lemma~\ref{lem:O=L}(a) above and 
Theorem~\ref{th:first_duality} below, for dominant $\lambda$.
Note that our dual $K_T$-Chevalley formula is just the usual $K_T$-Chevalley 
formula (Theorem~{\rm \ref{th:main-combinatorial-formulation}})
with $\lambda$ and $\mu$ replaced by $-\lambda$ and $-\mu$.
Since the Pittie-Ram formula does not work for nondominant weights,
Griffeth and Ram had to derive their dual version separately.
The symmetry between the Pittie-Ram formula and its dual version
given in~\cite{GR} is not so transparent as the symmetry in our construction.
In fact, Griffeth and Ram gave four different formulas for 
the products $e^\lambda\cdot [\O_w]$, $e^{-\lambda}\cdot [\O_{w}]$,
$e^{w_\circ(\lambda)}\cdot [\O_w]$, and $[\O_{w_\circ s_i}]\cdot [\O_{w}]$,
for a dominant weight $\lambda$, using LS-paths.
>From our point of view, these four products are given by various 
specializations of the $K_T$-Chevalley formula, for arbitrary $\lambda$. 
\end{remark}

Let us now discuss symmetries of the $K_T$-Chevalley coefficients.
In order to make our notation compatible with that in \cite{Bri}, 
we define the coefficients $c_u^w(\lambda)$ in $\Z[X]$ by
$$
e^\lambda\cdot [\O_u] = \sum_{w\in W} 
c_{u}^{w}(\lambda) \,[\O_w]\,.
$$
In other words, the $c_u^w(\lambda)$ are expressed in terms of the 
$K_T$-Chevalley coefficients, as follows: 
$c_u^w(\lambda)=\sum_{\mu\in\Lambda}c_{u,w}^{\lambda,\mu}\,x^\mu$,
see~(\ref{eq:LlOu}). 

\begin{theorem}\label{duality} {\rm \cite[Theorem~4]{Bri} }
We have the following duality formula for an arbitrary weight $\lambda$:
$$
c_u^w(\lambda)=
(-1)^{\ell(u)-\ell(w)}c_{ww_\circ}^{uw_\circ}(w_\circ\lambda)\,.
$$
\label{th:first_duality}
\end{theorem}

\begin{proof}
Let $(\beta_1,\ldots,\beta_l)$ and $(r_1,\ldots,r_l)$ be the $\lambda$-chain of
roots and the $\lambda$-chain of reflections associated with some alcove
path.  Let us translate this alcove path by $\lambda$, 
reverse its direction (cf.\ Remark~\ref{remtransl}), and then apply the map 
$A\mapsto -w_\circ(A)$ to the corresponding alcoves.  Note that
$-w_\circ(A_\circ)=A_\circ$.
The resulting alcove path corresponds to
the $(w_\circ \lambda)$-chain of roots 
$(w_\circ\beta_l,\ldots,w_\circ\beta_1)$ and a certain 
$w_\circ(\lambda)$-chain of reflections $(r_l',\ldots,r_1')$.
We can express the affine reflections $r_j'$, as follows.
Let $\gamma$ and $t_\lambda$ be  
the operators on $\hR$ given by $\gamma:\mu\mapsto -\mu$ and
$t_\lambda:\mu\mapsto\mu+\lambda$.
Then $r_j' = w_\circ \gamma\, t_\lambda r_j t_{-\lambda} \gamma\, w_\circ $.
Thus $\bar r_j' = w_\circ \bar r_j w_\circ$.

Clearly, to each sequence $J=(j_1,j_2,\dots,j_s)$ with 
$$
u\gtrdot u\bar{r}_{j_1}\gtrdot u\bar{r}_{j_1}
\bar{r}_{j_2}\gtrdot\cdots\gtrdot u\bar{r}_{j_1}\bar{r}_{j_2}
\cdots\bar{r}_{j_s}=w\,,
$$
corresponds the sequence $J'=(j_s,j_{s-1},\dots,j_1)$ with 
$$
ww_\circ\gtrdot ww_\circ\bar{r}_{j_s}'\gtrdot ww_\circ\bar{r}_{j_s}'\bar{r}_{j_{s-1}}'\gtrdot\cdots\gtrdot 
ww_\circ\bar{r}_{j_s}'\bar{r}_{j_{s-1}}'
\cdots\bar{r}_{j_1}'=uw_\circ\,.
$$
This correspondence is a bijection. 
Since $w_\circ$ maps positive roots to negative roots, we have 
$n(J')=s-n(J)=\ell(u)-\ell(w)-n(J)$, so
$(-1)^{n(J)}=(-1)^{\ell(u)-\ell(w)}(-1)^{n(J')}$. This takes care of the sign
in the duality formula. 

It remains to check that the sequences $J$ and $J'$ produce the same
weight, see condition (b) in Theorem~\ref{th:main-combinatorial-formulation}.
It suffices to show that
$$
r_{j_1}r_{j_2}\ldots r_{j_s}(-\lambda)=\bar{r}_{j_1}\bar{r}_{j_2}\ldots \bar{r}_{j_s}w_\circ r_{j_s}'r_{j_{s-1}}'\ldots r_{j_1}'w_\circ (-\lambda)\,.
$$
Let us denote $v=r_{j_1}\cdots r_{j_s}\in\Waff$.
Then the left-hand side of this expression is $v(-\lambda)$.
We can write the right-hand side of this expression as
$$
\bar r_{j_1}\cdots \bar r_{j_s} \,\gamma\, t_\lambda 
r_{j_s}\cdots r_{j_1} t_{-\lambda}\gamma\,(-\lambda)
= - \bar v \,t_\lambda v^{-1}(0).
$$
We claim that 
\begin{equation}
v(-\lambda)=-\bar v\, t_\lambda \, v^{-1}(0),
\label{eq:vlambdabarv}
\end{equation}
for any $v\in\Waff$ and $\lambda\in\Lambda$.
Indeed, if $v(-\lambda) = \bar v(-\lambda) +\mu$, then
$v^{-1}(0) = \bar v^{-1}(0-\mu) = -\bar v^{-1} (\mu)$.
Thus $\bar v\, t_\lambda v^{-1}(0) = \bar v(\lambda)- \mu$, as needed.
\end{proof}

Let us also present a new duality formula. We denote by $\iota$ the involutory
automorphism of $\Z[X]$ given by $\iota: x^\mu\mapsto x^{-w_\circ\mu}$.

\begin{theorem}
We have the following duality formula for an arbitrary weight $\lambda$:
\[c_u^w(\lambda)=(-1)^{\ell(u)-\ell(w)}\iota(c_{w_\circ w}^{w_\circ u}(-\lambda))\,.\]
\label{th:second_duality}
\end{theorem}

\begin{proof}
Let $(\beta_1,\ldots,\beta_l)$ and $(r_1,\ldots,r_l)$ be the $\lambda$-chain of
roots and the $\lambda$-chain of reflections associated with some alcove
path.  Let us translate the alcove path and reverse its direction,
as discussed in Remark~\ref{remtransl}.  We obtain the $(-\lambda)$-chain 
of roots $(-\beta_l,\ldots,-\beta_1)$
and the corresponding $(-\lambda)$-chain of roots
$(r_l',\ldots,r_1')$. 
Let $t_\lambda$ be the operator of translation by $\lambda$, as before.
Then $r_j' = t_\lambda r_j t_{-\lambda}$.
Thus $\bar r_j' = \bar r_j$.  In an almost identical way to the proof of
Theorem~\ref{duality}, we can now construct a bijection between the appropriate
decreasing saturated chains from $u$ to $w$, and those from $w_\circ w$ to
$w_\circ u$. The discussion about the signs is also similar. 
It remains to verify the weight condition: 
$$
r_{j_1}r_{j_2}\cdots r_{j_s}(-\lambda)=-\bar{r}_{j_1}\bar{r}_{j_2}\cdots
\bar{r}_{j_s} r_{j_s}'r_{j_{s-1}}'\cdots r_{j_1}' (\lambda)\,.  
$$
This identity can be written as $v(-\lambda) = - \bar v\, t_\lambda v^{-1}
t_{-\lambda} (\lambda)$, 
for $v = r_{j_1}\cdots r_{j_s}$,
which is equivalent 
to~(\ref{eq:vlambdabarv}).
\end{proof}

The two duality formulas above imply the following formula. 

\begin{corollary} 
Given an arbitrary weight $\lambda$, we have 
\[c_u^w(\lambda)=\iota(c_{w_\circ u w_\circ}^{w_\circ w w_\circ}
(-w_\circ\lambda))\,.\]
\label{cor:third_duality}
\end{corollary}

Note each of the two duality formulas in Theorems~\ref{th:first_duality}
and~\ref{th:second_duality} can be obtained from the other one
combined with Corollary~\ref{cor:third_duality}.

Kumar provided us with the following geometric explanation of 
Corollary~\ref{cor:third_duality}.
This duality in equivariant $K$-theory is induced by the standard involution 
on $G/B$, which interchanges the Schubert varieties $X_w$ and 
$X_{w_\circ w w_\circ}$.  
Let us denote by $\theta$ the canonical isomorphism (\ref{eq:KGB=ZL/I}) from
$(\Z[X]\otimes\Z[\Lambda])/\mathcal{I}$ to $K_T(G/B)$.

\begin{proposition}
\label{kumar}
There is an involutive automorphism $\omega$ on $K_T(G/B)$ 
such that
\begin{itemize}
\item[(a)] the involution $\omega$
maps each class $[\O_{w}]$ to $[\O_{w_\circ w w_\circ}]$; 
\item[(b)] 
under the isomorphism $\theta$, the involution $\omega$
maps $x^\mu\otimes e^\lambda$ to  $x^{-w_\circ\mu} \otimes e^{-w_\circ\lambda}$,
for $\lambda,\mu\in\Lambda$.
\end{itemize}
\end{proposition}

\begin{proof}[Algebraic proof]  
The involutive automorphism of $\Z[X]\otimes\Z[\Lambda]$ given by 
$x^\mu\otimes e^\lambda \mapsto x^{-w_\circ(\mu)}\otimes e^{-w_\circ(\lambda)}$
preserves the ideal $\mathcal{I}$ and, thus, induces
an involutive automorphism $\omega$ on $K_T(G/B)\simeq 
(\Z[X]\otimes\Z[\Lambda])/\mathcal{I}$.  Applying this involution to
the definition of the elementary Demazure operators $T_i$ 
in~(\ref{eq:isobaric-div-diff}), we deduce that
$\omega \, T_i \, \omega  = T_j$, where $j$ is given by
$\alpha_j = - w_\circ(\alpha_i)$, or equivalently, $s_j = w_\circ s_i w_\circ$.
Thus $\omega\, T_w \,\omega = T_{w_\circ w w_\circ}$, for any
$w\in W$.  Kostant-Kumar's formula~(\ref{eq:KK}) implies 
that $\omega:[\O_w]\mapsto [\O_{w_\circ w w_\circ}]$.
\end{proof}

\begin{proof}[Geometric proof
{\rm (due to Kumar~\cite{Kum})}]
Let $c\::\:G\rightarrow G$ be the Chevalley isomorphism. 
This is an algebraic
group isomorphism mapping $t\mapsto t^{-1}$ for $t$ in $T$, and $B\mapsto B^-$,
where $B^-$ is the {\em opposite Borel} subgroup. Also let
 $c_{w_\circ}\::\:G\rightarrow G$ be the automorphism given by $g\mapsto
\overline{w}_\circ g \overline{w}_\circ^{-1}$, where $\overline{w}_\circ$ in
$N(T)$ is a representative of $w_\circ$. Let $\phi\::\:G\rightarrow G$ be the
composite $c\circ c_{w_\circ}$. Then $\phi(B)=B$. Thus $\phi$ induces a variety
isomorphism $\overline{\phi}\::\:G/B\rightarrow G/B$. Moreover, since $c$
induces the identity map on the Weyl group, we see that
$\overline{\phi}(X_w)=X_{w_\circ w w_\circ}$. Thus $\overline{\phi}$ induces
the involution $\omega$ on $K_T(G/B)$ such that $\omega:[\O_w]\mapsto
[\O_{w_\circ w w_\circ}]$. 

To show that, under the isomorphism $\theta$, we have
$\omega:e^\lambda\mapsto e^{-w_\circ\lambda}$, 
we identify $G/B$ with $K/T$, where $K$
is a maximal compact subgroup of $G$. Let us consider the following bundle
morphism.  \[ \begin{diagram}
\node{K\times_T\C_{-w_\circ\lambda}}\arrow[2]{e,t}{\widehat{\phi}} \arrow{s}
\node[2]{K\times_T\C_\lambda} \arrow{s}\\
\node{K/T}\arrow[2]{e,b}{\overline{\phi}} \node[2]{K/T} \end{diagram} \] Here
we let $\widehat{\phi}(k,v_\circ):=(\phi(k),\overline{v}_\circ)$, where
$v_\circ$ is a generator of $\C_{-w_\circ\lambda}$, and $\overline{v}_\circ$ is
a generator of $\C_{\lambda}$. It is easy to see that $\widehat{\phi}$ is well
defined. Thus, we have 
$\omega\circ\theta(1\otimes e^\lambda)=\theta(1\otimes
e^{-w_\circ\lambda})\,$.
The proof of $\omega:x^\mu\mapsto x^{-w_\circ\mu}$ is
similar.  \end{proof}

Note that the map $\overline{\phi}$ in the above proof is not
$T$-equivariant, whence the involution $\omega$ is not a $\Z[X]$-linear map.

Let $c_{u,v}^w\in \Z[X]$ be the structure constants of $K_T(G/B)$ with 
respect to the basis of classes of structure sheaves of
Schubert varieties:
$$
[\O_u]\cdot[\O_v]=\sum_w c_{u,v}^w\,[\O_w]\,.
$$
The coefficients $c_u^w(\pm\omega_i)$ are related to certain structure 
constants $c_{u,v}^w$, as follows.

\begin{corollary} {\rm cf.~\cite{Bri}}
\label{strconst} 
For $v\ne w$, we have
\begin{itemize}
\item[(a)]
$c_u^w(-\omega_i) = - x^{-w_\circ(\omega_i)} \,
c_{w_\circ s_i,u}^w$\,;
\smallskip
\item[(b)]
$c_u^w(\omega_i)=(-1)^{\ell(u)-\ell(w)-1}
x^{\omega_i}\,
c_{s_iw_\circ,ww_\circ}^{uw_\circ}$\,;
\smallskip
\item[(c)]
$c_u^w(\omega_i)=
(-1)^{\ell(u)-\ell(w)-1} x^{\omega_i}\,
\iota(c_{w_\circ s_i,w_\circ w}^{w_\circ u})$\,.
\end{itemize}
Also, we have $c_{w_\circ s_i, u}^u = 1- x^{w_\circ(\omega_i) - u(\omega_i)}$.
\end{corollary}

The first two formulas (a) and (b) were given by Brion~\cite{Bri}
for $K(G/B)$ in a slightly  imprecise form.

\begin{proof}
Identity~(a) is obtained from the formula in Lemma~\ref{lem:O=L}(a)
by multiplying both sides by $[\O_u]$.  
Identity~(b) is obtained from~(a) and the duality
formula in Theorem~\ref{duality}, as follows:
\begin{align*}
c_u^w(\omega_i)&=
(-1)^{\ell(u)-\ell(w)}c_{ww_\circ}^{uw_\circ}(w_\circ(\omega_i))=
(-1)^{\ell(u)-\ell(w)}c_{ww_\circ}^{uw_\circ}(-\omega_j)\\&=
(-1)^{\ell(u)-\ell(w)-1} x^{-w_\circ(\omega_j)} \,
c_{w_\circ s_j,ww_\circ}^{uw_\circ}=
(-1)^{\ell(u)-\ell(w)-1} x^{\omega_i}\,
c_{s_iw_\circ,ww_\circ}^{uw_\circ}\,.
\end{align*}
Here we used the fact that $-w_\circ\alpha_i$ is the simple root $\alpha_j$ 
such that $s_j=w_\circ s_i w_\circ$.
Similarly, we obtain identity~(c) using the duality formula in
Theorem~\ref{th:second_duality}.
\end{proof}

\begin{remark}
We can easily expand the product $[\O_{w_\circ s_i}]\cdot[\O_{u}]$ using
our $K_T$-Chevalley formula, as shown in Corollary~\ref{cor:OiO=O}.
However, it is hard to apply the Pittie-Ram formula directly to the
calculation of this expansion, because the latter formula works for dominant 
weights only. 
In order to use this formula, one needs to invert the operator of 
multiplication by $e^{\omega_i}$ acting on the $|W|$-dimensional 
space $K_T(G/B)$.  Alternatively, one can use Brion's geometric argument
to derive the second formula in Corollary~\ref{strconst}.  But then, one
needs to apply the Pittie-Ram formula for computing {\it all\/} products
$e^{\omega_j}\cdot [\O_{ww_\circ}]$, for $w\in W$, and extract the
coefficient of $[\O_{uw_\circ}]$ in each result, where 
$j$ is given by $s_j=w_\circ s_i w_\circ$.
Indeed, we have no way of knowing in advance to which Weyl group element 
an LS-path leads, via Deodhar's lift operator.  In other words,
it is hard to ``invert'' the Pittie-Ram construction based on LS-paths
and Deodhar's lifts.
\end{remark}

\section{The Yang-Baxter Equation} 
\label{sec:yang-baxter}

Our construction is based on a certain 
$R$-matrix, that is, a collection of operators satisfying the Yang-Baxter equation.
In this section, we discuss the Yang-Baxter equation, following
the approach of Cherednik~\cite{Cher}. 
\medskip

For a pair of roots $\alpha,\beta\in\Phi$ 
such that $(\alpha,\beta)\leq 0$,
the subset of roots $\Delta\subset\Phi$ obtained from $\alpha$ and $\beta$ 
by a sequence of reflections $s_\alpha$ and $s_\beta$ is a rank 2 root 
system of type $A_1\times A_1$, $A_2$, $B_2$, or $G_2$.
The reflections $s_\alpha$ and $s_\beta$ generate a dihedral subgroup in $W$
of order $2m$, where $m=2,3,4,6$, for types 
$A_1\times A_1$, $A_2$, $B_2$, $G_2$, respectively.
The condition $(\alpha,\beta)\leq 0$ implies that $\alpha,\beta$ 
form a system of simple roots for $\Delta$.
The $m$ roots in $\Delta$ expressible as nonnegative linear combinations
of $\alpha$ and $\beta$ can be normally ordered as follows:
$\alpha, s_\alpha(\beta),s_\alpha s_\beta(\alpha),
\dots, s_\beta(\alpha),\beta$.

The following definition was given by Cherednik~\cite[Definition~2.1a]{Cher}
in a slightly different form.

\begin{definition}
We say that a collection of invertible operators $\{R_\alpha 
\mid \alpha\in\Phi\}$ 
labeled by roots
satisfies the {\it Yang-Baxter equation\/} if
$R_{-\alpha}=(R_\alpha)^{-1}$ and, for any pair of roots
$\alpha,\beta\in\Phi$ such that $(\alpha,\beta)\leq 0$, we have
\begin{equation}
R_{\alpha} R_{s_\alpha(\beta)} 
R_{s_\alpha s_\beta(\alpha)}\cdots
R_{s_\beta(\alpha)} R_\beta
=
R_\beta R_{s_\beta(\alpha)} \cdots
R_{s_\alpha s_\beta(\alpha)} R_{s_\alpha(\beta)} R_{\alpha}.
\label{eq:YB}
\end{equation}
A collection of operators $\{R_\alpha \mid \alpha\in\Phi\}$ satisfying
the Yang-Baxter equation is also called an {\it $R$-matrix}.
\end{definition}

For example, the operators $R_\alpha$ and $R_\beta$ 
commute whenever 
$(\alpha,\beta)=0$.  If $\Delta$ is of type $A_2$, then 
the Yang-Baxter equation~(\ref{eq:YB}) says that
$$
R_{\alpha} R_{\alpha+\beta} R_\beta =
R_{\beta} R_{\alpha+\beta} R_\alpha.
$$

The following two lemmas are implicit in~\cite{Cher}.

\begin{lemma} Consider a collection $\{R_\alpha \mid \alpha\in\Phi^+\}$ 
of invertible operators labeled by positive roots which satisfies 
the Yang-Baxter equation~{\rm(\ref{eq:YB})}, for any pair of positive roots 
$\alpha,\beta\in\Phi^+$ such that $(\alpha,\beta)\leq 0$.
Let us extend this collection to all roots $\alpha\in\Phi$ by
$R_{-\alpha}:=(R_\alpha)^{-1}$.
Then the collection $\{R_\alpha \mid \alpha\in\Phi\}$ is an $R$-matrix.
\label{lem:YB-positive-roots}
\end{lemma}

\begin{proof}
Let us multiply the Yang-Baxter equation~(\ref{eq:YB}) by $R_{-\beta}$ 
on the left and on the right. We get
$$
R_{-\beta}R_{\alpha} R_{s_\alpha(\beta)} R_{s_\alpha s_\beta(\alpha)}\cdots
R_{s_\beta(\alpha)} 
=
R_{s_\beta(\alpha)} \cdots
R_{s_\alpha s_\beta(\alpha)}
R_{s_\alpha(\beta)} R_{\alpha} R_{-\beta}.
$$
This is the same equation with $(\alpha,\beta)$ replaced
by the pair $(s_\beta(\beta), s_\beta(\alpha))$.
Applying this procedure repeatedly,
we can always transform the pair $(\alpha,\beta)$ into 
a pair of positive roots.
\end{proof}

For a decomposition $v=s_{i_1}\cdots s_{i_l}\in\Waff$, reduced
or not, of an affine Weyl group element $v$, let
$(\beta_1,\dots,\beta_l)$ be the corresponding $\lambda$-chain of roots.
For an $R$-matrix $\{R_\alpha \mid \alpha\in\Phi\}$, let us
define $R^{(s_{i_1}\cdots s_{i_l})} = R_{\beta_l} R_{\beta_{l-1}}\cdots 
R_{\beta_2} R_{\beta_1}$.

\begin{lemma}  
Let $\{R_\alpha \mid \alpha\in\Phi\}$ be an $R$-matrix.
Then the operator $R^{(s_{i_1}\cdots s_{i_l})}$ depends only on
the affine Weyl group element $v=s_{i_1}\cdots s_{i_l}$, 
not on the choice of the decomposition.
\label{lem:depends-only-on}
\end{lemma}

\begin{proof}
The Coxeter relations~(\ref{eq:coxeter-relations-affine}) imply that
any two decompositions of $v$ can be related by a sequence of local moves of the following two types:
(1) adding or removing segments $s_i s_i$; (2)
the Coxeter moves
\begin{equation}
s_{i_1} \cdots s_{i_a}
\stackrel{\textrm{\tiny $m_{ij}$ terms}}
{(s_i s_j s_i \cdots)}  
s_{i_b}\cdots s_{i_l}
\quad\longrightarrow\quad
s_{i_1} \cdots s_{i_a}
\stackrel{\textrm{\tiny $m_{ij}$ terms}}
{(s_j s_i s_j \cdots)} 
s_{i_b}\cdots s_{i_l}.
\label{eq:Coxeter-move}
\end{equation}

Adding or removing a segment $s_i s_i$ in a decomposition for $v$
results in adding or removing a segment $\beta,-\beta$ in 
the sequence of roots $(\beta_1,\dots,\beta_l)$.
This does not change the operator
$R_{\beta_l}\cdots R_{\beta_1}$, because $R_\beta R_{-\beta} = 1$.
A Coxeter move~(\ref{eq:Coxeter-move}) results in 
applying the Yang-Baxter transformation
$$
\alpha, s_\alpha(\beta), \dots ,s_\beta(\alpha),\beta
\quad\longrightarrow\quad
\beta, s_\beta(\alpha),\dots,s_\alpha(\beta),\alpha
$$
to the segment $(\beta_{a+1},\dots,\beta_{b-1})
=(\alpha,s_\alpha(\beta),\cdots,\beta)$
in the sequence $(\beta_1,\dots,\beta_l)$.
Here we have $\alpha=\bar s_{i_1}\cdots \bar s_{i_a}(\alpha_i)$ and
$\beta=\bar s_{i_1}\cdots \bar s_{i_a}(\alpha_j)$.
Note that $(\alpha,\beta) = (\alpha_i,\alpha_j)\leq 0$.
The Yang-Baxter equation~(\ref{eq:YB}) guarantees that 
this transformation of the sequence $(\beta_1,\dots,\beta_l)$ does 
not change the operator $R_{\beta_l}\cdots R_{\beta_1}$.
\end{proof}

\section{Bruhat Operators}
\label{sec:bruhat-operators}

In this section, we present a class of solutions of the Yang-Baxter 
equation.
\medskip

It will be convenient to extend the ring of coefficients $\Z[X]=R(T)$ 
in $K_T(G/B)$ as follows.
Let us shrink the weight lattice $h$ times by defining $\Lambda/h:= 
\{\lambda/h \mid \lambda \in\Lambda\}$,
where $h := (\rho,\theta^\vee)+1$ is the Coxeter number.
Let $\Z[\tilde X]$
be the group algebra of $\Lambda/h$, which has formal exponents 
$x^{\lambda/h}$, for $\lambda\in\Lambda$.
This is the algebra of Laurent polynomials  
$\Z[\tilde X]=\Z[x^{\pm\omega_1/h},\dots,x^{\pm\omega_r/h}]$.
Let 
$$
\tilde K_T(G/B) := K_T(G/B)\otimes_{\Z[X]} \Z[\tilde X].
$$
The space $\tilde K_T(G/B)$ has 
the $\Z[\tilde X]$-linear basis given by the classes $[\O_w]$, for $w\in W$.

For a positive root $\alpha\in\Phi^+$, let us define the 
{\it Bruhat operator\/} $B_{\alpha}$ acting $\Z[\tilde X]$-linearly on
$\tilde K_T(G/B)$ by
\begin{equation}
B_\alpha : [\O_w]\longmapsto\left\{\begin{array}{cl}
[\O_{ws_\alpha}] & \textrm{if }  \ell(ws_\alpha) = \ell(w)-1, \\[.05in]
0 & \textrm{otherwise.}
\end{array}\right.
\label{eq:Bruhat-operators}
\end{equation}
Also define $B_{\alpha} := - B_{-\alpha}$, if $\alpha$ is a negative root.  The operators $B_\alpha$ 
move Weyl group elements one step down in Bruhat order.   

For a weight $\lambda$, 
define the $\Z[\tilde X]$-linear operators $X^\lambda$ acting on
$\tilde K_T(G/B)$ by
\begin{equation}
X^\lambda: [\O_w] \mapsto x^{w(\lambda/h)} [\O_w].
\label{eq:Qlambda}
\end{equation}

For $\alpha\in\Phi$ and $\lambda,\mu\in \Lambda$,
these operators satisfy the following relations:
\begin{eqnarray}
&&(B_\alpha)^2=0\,,
\label{eq:MM=0}
\\[.05in]
&& X^\lambda \,X^\mu = X^{\lambda+\mu}\,,\\[.05in]
&&B_\alpha \, X^\lambda = X^{s_\alpha(\lambda)} \, B_\alpha\,.
\label{eq:MQ=QM}
\end{eqnarray}

For a fixed weight $\lambda$ and 
$k\in\Z$, we define a family 
of operators $\{R_{\alpha} \mid \alpha\in\Phi\}$
labeled by roots $\alpha\in\Phi$ acting on $\tilde K_T(G/B)$ as follows:
\begin{equation}
R_\alpha =  X^{k\alpha} +  
X^{(\lambda,\alpha^\vee)\,\alpha}\,B_\alpha=
X^\lambda\,(X^{k\alpha} +  B_\alpha)\,X^{-\lambda}.
\label{eq:R-alpha-lambda}
\end{equation}

Using relations~(\ref{eq:MM=0}) and~(\ref{eq:MQ=QM}),
we obtain
$$
R_{-\alpha} =  X^{-k\alpha} - 
X^{(\lambda,\alpha^\vee) \,\alpha}\, 
B_{\alpha}
= (R_{\alpha})^{-1}.
$$

\begin{theorem} 
Fix a weight $\lambda$ and $k\in\Z$.
The family of operators $\{R_\alpha \mid \alpha\in\Phi\}$ 
given by~{\rm(\ref{eq:R-alpha-lambda})} 
satisfies the Yang-Baxter equation~{\rm(\ref{eq:YB})}.
\label{th:YB}
\end{theorem}

\begin{proof}
Let us first assume that $\lambda=0$ and $k=0$. 
In this case $R_\alpha = 1+B_\alpha$.
In~\cite{BFP}, we proved the Yang-Baxter equation for a general class
of operators by checking it for all the rank $2$ root systems (that is, for types $A_1\times A_1$, $A_2$, $B_2$, and $G_2$).  
In particular, the results of~\cite{BFP} imply that the family of operators 
$\{1+B_\alpha \mid \alpha\in\Phi^+\}$ satisfies the Yang-Baxter 
equation~(\ref{eq:YB}).  Also $R_{-\alpha} = 1-B_\alpha = 
(1+B_\alpha)^{-1}= (R_\alpha)^{-1}$.
According to Lemma~\ref{lem:YB-positive-roots}, the collection
$\{1+B_\alpha \mid \alpha\in\Phi\}$ is an $R$-matrix.

Let us now consider the general case.
For $\alpha\in\Phi$ and $n\in\Z$, let us define
$$
\hat R_\alpha^n := 1 + X^{n\alpha}\, B_\alpha.
$$
Then 
$ R_\alpha = X^{k\alpha} \, 
\hat R_\alpha^{(\lambda,\alpha^\vee)-k}$.
For $\mu\in\Lambda$, we get, using~(\ref{eq:MQ=QM}), 
\begin{equation}
\hat R_\alpha^n \,X^\mu = 
X^\mu\, \hat R_\alpha^{n - (\mu,\alpha^\vee)}.
\label{eq:hat-RQ=QR}
\end{equation}

Let us write the left-hand side of the Yang-Baxter equation~(\ref{eq:YB}) 
as follows:
$$
R_{\gamma_1} \cdots  R_{\gamma_m}
= 
X^{k\gamma_1} \,\hat R_{\gamma_1}^{n_1} \,
X^{k\gamma_2} \, \hat R_{\gamma_2}^{n_2} \cdots
X^{k\gamma_m} \, \hat R_{\gamma_m}^{n_l},
$$
where $(\gamma_1,\dots,\gamma_m) = 
(\alpha,s_\alpha(\beta),\cdots,s_\beta(\alpha),\beta)$
and $n_i = (\lambda,\gamma_i^\vee)-k$.
Using (\ref{eq:hat-RQ=QR}) to commute all $X^{k\gamma_i}$
to the left, we obtain the expression 
$$
R_{\gamma_1} \cdots R_{\gamma_m} = 
X^{k(\gamma_1+\cdots+\gamma_m)} \,
\hat R_{\gamma_1}^{n_1'} \, \hat R_{\gamma_2}^{n_2'} \cdots
\hat R_{\gamma_m}^{n_l'},
$$
where
$$
n_i'=n_i-\sum_{j=i+1}^m k(\gamma_{j},\gamma_i^\vee)=
(\lambda-k(\gamma_{i+1}-\cdots-\gamma_m),\gamma_i^\vee)-k.
$$

Let us show that
$$(\gamma_1+\cdots+\gamma_{i-1},\gamma_i^\vee) =
(\gamma_{i+1}+\cdots+\gamma_{m},\gamma_i^\vee)\,,$$ 
for all $i=1,\dots,m$.  Suppose that $i\leq (m+1)/2$.
The reflection $s_{\gamma_i}$ sends the roots 
$\gamma_1,\dots,\gamma_{i-1}$ to 
$-\gamma_{2i-1},\dots,-\gamma_{i+1}$, 
and the roots $\gamma_{2i},\dots,\gamma_m$ to 
$\gamma_m,\dots,\gamma_{2i}$, respectively.
Thus 
$$(\gamma_1+\cdots+\gamma_{i-1},\gamma_i^\vee) =
(\gamma_{i+1}+\cdots+\gamma_{2i-1},\gamma_i^\vee)
\quad\textrm{and}\quad
(\gamma_{2i}+\cdots+\gamma_m,\gamma_i^\vee)=0,
$$
as needed.  Since $(\gamma_i,\gamma_i^\vee)=2$, we get 
$$
n_i' = (\lambda-k(\gamma_{i+1}+\cdots+\gamma_m),\gamma_i^\vee)-k =
(\lambda-k\varrho,\gamma_i^\vee),
$$
where $\varrho = \frac{1}{2}(\gamma_1+\cdots+\gamma_m)$ 
is the ``rho'' for the rank 2 root system $\Delta$ generated 
by $\alpha$ and $\beta$.

This shows that 
$$
R_{\gamma_1} \cdots R_{\gamma_m} = 
X^{2k\varrho} \hat R_{\gamma_1}^{(\mu,\gamma_1^\vee)}\cdots 
R_{\gamma_l}^{(\mu,\gamma_m^\vee)}
=X^{\mu+2k\varrho} \hat R_{\gamma_1}^0 \cdots \hat R_{\gamma_m}^0
X^{-\mu},
$$
where $\mu=\lambda-k\varrho$.
Analogously, the right-hand side of the Yang-Baxter equation~(\ref{eq:YB}) 
can be written as
$$
R_{\gamma_m} \cdots R_{\gamma_1} = 
X^{\mu+2k\varrho} \hat R_{\gamma_m}^0 \cdots \hat R_{\gamma_1}^0
X^{-\mu}\,.
$$
The fact that the operators $\hat R_\alpha^0 = 1+B_\alpha$
satisfy the Yang-Baxter equation implies that
the family $\{R_\alpha \mid \alpha\in\Phi\}$
satisfies the Yang-Baxter equation as well.  This concludes the proof.
\end{proof}

In the rest of the paper, we only use a special case of 
the operators $R_\alpha$ defined in ~(\ref{eq:R-alpha-lambda}), namely we set $\lambda:=\rho$ and $k:=1$, which leads to
\begin{equation}
R_\alpha = X^\alpha + X^{(\rho,\alpha^\vee)\,\alpha}\,B_\alpha
= X^\rho\, (X^\alpha + B_\alpha)\,X^{-\rho},	
\quad\textrm{for } \alpha\in\Phi.
\label{eq:R-alpha=Q+QM}
\end{equation}

\section{Commutation Relations}
\label{sec:commutation-relations}

Let $T_i$ be the operator on $\tilde K_T(G/B)$
induced by the elementary Demazure operator~(\ref{eq:isobaric-div-diff}),
for $i=1,\dots,r$.  In view of~(\ref{eq:T-reduce}) and~(\ref{eq:KK}), 
this operator acts $\Z[\tilde X]$-linearly on $\tilde K_T(G/B)$ as
$$
T_i : [\O_w]\longmapsto\left\{\begin{array}{cl}
[\O_{ws_i}] & \textrm{if }  \ell(ws_i) = \ell(w)+1, \\[.05in]
[\O_{w}] & \textrm{if }  \ell(ws_i) = \ell(w)-1.
\end{array}\right.
$$
Let $B_i := B_{\alpha_i}$ be the Bruhat operator for a simple reflection, which is the $\Z[\tilde X]$-linear operator on $\tilde K_T(G/B)$ defined by
$$\begin{array}{l}
B_i:[\O_w]\mapsto
\left\{\begin{array}{cl}
[\O_{ws_i}] &\textrm{if } \ell(ws_i) = \ell(w) - 1,\\[.05in]
0 &\textrm{if } \ell(ws_i) = \ell(w) + 1.
\end{array}\right.
\end{array}$$
Let us define a similar $\Z[\tilde X]$-linear operator $B_i^*$ by
$$\begin{array}{l}
B_i^*:[\O_w]\mapsto
\left\{\begin{array}{cl}
[\O_{ws_i}] &\textrm{if } \ell(ws_i) = \ell(w) + 1,\\[.05in]
0 &\textrm{if } \ell(ws_i) = \ell(w) - 1.
\end{array}\right.
\end{array}
$$
Since both operators $B_i^*$ and $B_i$ map $[\O_w]$ to 
$[\O_{ws_i}]$ or to zero, we have
\begin{equation}
X^\mu \,B_i^* = B_i^* \,X^{s_i(\mu)},
\quad\textrm{and}\quad
X^\mu \,B_i = B_i \,X^{s_i(\mu)},
\label{eq:QU=UQ}
\end{equation}
for any weight $\mu\in\Lambda$.

The operator $B_i^*$ can be expressed in terms of 
$T_i$ and $B_i$ as follows.

\begin{lemma} We have
$B_i^*= T_i\,(1-B_i)= (1+B_i)(T_i-1)$, for $i=1,\dots,r$.
\label{lem:UTTD}
\end{lemma}

\begin{proof}
It is enough to check this claim for restrictions of
the operators on the 2-dimensional invariant subspace spanned by 
$[\O_w]$ and $[\O_{ws_i}]$, for any $w\in W$ such that $\ell(ws_i)=\ell(w)+1$.
The required identity is 
$$
\begin{pmatrix} 0 & 0 \\ 1 & 0 \end{pmatrix}=
\begin{pmatrix} 0 & 0 \\ 1 & 1 \end{pmatrix} \,
\begin{pmatrix} 1 & -1 \\ 0 & 1 \end{pmatrix}  =
\begin{pmatrix} 1 & 1 \\ 0 & 1 \end{pmatrix}  \,
\begin{pmatrix} -1 & 0 \\ 1 & 0 \end{pmatrix},
$$  
which we leave to the reader as an exercise.
\end{proof}

Recall that $B_\beta$ are the Bruhat operators given 
by~(\ref{eq:Bruhat-operators}).

\begin{lemma}  {\rm cf.\ Deodhar~\cite[Lemma~2.1]{Deo1}  } 
We have $B_\beta \,B_i^* = B_i^*\, B_{s_i(\beta)}$,
for $i=1,\dots,r$ and $\beta\in\Phi$ such that $\beta\ne\pm\alpha_i$. 
\label{lem:MU=UM}
\end{lemma}

\begin{proof}
We may assume that $\beta\in\Phi^+$.  Let $\beta'=s_i(\beta)$.  
Then $\beta'\in\Phi^+$ and $\beta'\ne\alpha_i$.  
Both operators $B_\beta \,B_i^*$ and 
$B_i^*\, B_{\beta'}$ map $[\O_w]$ to
$[\O_{ws_i s_\beta}]= [\O_{w s_{\beta'} s_i}]$ or to zero.
Thus, we need to show that $B_\beta \,B_i^*([\O_w])$ is nonzero 
if and only if $B_i^*\, B_{\beta'}([\O_w])$ is nonzero.

Suppose that this is not true. 
One possibility is that we have
$B_\beta \,B_i^*([\O_w])=0$ and $B_i^*\, B_{\beta'}([\O_w])\ne 0$.
Then $\ell(w)=\ell(ws_{\beta'})+1 = \ell(w s_i) + 1 = \ell(ws_{\beta'} s_i)$.
Indeed, $B_i^*\, B_{\beta'}([\O_w])\ne 0$ implies that $\ell(w s_{\beta'}) = 
\ell(w)-1$ and $\ell(w s_{\beta'} s_i) = \ell(w s_{\beta'}) + 1$,
while $B_\beta \,B_i^*([\O_w])=0$ implies that $\ell(w s_i)\ne \ell(w)+1$,
and, thus, $\ell(w s_i)= \ell(w)-1$.

Let us choose a reduced decomposition for $w=s_{i_1}\cdots s_{i_l}$
such that $i_l=i$.  By the Strong Exchange Condition \cite[Theorem~5.8]{Hum}, the fact that $\ell(w)=\ell(ws_{\beta'})+1$ implies
that there exists $k\in\{1,\dots,l\}$ such that 
$s_{i_1}\cdots \widehat{s_{i_k}} \cdots s_{i_l}$ is a reduced
decomposition for $ws_{\beta'}$. 
Furthermore, we have $\beta'=s_{i_l}\cdots s_{i_{k+1}}(\alpha_{i_k})$.
Since $\beta'\ne \alpha_i$, we have $k\ne l$. 
We obtain a reduced decomposition for 
$ws_{\beta'}$ that ends with $s_i$.   Thus $\ell(w s_{\beta'} s_i) = 
\ell(w s_{\beta'}) - 1$, 
which is a contradiction.

Now suppose that we have $B_\beta \,B_i^*([\O_w])\ne 0$ and 
$B_i^*\, B_{\beta'}([\O_{w}])= 0$.
Then $\ell(w)=\ell(ws_i)-1= \ell(ws_{\beta'})-1 = \ell(ws_{\beta'} s_i)$
or, equivalently, 
$\ell(w') = \ell(w's_i) + 1 = \ell(w's_{\beta})+1 = \ell(w' s_\beta s_i)$,
for $w'=ws_i$.  The above argument shows that this is impossible.
\end{proof}

\begin{remark}\label{deodhar} The contradictions derived in the above proof are
essentially the content of Lemma~2.1 in~\cite{Deo1}, 
which is proved in a similar way.
\end{remark}

Let $\{R_\alpha \mid\ \alpha\in\Phi\}$ be the $R$-matrix 
given by~(\ref{eq:R-alpha=Q+QM}).
The main technical result of this section is the following statement that
gives a commutation relation between this $R$-matrix and 
the Demazure operators $T_i$.

\begin{proposition}  For any $\beta\in\Phi$ and $i=1,\dots,r$, we have
\begin{enumerate}
\item[(a)] 
$R_{\alpha_i} \, T_i = T_i \,R_{-\alpha_i} + R_{\alpha_i}$,
\item[(b)] 
$R_{-\alpha_i} \,T_i = T_i \,R_{\alpha_i} - R_{\alpha_i}$,
\item[(c)]
$R_\beta \, T_i = T_i \,R_{-\alpha_i} \,R_{s_i(\beta)} \,R_{\alpha_i}$
if $\beta\ne \pm\alpha_i$.
\end{enumerate}
\label{prop:RT=TR}
\end{proposition}

\begin{proof}
We have $R_{\alpha_i}=X^{\alpha_i}\,(1+B_i)$ and 
$R_{-\alpha_i}=(1-B_i)\,X^{-\alpha_i}$.

(a)  \ 
By Lemma~\ref{lem:UTTD},  
$(1+B_i)\,(T_i-1) = T_i\,(1-B_i)$. 
Thus  
$$
X^{\alpha_i}\,(1+B_i)\,T_i = X^{\alpha_i}\,T_i\,(1-B_i) + X^{\alpha_i}\,(1+B_i).
$$
Then use~(\ref{eq:QU=UQ}) to commute $X^{\alpha_i}$ 
with $T_i\,(1-B_i) = B_i^*$ in the first term in the 
right-hand side.  This produces~(a).

(b) \  Multiply (a) by $R_{-\alpha_i}$ on the left and by 
$R_{\alpha_i}$ on the right.

(c) \  Let $\beta'=s_i(\beta)$. 
Identity~(c) can be written as
$$
(X^\beta + X^{k\beta}\,B_\beta)\,T_i = T_i \, (1-B_i)\,X^{-\alpha_i}\,
(X^{\beta'}+X^{k'\beta'}\,B_{\beta'})\,
X^{\alpha_i}\,(1+B_i),
$$
where $k=(\rho,\beta^\vee)$ and
$k'=(\rho,(\beta')^\vee)=(s_i(\rho),\beta^\vee)=
(\rho-\alpha_i,\beta^\vee)$.
The right-hand side of this identity can be written as
$$ 
T_i\,(1-B_i)\,
(X^{\beta'}+X^{k\beta'}\,B_{\beta'})\, (1+B_i).
$$
Indeed, $ X^{k'\beta'-\alpha_i}\,B_{\beta'} X^{\alpha_i} = 
X^{k\beta'} B_{\beta'}$, because
$k'\beta'-\alpha_i+ s_{\beta'}(\alpha_i)
=(\rho-\alpha_i,\beta^\vee)\,\beta'
-(\alpha_i,(\beta')^\vee)\,\beta'=
(\rho,\beta^\vee)\,\beta' = k \beta'$.
Commuting $X^{\beta'}$ and $X^{k\beta'}\,B_{\beta'}$ 
with $T_i\,(1-B_i)=B_i^*$ using~(\ref{eq:QU=UQ}) 
and Lemma~\ref{lem:MU=UM}, we can rewrite this as
$$
(X^{\beta}\,+X^{k\beta}\,B_{\beta})\,B_i^*(1+B_i) =
(X^{\beta}\,+X^{k\beta}\,B_{\beta})\,T_i,
$$
which is equal to the left-hand side of required identity.
\end{proof}

\section{Path Operators}
\label{sec:path-operators}

Recall that $v_{-\lambda}\in\Waff$, $\lambda\in\Lambda$, is the unique element 
of the affine Weyl group such that 
$v_{-\lambda}(A_\circ) = A_{-\lambda}=A_\circ-\lambda$. 
Each decomposition $v_{-\lambda} = s_{i_1}\cdots s_{i_l}$ in $\Waff$ 
corresponds to an alcove path 
$A_\circ\stackrel{-\beta_1}\longrightarrow \cdots
\stackrel{-\beta_l}\longrightarrow A_{-\lambda}$;
and the sequence of roots $(\beta_1,\dots,\beta_l)$ is called 
a $\lambda$-chain, see Definition~\ref{def:lambda-chain}.
Also recall that there is an associated alcove path
$A_\circ\stackrel{\beta_l}\longrightarrow \cdots
\stackrel{\beta_1}\longrightarrow A_{\lambda}$, as discussed in
Remark~\ref{remtransl}.

For $\lambda\in\Lambda$, let us define the operator
$R^{[\lambda]}$ acting on 
$\tilde K_T(G/B)$ by
\begin{equation}
R^{[\lambda]} := R_{\beta_l} R_{\beta_{l-1}}\cdots R_{\beta_2} R_{\beta_1},
\label{eq:R[lambda]}
\end{equation}
where $(\beta_1,\dots,\beta_l)$ is a $\lambda$-chain,
and the $R$-matrix $\{R_\alpha \mid \alpha\in\Phi\}$ 
is given by~(\ref{eq:R-alpha=Q+QM}).

\begin{remark}\label{remybe}
Theorem~\ref{th:YB} and Lemma~\ref{lem:depends-only-on}
imply that the operator $R^{[\lambda]}$
depends only on the weight $\lambda$ and does not depend on the 
choice of a $\lambda$-chain.

\end{remark}

The following result is not used in subsequent proofs. We state it because it exhibits the commutativity of the operators $E^\lambda$ and $E^\mu$ in our combinatorial model, based on Remark~\ref{remybe}.

\begin{proposition}
For any $\lambda,\mu\in\Lambda$, we have
$R^{[\lambda]}\cdot R^{[\mu]} = R^{[\lambda+\mu]}$.
\end{proposition}

\begin{proof}
Let us choose a  $\lambda$-chain $(\beta_1,\dots,\beta_l)$
and a $\mu$-chain $(\beta_1',\dots,\beta_m')$.
They  correspond
to alcove paths
$A_\circ\stackrel{\beta_l}\longrightarrow \cdots
\stackrel{\beta_1}\longrightarrow A_{\lambda}$
and $A_\circ\stackrel{\beta_m'}\longrightarrow \cdots
\stackrel{\beta_1'}\longrightarrow A_{\mu}$.
If we translate all alcoves in the second path $\lambda$, we obtain
the alcove path $A_{\lambda}\stackrel{\beta_m'}\longrightarrow \cdots
\stackrel{\beta_1'}\longrightarrow A_{\lambda+\mu}$.
Let us concatenate the first path from $A_\circ$ to $A_{\lambda}$
with the translated path from  $A_{\lambda}$ to $A_{\lambda+\mu}$.
We obtain the alcove path
$$
A_\circ\stackrel{\beta_l}\longrightarrow \cdots
\stackrel{\beta_1}\longrightarrow A_{\lambda}
\stackrel{\beta_m'}\longrightarrow \cdots
\stackrel{\beta_1'}\longrightarrow A_{\lambda+\mu}.
$$
This shows that the sequence 
$(\beta_1',\dots,\beta_m',\beta_1,\dots,\beta_l)$ is a
$(\lambda+\mu)$-chain.
Thus
$$
R^{[\lambda]}\cdot R^{[\mu]}= R_{\beta_l}\cdots R_{\beta_1}
R_{\beta_m'}\cdots R_{\beta_1'} = R^{[\lambda+\mu]}, 
$$
as needed.
\end{proof}

\begin{lemma}
Let $(\beta_1,\dots,\beta_l)$  be a $\lambda$-chain.  Then,
for any $i=1,\dots,r$, the sequence of roots \linebreak 
$(\alpha_i,s_i(\beta_1),\dots, s_i(\beta_l),-\alpha_i)$
is an $s_i(\lambda)$-chain.
\label{lem:s_i(beta)}
\end{lemma}

\begin{proof}  
Applying the reflection $s_i$ to the alcove path
$A_\circ\stackrel{\beta_l}\longrightarrow \cdots
\stackrel{\beta_1}\longrightarrow A_{\lambda}$, 
we obtain the alcove path
$s_i(A_\circ)\stackrel{s_i(\beta_l)}\longrightarrow \cdots
\stackrel{s_i(\beta_1)}\longrightarrow s_i(A_{\lambda})$. 
We have $A_\circ\stackrel{-\alpha_i}\longrightarrow s_i(A_\circ)$.
Translating this relation by $s_i(\lambda)$, we obtain
$(s_i(A_\circ) + s_i(\lambda)) 
\stackrel{\alpha_i}\longrightarrow
(A_\circ + s_i(\lambda))$, or, equivalently,
$s_i(A_{\lambda})\stackrel{\alpha_i}\longrightarrow A_{s_i(\lambda)}$.
Thus
$$
A_\circ\stackrel{-\alpha_i}\longrightarrow s_i(A_\circ)
\stackrel{s_i(\beta_l)}\longrightarrow \cdots
\stackrel{s_i(\beta_1)}\longrightarrow s_i(A_{\lambda}) 
\stackrel{\alpha_i}\longrightarrow A_{s_i(\lambda)}
$$
is an alcove path, and $(\alpha_i,s_i(\beta_1),\dots,
s_i(\beta_l),-\alpha_i)$ is an $s_i(\lambda)$-chain.
\end{proof}

\begin{lemma}
Let $(\beta_1,\dots,\beta_l)$ be a $\lambda$-chain, and let 
$A_0\stackrel{\beta_l}\longrightarrow \cdots \stackrel{\beta_1}
\longrightarrow A_{l}$ be the corresponding 
alcove path from $A_0=A_\circ$ to $A_l=A_\lambda$.
Assume that $\pm\beta_j = \alpha_i$ is a simple root, for some 
$i\in\{1,\dots,r\}$ and $j\in\{1,\dots,l\}$.  
Then
$$
(\alpha_i,s_i(\beta_1),\dots, s_i(\beta_{j-1}), \beta_{j+1},\dots,\beta_l)
$$
is an $s(\lambda)$-chain, where
$s=s_{\alpha_i,k}$ denotes the affine reflection with respect 
to the common wall of the alcoves 
$A_{l-j}\stackrel{\beta_j}\longrightarrow A_{l-j+1}$.
\label{lem:flip-tail}
\end{lemma}

\begin{proof}  
Let us apply the following tail-flip to the alcove path 
$A_0\stackrel{\beta_l}\longrightarrow \cdots \stackrel{\beta_1} 
\longrightarrow A_{l}$.
We leave the initial segment 
$A_0\stackrel{\beta_l}\longrightarrow\cdots
\stackrel{\beta_{j+1}}\longrightarrow A_{l-j}$ unmodified
and apply the affine reflection $s$
to the remaining tail:
$s(A_{l-j+1}) \stackrel{\bar s(\beta_{j-1})} \longrightarrow s(A_{l-j+2})
\stackrel{\bar s(\beta_{j-2})} \longrightarrow \cdots
\stackrel{\bar s(\beta_{1})} \longrightarrow s(A_l)$.
Note that $A_{l-j} = s(A_{l-j+1})$ and $\bar s = s_i$. 
Also note that $s(A_l) = s(A_\circ+\lambda) = s_i(A_\circ) + s(\lambda)$,
and, thus, $s(A_l)\stackrel{\alpha_i} \longrightarrow A_{s_i(\lambda)}$.
Let us add the step 
$s(A_l)\stackrel{\alpha_i} \longrightarrow A_{s_i(\lambda)}$ at the end
of the alcove path with flipped tail.
We obtain the alcove path
$$
A_0\stackrel{\beta_l}\longrightarrow 
\cdots \stackrel{\beta_{j+1}} \longrightarrow A_{l-j}
\stackrel{s_i(\beta_{j-1})} \longrightarrow s(A_{l-j+2})
\stackrel{s_i(\beta_{j-2})} \longrightarrow \cdots
\stackrel{s_i(\beta_{1})} \longrightarrow s(A_l)
\stackrel{\alpha_i} \longrightarrow A_{s_i(\lambda)}.
$$
from $A_\circ$ to $A_{s_i(\lambda)}$.   Thus
$(\alpha_i,s_i(\beta_1),\dots,s_i(\beta_{j-1}), \beta_{j+1},\dots,\beta_l)$
is an $s(\lambda)$-chain.
\end{proof}

\begin{proposition}
For any $\lambda\in\Lambda$ and $i\in\{1,\dots,r\}$, we have
$$
R^{[\lambda]}\cdot T_i = T_i \cdot R^{[s_i(\lambda)]}  + 
\sum_{0\leq k < (\lambda,\alpha_i^\vee)}  R^{[\lambda - k \alpha_i]} 
-\sum_{(\lambda,\alpha_i^\vee)\leq k < 0}  R^{[\lambda - k \alpha_i]}. 
$$
\label{prop:same-relations}
\end{proposition}

\begin{proof}
Let us choose a $\lambda$-chain $(\beta_1,\dots,\beta_l)$.
Let $A_0\stackrel{\beta_l}\longrightarrow \cdots \stackrel{\beta_1}
\longrightarrow A_{l}$ be the corresponding 
alcove path from $A_0=A_\circ$ to $A_l=A_\lambda$.
And let $r_j$ be the affine reflection with respect to the common wall 
of the alcoves $A_{l-j}\stackrel{\beta_j}\longrightarrow A_{l-j+1}$.

Then $R^{[\lambda]}= R_{\beta_l}\cdots R_{\beta_1}$.  
Using the relations in Proposition~\ref{prop:RT=TR} repeatedly
to commute $T_i$ with $R_{\beta_l}\cdots R_{\beta_1}$,  
we obtain
$$
\begin{array}{rcl}
R_{\beta_l}\cdots R_{\beta_1} \,T_i &=& 
\ds
T_i\, R_{-\alpha_i} R_{s_i(\beta_l)} \cdots
R_{s_i(\beta_1)} R_{\alpha_i}  \\[.2in]
&+&
\ds
\sum_{j:\,\beta_j=\alpha_i}
R_{\beta_l}\cdots R_{\beta_{j+1}} R_{s_i(\beta_{j-1})} \cdots
R_{s_i(\beta_{1})} R_{\alpha_i} \\[.2in]
&-&
\ds
\sum_{j:\,\beta_j=-\alpha_i}
R_{\beta_l}\cdots R_{\beta_{j+1}} R_{s_i(\beta_{j-1})} \cdots
R_{s_i(\beta_{1})} R_{\alpha_i}.
\end{array}
$$
According to Lemmas~\ref{lem:s_i(beta)} and~\ref{lem:flip-tail},
the right-hand side of this expression can be written as
$$
R^{[\lambda]} \cdot T_i =  T_i\cdot R^{[s_i(\lambda)]} +
\sum_{j:\,\beta_j = \alpha_i} R^{[r_j(\lambda)]}
-\sum_{j:\,\beta_j = -\alpha_i} R^{[r_j(\lambda)]}.
$$

For a hyperplane $H$ of the form $H_{\alpha_i,k}$, $k\in\Z$, 
let $p_k$ be the number of times the alcove path 
$A_\circ\stackrel{\beta_l}\longrightarrow \cdots \stackrel{\beta_1} 
\longrightarrow A_{\lambda}$ crosses $H$ in the positive direction, 
and $n_k$ be the number of times the path crosses $H$ in the negative
direction.  In other words, $p_k = \#\{j \mid \beta_j = \alpha_i, \
r_j = s_{\alpha_i,k}\}$ 
and $n_k = \#\{j \mid \beta_j = - \alpha_i, \
r_j = s_{\alpha_i,k}\}$. 
Then $p_k - n_k$ is nonzero if and only if $H$ separates the alcoves
$A_\circ$ and $A_\lambda$.  More specifically,
$$
p_k-n_k = 
\left\{
\begin{array}{cl}
1  & \textrm{if }  0 < k \leq (\lambda,\alpha_i^\vee) , \\[.05in]
-1 & \textrm{if }  0 \geq k >  (\lambda,\alpha_i^\vee) , \\[.05in]
0  & \textrm{otherwise}.
\end{array}\right.
$$
This shows that 
$$
R^{[\lambda]} \cdot T_i =  T_i\cdot R^{[s_i(\lambda)]} +
\sum_{0<k\leq (\lambda,\alpha_i^\vee)}
R^{[s_{\alpha_i,k}(\lambda)]}  
- \sum_{(\lambda,\alpha_i^\vee)<k\leq 0} R^{[s_{\alpha_i,k}(\lambda)]} ,
$$
which is equivalent to the claim of the proposition.
\end{proof}

\section{The $K_T$-Chevalley Formula: Operator Notation}
\label{sec:k-chevalley-operator-notation}

We can formulate and prove our 
main result---the equivariant $K$-theory Chevalley 
formula---using the operator notation, as follows.
Recall that 
$$
R^{[\lambda]} = R_{\beta_l}\cdots R_{\beta_1} 
= X^\rho\,(X^{\beta_l} + B_{\beta_l})\cdots
(X^{\beta_2} + B_{\beta_2})\,
(X^{\beta_1} + B_{\beta_1})\, X^{-\rho},
$$
where $(\beta_1,\dots,\beta_l)$ is a $\lambda$-chain.

\begin{theorem}  
For any weight $\lambda$, the operator
$R^{[\lambda]}$ preserves the space $K_T(G/B)$.
For any $u\in W$, we have
$$
e^\lambda\cdot [\O_u] = R^{[\lambda]}([\O_u]),
$$
i.e., the operator $R^{[\lambda]}$ acts on the space $K_T(G/B)$ 
as the operator of multiplication by the class $e^\lambda$
of the corresponding  line bundle.
\label{thm:K-Chevalley}
\end{theorem}

\begin{proof}  Proposition~\ref{prop:same-relations} says
that the operators $R^{[\lambda]}$ satisfy the same commutation
relations with the elementary Demazure operators $T_i$ as the 
operators $E^\lambda$, see~(\ref{eq:YT=TY}). 
Also $R^{[\lambda]}  ([\O_1])  = x^\lambda\,[\O_1]$,
by Proposition~\ref{prop:sum-over-subsets}.
Now Lemma~\ref{lem:simple-observation} implies that 
the operator $R^{[\lambda]}$ preserves $K_T(G/B)\subset \tilde K_T(G/B)$
and acts as the operator $E^\lambda$ of multiplication by the class 
$e^\lambda$ of a line bundle.
\end{proof}

In Section~\ref{sec:central-points}, we show that
Theorem~\ref{thm:K-Chevalley} is equivalent to 
Theorem~{\rm \ref{th:main-combinatorial-formulation}}.
In Sections~\ref{sec:Type-A} and~\ref{sec:other-types}, we illustrate
Theorems~{\rm \ref{th:main-combinatorial-formulation}} 
and~\ref{thm:K-Chevalley} by several examples.

\begin{remark}\label{domwgt}
If $\lambda$ is a dominant weight, then, according to
Lemma~\ref{lem:reduced-dominant-antidominant},
the operator $R^{[\lambda]}$ expands as a
positive expression in the Bruhat operators $B_\alpha$, $\alpha\in\Phi^+$,
and the operators $X^\mu$.  
Indeed, a reduced $\lambda$-chain involves only positive roots.
In this case, Theorem~\ref{thm:K-Chevalley}
gives a positive formula for $e^\lambda\cdot [\O_u]$.
\end{remark}

Specializing $x^\mu\mapsto 1$, we obtain the 
nonequivariant $K$-theory Chevalley formula.
By a slight abuse of notation, we will use the same symbols $e^\lambda$ and $[\O_w]$ for the obvious classes in $K(G/B)$ as in $K_T(G/B)$.

\begin{corollary}
Let $\lambda\in\Lambda$ and 
$(\beta_1,\dots,\beta_l)$ be a $\lambda$-chain.
Then the operator
$$
R^{[\lambda]}_{x=1}= (1+ B_{\beta_l})\cdots (1+ B_{\beta_1})
$$
acts on the Grothendieck ring $K(G/B)$ as the operator
of multiplication by the class $e^\lambda$ of the corresponding 
line bundle.
\label{cor:non-equivariant-K-chevalley}
\end{corollary}

\begin{remark}\label{classchev} 
We claim that Corollary~\ref{cor:non-equivariant-K-chevalley} 
implies the classical Chevalley formula~(\ref{eq:chevalley}).   
In order to derive this formula, we need to collect linear terms 
in the expansion of the product $(1+ B_{\beta_l})\cdots (1+ B_{\beta_1})$.
Indeed, the coefficient $c_{u, u s_\alpha}^\lambda$, for  
$\ell(us_\alpha)=\ell(u)-1$, equals to the number of times the term
$B_\alpha$ appears in the expansion minus the number of times 
$B_{-\alpha}$ appears in the expansion.
According to Lemma~\ref{lem:chains=decompostions}, for any $\alpha\in\Phi^+$,
this coefficient is
$$ 
\#\{j \mid \beta_j = \alpha\} - \#\{j \mid \beta_j= -\alpha\} =
-m_{\alpha}(A_{-\lambda})= (\lambda,\alpha^\vee),  
$$
which is exactly the coefficient in the Chevalley formula.
Thus, (\ref{eq:K-theory->Chevalley}) 
and~(\ref{eq:chevalley})  follow. 
\end{remark}

\section{Central Points of Alcoves}
\label{sec:central-points}

In this section, we show that 
Theorem~{\rm \ref{th:main-combinatorial-formulation}}
is equivalent to Theorem~\ref{thm:K-Chevalley}. 
In order to do this, we show explicitly the way in which 
the operator $R^{[\lambda]}$ acts on basis elements $[\O_u]$.  
It is convenient to do this using central points of alcoves.
\medskip

Let us define the set $Z\subset\hR$ as
$$
Z := \{\zeta\in \Lambda/h \mid 
(\zeta,\alpha^\vee)\not\in\Z \textrm{ for any } \alpha\in\Phi\},
$$
i.e., $Z$ is the set of the elements of the lattice $\Lambda/h$ 
that do not belong to any hyperplane $H_{\alpha,k}$, where $h$
is the Coxeter number.
Then every element of $Z$ belongs to some alcove.
The affine Weyl group $\Waff$ preserves the set $Z$.  
This set was considered by Kostant~\cite{Kost}.

\begin{lemma} {\rm \cite{Kost}} \
Each alcove contains precisely one element of the set $Z$.
The only element of $Z$ in the fundamental alcove $A_\circ$ 
is $\rho/h$.
\label{lem:rho/hvee}
\end{lemma}

\begin{proof}
It is enough to prove the statement only for the fundamental alcove,
because $\Waff$ acts transitively on the alcoves. 
Let us express the highest coroot as a linear combination of simple coroots: 
$\theta^\vee=c_1\,\alpha_1^\vee+\cdots+c_r\,\alpha_r^\vee$.
Then $c_i$ are strictly positive integers and $h=c_1+\cdots+c_r+1$.
Every element $\zeta$ of $Z$ can be written as
$\zeta=(a_1 \,\omega_1+\cdots+a_r \, \omega_r)/h$, where
$a_1,\dots,a_r\in\Z$.
The condition that 
$\zeta \in Z\cap A_\circ$ can be written as $a_1,\dots,a_r>0$ and 
$(a_1\,c_1+\cdots + a_r\,c_r)/(c_1+\cdots+c_r+1)<1$, 
see~(\ref{eq:fund-alcove}). 
The only sequence of integers $(a_1,\dots,a_r)$ 
that satisfies these conditions is $(1,\dots,1)$.
Thus $Z\cap A_\circ$ consists of the single element 
$(\omega_1+\cdots+\omega_r)/h= \rho/h$.
\end{proof}

For an alcove $A$, the only element $\zeta_A$ of $Z\cap A$ 
is called the {\it central point\/} of the alcove $A$.
In particular, $\zeta_{A_\circ} = \rho/h$.
The map $A\mapsto\zeta_A$ is a one-to-one correspondence between
the set of all alcoves and $Z$.

\begin{lemma}
For a pair of adjacent alcoves $A\stackrel{\alpha}\longrightarrow B$,
we have $\zeta_B -\zeta_A = \alpha/h$.
\label{lem:paths-zeta}
\end{lemma}

\begin{proof}
It is enough
to prove this lemma for the fundamental alcove $A=A_\circ$.
All alcoves adjacent to $A_\circ$ are obtained from $A_\circ$
by the reflections $s_0,s_1,\dots, s_r$;  and
$A_\circ\stackrel{-\alpha_i}\longrightarrow s_i(A_\circ)$.
Applying these reflections to the central point 
$\zeta_{A_\circ}=\rho/h$, we obtain 
$s_i(\zeta_{A_\circ})-\zeta_{A_\circ} = -\alpha_i/h$,
for $i=0,\dots,r$.
\end{proof}

In fact, in the simply-laced case, the converse statement is true
as well.

\begin{lemma}  Suppose that $\Phi$ is a root system of type $A$-$D$-$E$.
Then $A\stackrel{\alpha}\longrightarrow B$ if and only if
$\zeta_B -\zeta_A = \alpha/h$.
\label{lem:zeta-alpha/hvee-inverse}
\end{lemma}

\begin{proof}
Again, we can assume that $A=A_\circ$ is the fundamental alcove.
In view of Lemma~\ref{lem:paths-zeta}, it remains
to show that $\mu=\rho/h+\alpha/h\not\in Z$,
for any root $\alpha\in\Phi\setminus\{-\alpha_1,\dots,-\alpha_r,\theta\}$.
For any such $\alpha$, there is a simple root $\alpha_i$ such that
$\alpha+\alpha_i$ is a root.  Thus $(\alpha,\alpha_i^\vee)=-1$ and
$(\mu,\alpha_i^\vee)=0$.  This implies that $\mu$ belongs to 
the hyperplane $H_{\alpha_i,0}$ and, thus, $\mu\not\in Z$.
\end{proof}

\begin{remark}
In the case of a nonsimply-laced root system, the statement converse to
Lemma~\ref{lem:paths-zeta} is not true.
In other words, there are nonadjacent alcoves $A$ and $B$ such that 
$\zeta_B-\zeta_A=\alpha/h$ for some root $\alpha$.
\end{remark}

Let us now fix an alcove path 
$A_\circ\stackrel{-\beta_1}\longrightarrow \cdots \stackrel
{-\beta_l} \longrightarrow A_{-\lambda}$ and the associated
$\lambda$-chain $(\beta_1,\dots,\beta_l)$.
By the definition, the operator $R^{[\lambda]}$ can be expressed as
\begin{equation}
R^{[\lambda]} = X^\rho\,(X^{\beta_l} + B_{\beta_l})\cdots
(X^{\beta_2} + B_{\beta_2})\,
(X^{\beta_1} + B_{\beta_1})\, X^{-\rho}.
\label{eq:R-rho-XD-rho}
\end{equation}
We can 
expand $R^{[\lambda]}$ as a sum of $2^l$ terms.  For a subset
$J\subset\{1,\dots,l\}$, let $R^{[\lambda]}_J$ be the term that contains
$B_{\beta_j}$, if $j\in J$, and $X^{\beta_j}$, otherwise.  
It is convenient to give the following interpretation 
for the term $R^{[\lambda]}_J$ using tail-flips.

Let $\pi=(0,\pi_0,\pi_1,\dots,\pi_l,\mu)$ be a collection of points in 
$\hR$.  We can think of this collection as a continuous piecewise-linear path 
in $\hR$ from $0$ to $\mu$.  Let $j$ be an index such that 
$\pi_{j-1}\ne\pi_{j}$, and let $r_j$ be the affine reflection with respect to 
the perpendicular bisector of the segment $[\pi_{j-1},\pi_j]$. 
In other words, the affine reflection $r_j$ is given by the condition
$r_j(\pi_{j-1})=\pi_j$.
For such an index $j$, we define the $j$-th {\it tail-flip\/} of $\pi$ as
$$
f_j(\pi) = (0,\pi_0,\dots,\pi_{j-1},r_j(\pi_{j+1}),\dots,r_j(\pi_l),r_j(\mu)).
$$
Then $f_j(\pi)$ corresponds to a path from $0$ to $r_j(\mu)$.
Let us associate with $\pi$ the following composition of operators
$$
X_\pi :=  X^{h(\pi_l-\mu)}
X^{h(\pi_{l-1}-\pi_{l})}
\cdots
X^{h(\pi_0-\pi_1)}
X^{h(0-\pi_0)} = X^{-h\mu}.
$$
Then $X_{f_j(\pi)} = X^{-h r_j(\mu)}$.

Let us now assume that $\pi=(0,\zeta_{A_0},\dots,\zeta_{A_l},-\lambda)$,
i.e., $\pi_i$'s are the central points of the alcoves $A_i$.
Then
$$
X_{\pi} = X^\rho\, X^{\beta_l}\cdots X^{\beta_1}\,X^{-\rho}
=X^{h\lambda}.
$$
Indeed, $h(0-\zeta_{A_\circ}) = -\rho$,
$h(\zeta_{A_{j-1}}-\zeta_{A_j})=\beta_j$, 
and $h(\zeta_{A_{-\lambda}}-(-\lambda)) = \rho$, 
see Lemmas~\ref{lem:rho/hvee} and \ref{lem:paths-zeta}.
The expression $X_\pi$ is precisely the term $R^{[\lambda]}_\emptyset$
in the expansion of~(\ref{eq:R-rho-XD-rho}).

In this case, $r_j$ is the affine reflection with respect to the common face of
$A_{j-1}$ and $A_j$ and $\bar r_j = s_{\beta_j}$,
for $j=1,\dots,l$.
Suppose that the subset $J$ consists of a single element $j$.
The corresponding term $R^{[\lambda]}_{\{j\}}$
in the expansion of~(\ref{eq:R-rho-XD-rho}) 
is obtained from the above expression $X_\pi$
by replacing the term $X^{\beta_j}$
with $B_{\beta_j}$.  Let us commute $B_{\beta_j}$ all the way 
to the left using relation~(\ref{eq:MQ=QM}).
We obtain
\begin{align*}
R^{[\lambda]}_{\{j\}} &= 
X^\rho\, X^{\beta_l}\cdots X^{\beta_{j+1}} B_{\beta_j} 
X^{\beta_{j-1}}\cdots X^{\beta_1}\,X^{-\rho}
 \\[.05in]
&=
B_{\beta_j}\,X^{\bar r_j(\rho)}\, 
X^{\bar r_j(\beta_l)}\cdots X^{\bar r_j(\beta_{j+1})} 
X^{\beta_{j-1}}\cdots X^{\beta_1}\,X^{-\rho}.
\end{align*}
The product of $X$'s in the last expression is precisely the 
operator $X_{f_j(\pi)}$ for the $j$-th tail-flip $\pi$.
In other words, $R^{[\lambda]}_{\{j\}} = B_{\beta_j}\, X_{f_j(\pi)}$.

In general, for a subset $J=\{j_1<\cdots< j_s \}\subset \{1,\dots,l\}$,
we have
$$
R^{[\lambda]}_{J} = B_{\beta_{j_s}}\cdots B_{\beta_{j_1}} 
X_{f_{j_1}\cdots f_{j_s}(\pi)}.
$$
Indeed, let us start with the expression $X_\pi$.  Replace the term 
$X^{\beta_{j_s}}$ in it with 
$B_{\beta_{j_s}}$, and commute it all the way to the left.  This 
leads to the expression $B_{\beta_{j_s}}\,X_{f_{j_s}(\pi)}$.
Then replace the term $X^{\beta_{j_{s-1}}}$ with $B_{\beta_{j_{s-1}}}$ 
and commute it to the left.
This leads to the expression 
$B_{\beta_{j_s}}\,B_{\beta_{j_{s-1}}}\,X_{f_{j_{s-1}}f_{j_s}(\pi)}$, etc.

We have
$$
X_{f_{j_{1}}\cdots f_{j_s}(\pi)} = 
X^{-h\, r_{j_1}\cdots r_{j_s}(-\lambda)}.
$$
According to (\ref{eq:Qlambda}),
this operator is explicitly given by
$$
X_{f_{j_{1}}\cdots f_{j_s}(\pi)}:[\O_u]\longmapsto
x^{-u\,r_{j_1}\cdots r_{j_s}(-\lambda)}\,[\O_u].
$$

Let us summarize our calculations.

\begin{proposition}
Let $\lambda\in\Lambda$ be a weight.  
Let $(r_1,\dots,r_l)$ and $(\beta_1,\dots,\beta_l)$ be the $\lambda$-chain of reflections and the $\lambda$-chain of roots associated with a decomposition $v_{-\lambda} = s_{i_1} \cdots s_{i_l}$.
Then the operator $R^{[\lambda]}$ is given by
$$
R^{[\lambda]} : [\O_u] \longmapsto
\sum_J  x^{-u\,r_{j_1}\cdots r_{j_s}(-\lambda)}\,
B_{\beta_{j_s}}\cdots B_{\beta_{j_1}}([\O_u]), 
$$
over all subsets $J=\{j_1<\cdots< j_s\}\subset\{1,\dots,l\}$.
\label{prop:sum-over-subsets}
\end{proposition}

We can now finish the proof 
Theorem~{\rm \ref{th:main-combinatorial-formulation}}.
\begin{proof}[Proof of Theorem~{\rm \ref{th:main-combinatorial-formulation}}]
This follows from Theorem~\ref{thm:K-Chevalley} and
Proposition~\ref{prop:sum-over-subsets}. 
\end{proof}

\section{Examples for Type $A$}
\label{sec:Type-A}

In this and the next sections we illustrate our results by 
presenting several examples.

Suppose that $G=\SL_n$.  Then the root system $\Phi$ is of type $A_{n-1}$
and the Weyl group $W$ is the symmetric group $S_n$.
We can identify the space $\h_\R^*$ with the quotient space 
$V:=\R^n/\R(1,\dots,1)$,
where $\R(1,\dots,1)$ denotes the subspace in $\R^n$ spanned 
by the vector $(1,\dots,1)$.  
The action of the symmetric group $S_n$ on $V$ is obtained 
from the (left) $S_n$-action on $\R^n$ by permutation of coordinates.
Let $\varepsilon_1,\dots,\varepsilon_n\in V$ 
be the images of the coordinate vectors in $\R^n$.
The root system $\Phi$ can be represented as 
$\Phi=\{\alpha_{ij}:=\varepsilon_i-\varepsilon_j \mid i\ne j,\ 1\leq i,j\leq n\}$.
The simple roots are $\alpha_i=\alpha_{i\,i+1}$, 
for $i=1,\dots,n-1$.
The longest coroot
is $\theta^\vee = \alpha_{1n}^\vee$.
The fundamental weights are $\omega_i = \varepsilon_1+\cdots +\varepsilon_i$, 
for $i=1,\dots,n-1$. 
We have $\rho = n\varepsilon_1+(n-1)\varepsilon_2+\cdots + 2 \varepsilon_{n-1}+
\varepsilon_{n}$, and the Coxeter number is $h=(\rho,\theta^\vee)+1 = n$.
The weight lattice is $\Lambda=\Z^n/\Z(1,\dots,1)$. 
We use the notation $[\lambda_1,\dots,\lambda_n]$ for a weight, as the coset of $(\lambda_1,\dots,\lambda_n)$ in $\Z^n$.  

Let $nZ\subset\Lambda$ be the set $Z$ of central points of alcoves scaled 
by the factor $h=n$.  
The fundamental alcove corresponds to the point 
$\rho$ in $nZ$.  
According Lemma~\ref{lem:zeta-alpha/hvee-inverse}, two alcoves
are adjacent $A\stackrel{\alpha}\longrightarrow B$, $\alpha\in\Phi$,
if and only if the corresponding elements of $nZ$ 
are related by $n\zeta_B-n\zeta_A=\alpha$. 
In this case, we write $n\zeta_A \stackrel{\alpha}\longrightarrow
n\zeta_B$.  Thus, we have the structure of a directed graph with labeled edges
on the set $nZ$.  Alcove paths correspond to paths in this graph.
The set $nZ$ can be explicitly described as
$$
nZ=\{[\mu_1,\dots,\mu_n]\in\Lambda \mid
 \mu_1,\dots,\mu_n 
\textrm{ have distinct residues modulo $n$}\}.
$$
For an element $\mu=[\mu_1,\dots,\mu_n]\in nZ$, there exists an edge
$\mu\stackrel{\alpha_{ij}}\longrightarrow (\mu+\alpha_{ij})$
if and only if $\mu_i+1\equiv \mu_j\ \textrm{mod} \ n$. 
Given a weight $\lambda$, the corresponding $\lambda$-chains are in one-to-one
correspondence with directed paths in the graph $nZ$ from $\rho$ to 
$\rho-n\lambda$.

\begin{example} 
Suppose that $n=4$ and $\lambda=\omega_2=[1,1,0,0]$.
The directed path
$$
[4,3,2,1]\stackrel{-\alpha_{23}}\longrightarrow [4,2,3,1] 
\stackrel{-\alpha_{13}}\longrightarrow [3,2,4,1] 
\stackrel{-\alpha_{24}}\longrightarrow [3,1,4,2] 
\stackrel{-\alpha_{14}}\longrightarrow [2,1,4,3]
$$
from $\rho=[4,3,2,1]$ to $\rho-n\,\omega_2= [0,-1,2,1]=[2,1,4,3]$
produces the $\omega_2$-chain 
$(\alpha_{23},\alpha_{13},\alpha_{24},\alpha_{14})$.
\end{example}

\begin{example}
For an arbitrary $n$, we have $\omega_1=\varepsilon_1=[1,0,\dots,0]$.
The path
$$
\begin{array}{l}
[n,n-1,\dots,1] \stackrel{-\alpha_{12}}\longrightarrow [n-1,n,n-2,\dots,1] 
\stackrel{-\alpha_{13}}\longrightarrow [n-2,n,n-1,n-3,\dots,1] 
\\[.05in]
\qquad 
\qquad 
\stackrel{-\alpha_{14}}\longrightarrow 
[n-3,n,n-1,n-2,n-4,\dots,1] \stackrel{-\alpha_{15}}\longrightarrow 
\cdots \stackrel{-\alpha_{1n}}\longrightarrow 
[1,n,n-1,\dots,2].
\end{array}
$$
from $\rho$ to $\rho-n\,\omega_1$ gives the $\omega_1$-chain
$(\alpha_{12},\alpha_{13},\alpha_{14},\dots,\alpha_{1n})$.
In general, for any $k=1,\dots,n$,
we have the $\varepsilon_k$-chain
\begin{equation}
(\alpha_{k\,k+1},\alpha_{k\,k+2},\dots,\alpha_{k\,n},\alpha_{k\,1},
\alpha_{k\,2},\dots,\alpha_{k\,k-1})
\label{eq:epsilon-k-chain}
\end{equation}
given by the corresponding path from $\rho$ to $\rho-n\varepsilon_k$.
\end{example}

Recall that $v_{-\lambda}$ is the unique element of $\Waff$ such that 
$v_{-\lambda}(A_\circ)=A_{-\lambda}$.  Equivalently, we can define 
$v_{-\lambda}$ in terms of central points of alcoves by the condition
$v_{-\lambda} (\rho/h) = \rho/h - \lambda$.  

\begin{lemma}  
Suppose that $\Phi$ is of type $A_{n-1}$.
Then, for $k=1,\dots,n-1$, the affine Weyl group element
$v_{-\omega_k}$ belongs, in fact, to $S_n\subset \Waff$.
This permutation is given by 
$$
v_{-\omega_k} = \left(\begin{array}{ccccccc}
1 & 2 & \cdots & n-k & n-k+1 & \cdots & n \\
k+1 & k+2 & \cdots & n & 1 & \cdots & k
\end{array}
\right)\in S_n\subset \Waff.
$$
\label{lem:v-omega-k}
\end{lemma}

\begin{proof}
This permutation maps $\rho=[n,\dots,1]$ to
$[k, k-1, \dots, 1, n, n-1,\dots,k+1]= [0,-1,\dots,-k+1,n-k,n-k-1,\dots,1] 
= \rho-n\,\omega_k$, as needed.
\end{proof}

Let $R_{ij}:=R_{\alpha_{ij}}$.
Theorem~\ref{thm:K-Chevalley} implies the following statement.

\begin{corollary}  For $k=1,\dots,n$, 
the operator of multiplication
by $e^{\varepsilon_k}$ in the Grothendieck ring $K_T(\SL_n/B)$ is given by
$$
R^{[\varepsilon_k]} =
R_{k\,k-1} R_{k\,k-2} \cdots R_{k\,1} R_{k\,n} R_{k\,n-1}\cdots R_{k\,k+1}.
$$
For $k=1,\dots,n-1$, the operator of multiplication by the class $e^{\omega_k}$ corresponding to the $k$-th fundamental weight
$\omega_k$ is given by
\begin{equation}
R^{[\omega_k]}  = R^{[\varepsilon_1]}  \cdots R^{[\varepsilon_k]}   = 
\prod_{i=1,\dots,k}^{\longrightarrow}\quad
\prod_{j=k+1,\dots,n}^{\longleftarrow} R_{ij}.
\label{eq:type-A-fundamental}
\end{equation}
\label{cor:type-A-epsilon-omega}
\end{corollary}

The combinatorial formula for multiplication by $e^{\omega_k}$ 
in the Grothendieck ring $K(\SL_n/B)$ 
that follows from formula~(\ref{eq:type-A-fundamental}) was
originally found in~\cite{Len}.

\begin{proof}
The expression for $R^{[\varepsilon_k]}$ is given by the
$\varepsilon_k$-chain~(\ref{eq:epsilon-k-chain}).
The expression for $R^{[\omega_k]}$ can be obtained by simplifying 
$R^{[\varepsilon_1]}\cdots R^{[\varepsilon_k]}$, as shown in~\cite{Len}. 
Alternatively, the reduced decomposition 
$v_{-\omega_k} = (s_k \cdots s_{n-1})(s_{k-1}\cdots s_{n-2})\cdots
(s_1 \cdots s_{n-k})$
for the permutation $v_{-\omega_k}$ given by Lemma~\ref{lem:depends-only-on} 
corresponds to an $\omega_k$-chain, see Definition~\ref{def:lambda-chain}.
This $\omega_k$-chain produces the needed expression for
$R^{[\omega_k]}$.
\end{proof}

\begin{example}
For $n=3$, Corollary~\ref{cor:type-A-epsilon-omega} says that
$$
R^{[\omega_1]}=R_{13} \,R_{12}\quad\textrm{and}\quad
R^{[\omega_2]}=R_{13} \, R_{23}.
$$
\end{example}

For a weight $\lambda=a_1\omega_1+\cdots + a_r \omega_r$,
we can obtain an expression for $R^{[\lambda]}$ by concatenation
of $a_1$ copies of $R^{[\omega_1]}$, $a_2$ copies of $R^{[\omega_2]}$,
etc.

Theorem~\ref{th:main-combinatorial-formulation} says that
that the coefficient of $[\O_w]$ in the product $e^\lambda\cdot [\O_u]$
in $K_T(G/B)$ is given by the sum over subsequences in 
the $\lambda$-chain $(\beta_1,\dots,\beta_l)$
that give saturated decreasing chains $u\gtrdot \cdots \gtrdot w$ 
in the Bruhat order on $W$.
Let us illustrate this theorem by the following two examples.

\begin{example} Suppose that $n=3$, $\lambda=\omega_1$, and 
$u=w_\circ=s_1s_2s_1\in W$.  
Let us calculate the product $e^{\lambda}\cdot [\O_{u}]$
in $K_T(\SL_n/B)$
using Theorem~\ref{th:main-combinatorial-formulation}.
The $\omega_1$-chain $(\beta_1,\beta_2)=(\alpha_{12},\alpha_{13})$ 
is associated with the reduced decomposition $s_1s_2=v_{-\omega_1}$.  
The corresponding $\omega_1$-chain of reflections is
$(r_1,r_2)=(s_1,s_1s_2s_1)= (s_{\alpha_{12},0},s_{\alpha_{13},0})$.
Three out of four subsequences in $(\beta_1,\beta_2)$ correspond
to decreasing chains in Bruhat order starting at $w_\circ$:
(empty subsequence), $(\alpha_{12})$, and $(\alpha_{12},\alpha_{13})$.
Thus we have
$$
e^{\omega_1}\cdot [\O_{w_\circ}] = x^{-w_\circ(-\omega_1)} [\O_{w_\circ}] +
x^{-w_\circ r_1(-\omega_1)}[\O_{s_1s_2}]
+ x^{-w_\circ r_1r_2(-\omega_1)}[\O_{s_2}].
$$
We can write this expression as
$$
e^{[1,0,0]}\cdot [\O_{w_\circ}] = x^{[0,0,1]} [\O_{w_\circ}] +
x^{[0,1,0]}[\O_{s_1 s_2}]
+ x^{[1,0,0]}[\O_{s_2}].
$$
The character of the irreducible representation $V_{\omega_1}$
is obtained from the
right-hand side of this expression by replacing each term $x^\mu[\O_w]$
with $e^\mu$:
$$
ch(V_{\omega_1})= e^{[0,0,1]} + e^{[0,1,0]} + e^{[1,0,0]}.
$$
\end{example}

Let us give a less trivial example.
\begin{example}
Suppose $n=3$ and $\lambda = 2\,\omega_1+\omega_2=[3,1,0]$.  
The path 
$$
\begin{array}{l}
[3,2,1] \stackrel{-\alpha_{12}}\longrightarrow [2,3,1]
\stackrel{-\alpha_{13}}\longrightarrow [1,3,2]
\stackrel{-\alpha_{23}}\longrightarrow [1,2,3]
\qquad\qquad
\\[.05in]
\qquad\qquad\qquad\qquad\stackrel{-\alpha_{13}}\longrightarrow [0,2,4]
\stackrel{-\alpha_{12}}\longrightarrow [-1,3,4]
\stackrel{-\alpha_{13}}\longrightarrow [-2,3,5]
\end{array}
$$
from $\rho=[3,2,1]$ to $\rho-n\lambda=[-2,3,5]$ gives
the $\lambda$-chain
$$
(\beta_1,\dots,\beta_6) = 
(\alpha_{12},\ \alpha_{13},\ \alpha_{23},\ \alpha_{13},\ 
\alpha_{12},\ \alpha_{13}),
$$
which is associated with 
the reduced decomposition $v_{-\lambda}=s_1 s_2 s_1 s_0 s_1 s_2$ 
in the affine Weyl group.
We have
$$
R^{[\lambda]} = 
R_{\beta_6}\cdots R_{\beta_1}=
R_{13}\,R_{12}\,R_{13} \,R_{23}\, R_{13}\, R_{12} =
R^{[\omega_1]} \,R^{[\omega_2]} \,R^{[\omega_1]}.
$$
The corresponding $\lambda$-chain of reflections is 
$$
(r_1,\dots,r_6) = (s_{\alpha_{12},0},\  
s_{\alpha_{13},0},\  s_{\alpha_{23},0},\ s_{\alpha_{13},-1},\  
s_{\alpha_{12},-1},\  s_{\alpha_{13},-2}).
$$

Suppose that $u=s_2 s_1$.  There are five 
saturated chains in Bruhat order descending from $u$: 
(empty chain), 
$(u\gtrdot u s_{\alpha_{12}} = s_2)$,
$(u\gtrdot u s_{\alpha_{13}} = s_1)$,
$(u\gtrdot u s_{\alpha_{12}} \gtrdot u s_{\alpha_{12}} s_{\alpha_{23}} = 1)$,
$(u\gtrdot u s_{\alpha_{13}} \gtrdot u s_{\alpha_{13}} s_{\alpha_{12}} = 1)$.
Thus, the expansion of $e^\lambda\cdot [\O_u]$ is given by
the sum over the following subsequences in the $\lambda$-chain $(\beta_1,\dots,\beta_6)$:
$$
\textrm{(empty subsequence)}, \ (\alpha_{12}), \  (\alpha_{13}),
(\alpha_{12},\,\alpha_{23}), \ 
(\alpha_{13},\,\alpha_{12}).
$$
The sequence $(\beta_1,\dots,\beta_6)$ contains one empty subsequence, 
two subsequences of the form $(\alpha_{12})$,
three subsequences of the form $(\alpha_{13})$,
one subsequence of the form $(\alpha_{12},\alpha_{23})$, and 
two subsequence of the form $(\alpha_{13},\alpha_{12})$.
Hence, we have
$$
\begin{array}{l}
e^\lambda\cdot [\O_{s_2s_1}] = x^{-u(-\lambda)} \,[\O_{s_2s_1}] +
\left(x^{-u r_1(-\lambda)} +  x^{-u r_5(-\lambda)}\right) \,[\O_{s_2}] \; +
\\[.05in]
\qquad + \;
\left(x^{-u r_2(-\lambda)} +  x^{-u r_4(-\lambda)} +
x^{-u r_6(-\lambda)}\right)\,[\O_{s_1}] \;+
\\[.05in]
\qquad + \;
 x^{-u r_1 r_3 (-\lambda)}\,[\O_{1}] +
 \left(x^{-u r_2 r_5 (-\lambda)} + x^{-u r_4 r_5 (-\lambda)}\right)\,[\O_{1}].
\end{array}
$$
We can explicitly write this expression as
$$
\begin{array}{l}
e^{[3,1,0]}\cdot [\O_{s_2s_1}] = x^{[1,0,3]} \,[\O_{s_2s_1}] +
\left(x^{[3,0,1]} +  x^{[2,0,2]}\right) \,[\O_{s_2}] \; +
\\[.05in]
\qquad + \;
\left(x^{[1,3,0]} +  x^{[1,2,1]} +
x^{[1,1,2]}\right)\,[\O_{s_1}] 
+ x^{[3,1,0]}\,[\O_{1}] 
+ \left(x^{[2,2,0]} + x^{[2,1,1]}\right)\,[\O_{1}].
\end{array}
$$

The corresponding Demazure character is
$$
\begin{array}{l}
ch(V_{[3,1,0],s_2s_1}) = 
\\[.05in]
\quad e^{[1,0,3]}  +
e^{[3,0,1]} +  e^{[2,0,2]} +
e^{[1,3,0]} +  e^{[1,2,1]} + e^{[1,1,2]} 
+ e^{[3,1,0]} + e^{[2,2,0]} + e^{[2,1,1]}.
\end{array}
$$
\end{example}

\section{Examples for Other Types}
\label{sec:other-types}

For an arbitrary root system, we can use the explicit construction
of the $\lambda$-chain of reflections $(r_1,\dots,r_l)$ and the 
$\lambda$-chain of roots $(\beta_1,\dots,\beta_l)$ given by
Proposition~\ref{prop:construction-for-lambda-chain}.

\begin{example}  Suppose that the root system $\Phi$ is of type $G_2$.
Let us find $\lambda$-chains for $\lambda=\omega_1$ and 
$\lambda=\omega_2$ using
Proposition~\ref{prop:construction-for-lambda-chain}.
The positive roots are $\gamma_1=\alpha_1,\ 
\gamma_2=3\alpha_1+\alpha_2,\ \gamma_3=2\alpha_1+\alpha_2,\ 
\gamma_4=3\alpha_1+2\alpha_2,\ \gamma_5=\alpha_1+\alpha_2,\
\gamma_6=\alpha_2$.
The corresponding coroots are
$\gamma_1^\vee=\alpha_1^\vee, \
\gamma_2^\vee= \alpha_1^\vee+\alpha_2^\vee,\
\gamma_3^\vee = 2\alpha_1^\vee+3\alpha_2^\vee,\
\gamma_4^\vee = \alpha_1^\vee+2\alpha_2^\vee,\
\gamma_5^\vee = \alpha_1^\vee+3\alpha_2^\vee,\ 
\gamma_6^\vee=\alpha_2^\vee$.

Suppose that $\lambda=\omega_1$.
The set $\mathcal{R}_{\omega_1}$ of affine reflections with respect 
to the hyperplanes separating the alcoves $A_\circ$ and $A_{-\omega_1}$ is
$$
\mathcal{R}_{\omega_1} =\{s_{\gamma_1,0},\ s_{\gamma_2,0},\ s_{\gamma_3,0},\ 
s_{\gamma_3,-1},\ s_{\gamma_4,0},\ s_{\gamma_5,0}\}.
$$
The map $h:\mathcal{R}_{\omega_1}\to\R^{r+1}$ 
given by~(\ref{eq:map-h}) sends these affine reflections to the vectors
$$
\begin{array}{l}
(0,1,0),\ (0,1,1),\ (0,1,\frac32),\ (\frac12,1,\frac32),\  (0,1,2),\ (0,1,3),
\end{array}
$$
respectively.  The lexicographic order on vectors in $\R^3$ induces
the following total order on the set
$\mathcal{R}_{\omega_1}$: 
$$
s_{\gamma_1,0}<s_{\gamma_2,0}<s_{\gamma_3,0}<
s_{\gamma_4,0}<s_{\gamma_5,0}< s_{\gamma_3,-1}\,.
$$

Suppose now that $\lambda=\omega_2$.
The set $\mathcal{R}_{\omega_2}$ of affine reflections with respect 
to the hyperplanes separating $A_\circ$ and $A_{-\omega_2}$ is
$$
\mathcal{R}_{\omega_2} =\{s_{\gamma_2,0},\ s_{\gamma_3,0},\ s_{\gamma_3,-1},\ 
s_{\gamma_3,-2},\ s_{\gamma_4,0},\ s_{\gamma_4,-1},\
s_{\gamma_5,0},\ s_{\gamma_5,-1},\  s_{\gamma_5,-2}, s_{\gamma_6,0}
\}.
$$
The map $h:\mathcal{R}_{\omega_2}\to\R^{r+1}$ 
sends these affine reflections to the vectors
$$
\begin{array}{l}
(0,1,1),\ (0,\frac23,1),\ (\frac13,\frac23,1),\ (\frac23,\frac23,1),\  
(0,\frac12,1),\ (\frac12,\frac12,1),\ 
\\[.05in]
\qquad\qquad 
\qquad\qquad 
\qquad\qquad 
(0,\frac13,1),\ (\frac13,\frac13,1),\ 
(\frac23,\frac13,1),\ (0,0,1),
\end{array}
$$
respectively.  The lexicographic order on vectors in $\R^3$ induces
the following total order on $\mathcal{R}_{\omega_2}$:
$$
s_{\gamma_6,0}<s_{\gamma_5,0}<s_{\gamma_4,0}<s_{\gamma_3,0}<s_{\gamma_2,0}
<s_{\gamma_5,-1}<s_{\gamma_3,-1}<
s_{\gamma_4,-1}< s_{\gamma_5,-2}<s_{\gamma_3,-2}\,.
$$

The total orders on $\mathcal{R}_{\omega_1}$ and 
$\mathcal{R}_{\omega_2}$ correspond to
the $\omega_1$-chain 
$(\gamma_1,\gamma_2,\gamma_3,\gamma_4,\gamma_5,\gamma_3)$
and  the $\omega_2$-chain
$(\gamma_6,\gamma_5,\gamma_4,\gamma_3,\gamma_2,\gamma_5,\gamma_3,\gamma_4,
\gamma_5,\gamma_3)$.
Thus, the operators of multiplication by the classes $e^{\omega_1}$ 
and $e^{\omega_2}$ in $K_T(G/B)$ are given by
$$
\begin{array}{l}
R^{[\omega_1]} = 
R_{\gamma_3}\,R_{\gamma_5}\,R_{\gamma_4}\,R_{\gamma_3}\,
R_{\gamma_2}\,R_{\gamma_1},
\\[.05in]
R^{[\omega_2]} = 
R_{\gamma_3}\,R_{\gamma_5}\,R_{\gamma_4}\,R_{\gamma_3}\,R_{\gamma_5}\,
R_{\gamma_2}\,R_{\gamma_3}\,R_{\gamma_4}\,R_{\gamma_5}\,R_{\gamma_6}.
\end{array}
$$
\label{example:G-2}
\end{example}

By Lemma~\ref{lem:v-omega-k}, the element $v_{-\omega_k}$ belongs
to the (nonaffine) Weyl group $W$, for all fundamental
weights $\omega_k$ in type $A$.  Let us show that a similar phenomenon occurs 
for minuscule weights in other types as well.  
Recall that a dominant weight $\lambda$ is {\it minuscule\/} if the set
of weights in the $G$-module $V_\lambda$ is in the orbit 
$W\cdot \lambda$ of the Weyl group. 

\begin{lemma}
Let $\lambda\in\Lambda^+$.  Then 
$v_{-\lambda}\in W$ if and only if 
$\lambda$ is a minuscule weight. 
\end{lemma}

\begin{proof}  Let 
$(\beta_1,\dots,\beta_l)$ be a reduced $\lambda$-chain of roots, and 
let $(r_1,\dots,r_l)$ be the corresponding $\lambda$-chain
of reflections.  
According to Lemmas~\ref{lem:chains=decompostions} 
and~\ref{lem:reduced-dominant-antidominant}, the following statements
are equivalent:
(1) $v_{-\lambda}\in W$;
(2) $r_1,\dots,r_l\in W$;
(3) all (positive) roots $\beta_1,\dots,\beta_l$ are distinct;
(4) $(\lambda,\alpha^\vee)=0\textrm{ or }1$, 
for any $\alpha\in\Phi^+$. 
According to Corollary~\ref{cor:model-for-character},
the condition $r_1,\dots,r_l\in W$ implies that all weights 
in $V_\lambda$ are in the $W$-orbit $W\cdot \lambda$ and, thus, $\lambda$
is minuscule.
On the other hand, if $\lambda$ is minuscule, then 
$(\lambda,\alpha^\vee)=0\textrm{ or }1$, 
for any $\alpha\in\Phi^+$. Otherwise, if 
$(\lambda,\alpha^\vee)\geq 2$, then $V_\lambda$ contains the weight 
$\lambda-\alpha\not\in W\cdot \lambda$.
\end{proof}

The last two examples concern minuscule weights in types $B$ and $C$.
Recall that the element $v_{-\lambda}$ is uniquely defined by the condition
$v_{-\lambda}(\rho/h) = \rho/h-\lambda$.
If $v_{-\lambda}\in W$.  We can rewrite this condition as 
$v_{-\lambda}(\rho) = \rho - h\,\lambda$.

\begin{example}  
Suppose that $\Phi$ is of type $C_r$.  This root system can be 
embedded into $\R^r$ as follows:  $\Phi=\{\pm\varepsilon_i\pm\varepsilon_j,\ 
\pm 2\varepsilon_i \mid i\ne j\}$, where $\varepsilon_1,\dots,\varepsilon_r$
are the coordinate vectors in $\R^r$.  The simple roots
are $\alpha_1=\varepsilon_1-\varepsilon_2$,
$\alpha_2=\varepsilon_2-\varepsilon_3$, \dots
$\alpha_{r-1}=\varepsilon_{r-1}-\varepsilon_r$,
$\alpha_{r}=2\varepsilon_{r}$.
The Weyl group $W$ is the semidirect product of $S_r$ and $(\Z/2\Z)^r$.
It acts on $\R^r$ by permuting the coordinates and changing their signs.
The fundamental weights are $\omega_k= \varepsilon_1+\cdots + \varepsilon_k$, 
$k=1,\dots,r$. We have $\rho = (r,\dots,1)\in\R^r$, and 
the Coxeter number is $h = (\rho,\theta^\vee)+1=2r$.

Suppose that $\lambda=\omega_1$.  Then $\rho-h\omega_1 = 
(-r,r-1,r-2,\dots,1)\in\R^r$.  This weight is obtained from $\rho$
by applying the Weyl group element $s_{2\varepsilon_1}$ that changes 
the sign of the 
first coordinate.  Thus $v_{-\omega_1} = s_{2\varepsilon_1}
\in W\subset \Waff$. 
The only reduced decomposition of this element is
$v_{-\omega_1} = s_1 \cdots s_{r-1} \, s_r \, s_{r-1} \cdots s_1$,
so $\ell(v_{-\omega_1}) = 2r-1$.
This reduced decomposition corresponds to the $\omega_1$-chain 
$$
\begin{array}{l}
(\alpha_1,\ s_1(\alpha_2),\ s_1 s_2(\alpha_3),\ \dots,\ 
s_1 \dots s_{r-1}(\alpha_r),\ \dots,\ 
s_1 \dots s_r \dots s_2 (\alpha_1))=\\[.05in]
(\varepsilon_1-\varepsilon_2,\
\varepsilon_1-\varepsilon_3,\ \cdots,\
\varepsilon_1-\varepsilon_r, \ 
2\varepsilon_{1},\
\varepsilon_1+\varepsilon_r,\ \cdots,\
\varepsilon_1+\varepsilon_{3},\
\varepsilon_1+\varepsilon_{2}),\
\end{array}
$$
cf. Definition~\ref{def:lambda-chain}.
 The operator $R^{[\omega_1]}$ is given by
$$
R^{[\omega_1]} = 
R_{\varepsilon_1+\varepsilon_{2}}
R_{\varepsilon_1+\varepsilon_{3}}
\cdots
R_{\varepsilon_1+\varepsilon_r}
R_{2\varepsilon_{1}}
R_{\varepsilon_1-\varepsilon_r}
\cdots
R_{\varepsilon_1-\varepsilon_3}
R_{\varepsilon_1-\varepsilon_2}\,.
$$
\end{example}

\begin{example}  
Suppose that $\Phi$ is of type $B_r$.  This root system can be 
embedded into $\R^r$ as follows:  $\Phi=\{\pm\varepsilon_i\pm\varepsilon_j,\ 
\pm\varepsilon_i \mid i\ne j\}$, where $\varepsilon_1,\dots,\varepsilon_r$
are the coordinate vectors in $\R^r$.  The simple roots
are $\alpha_1=\varepsilon_1-\varepsilon_2$,
$\alpha_2=\varepsilon_2-\varepsilon_3$, \dots
$\alpha_{r-1}=\varepsilon_{r-1}-\varepsilon_r$,
$\alpha_{r}=\varepsilon_{r}$.
The Weyl group $W$ and its action on $\R^r$ are the same as in type $C_r$.
The fundamental weights are $\omega_k= \varepsilon_1+\cdots + \varepsilon_k$, 
$k=1,\dots,r-1$,
and $\omega_r=\frac12 (\varepsilon_1+\cdots+\varepsilon_r)$. 
We have $\rho = (r-\frac12,\dots,1-\frac12)\in\R^r$, and 
$h = (\rho,\theta^\vee)+1=2r$.

Suppose that $\lambda=\omega_r$ is the last fundamental weight.  
Then $\rho-h\omega_r = 
(-\frac 12,-1-\frac12,-2-\frac 12,\dots,-r+\frac12)\in\R^r$.  
This weight is obtained from $\rho$ by applying the Weyl group element
$v_{-\omega_r}\in W\subset\Waff$ that reverses the order of all 
coordinates and changes their signs.
The element $v_{-\omega_r}\in W$ has length 
$\ell(v_{-\omega_r})= r(r+1)/2$.  One of the reduced decompositions for 
this element is 
$$
v_{-\omega_r} = (s_r)(s_{r-1}\, s_r)(s_{r-2} \,s_{r-1}\, s_r)\cdots
(s_2\cdots s_r)(s_1\cdots s_r).
$$
The associated $\omega_r$-chain is
$(\alpha_r,\, s_r(\alpha_{r-1}),\, 
s_r s_{r-1}(\alpha_{r}),\, s_r s_{r-1} s_{r}(\alpha_{r-2}),\, \dots)$.
We can explicitly find the roots in this $\omega_r$-chain and write
the operator $R^{[\omega_r]}$ as 
$$
\begin{array}{l}
R^{[\omega_r]}=
(R_{\varepsilon_1}\, R_{\varepsilon_1+\varepsilon_2}\, 
R_{\varepsilon_1+\varepsilon_3}\cdots R_{\varepsilon_1+\varepsilon_r})
(R_{\varepsilon_2}\, R_{\varepsilon_2+\varepsilon_3}\, 
R_{\varepsilon_2+\varepsilon_4}\cdots R_{\varepsilon_2+\varepsilon_r})
\cdots \\[.05in]
\qquad \qquad \qquad 
\qquad \qquad 
\cdots (R_{\varepsilon_{r-2}}\, R_{\varepsilon_{r-2}+\varepsilon_{r-1}}\,
R_{\varepsilon_{r-2}+\varepsilon_r}\,)
(R_{\varepsilon_{r-1}}\, R_{\varepsilon_{r-1}+\varepsilon_r})
(R_{\varepsilon_r}).
\end{array}
$$
\end{example} 

\section{Quantum $K$-theory}
\label{sec:Quantum}

In this section, we conjecture a natural Chevalley-type formula in the
{\it quantum $K$-theory\/} of $G/B$. 
The quantum $K$-theory, which is a $K$-theoretic
version of {\it quantum cohomology}, was introduced by Lee~\cite{Lee}. 
The quantum $K$-theory of flag varieties, in particular, has been first 
studied by Givental and Lee~\cite{GiLe}. We recall a few basic facts below.

Let us denote by $QK(G/B)$ the quantum $K$-theory of $G/B$. In order to
describe it, we associate a variable $q_i$ to each simple root $\alpha_i$, 
and let $\Z[q]=\Z[q_1,\dots,q_r]$ be the polynomial ring in the $q_i$.
Given a collection of nonnegative integers $d=(d_1,\dots,d_r)$,
called multidegree, we let
$q^d:=q_1^{d_1}\ldots q_r^{d_r}$. 
As a $\Z[q]$-module, the quantum $K$-theory is defined as
$QK(G/B):= K(G/B)\otimes_\Z \Z[q]$. 
Let $[w]$ denote the class of the structure sheaf of the Schubert variety
$X_{w_\circ w}$.
Then the classes of $[w]$ form a $\Z[q]$-basis of $QK(G/B)$. 
The multiplication in $QK(G/B)$ is a deformation of the classical
multiplication:
$$
[u]\circ [v]=\sum_d q^d\sum_{w\in W} N_{uv}^w(d) \, [w]\,,
$$
where the first sum is over all multidegrees $d$, and $N_{uv}^w(d)$ is the 
{\it
quantum $K$-invariant of Gromov-Witten type} for $[u]$, $[v]$, and the
quantum dual of $[w]$. As defined in \cite{Lee}, this invariant is
the $K$-theoretic push-forward to $\mathrm{Spec}\, \C$ of some natural vector
bundle on the moduli space $\overline{M}_{3,0}(G/B,d)$ (via the orientation
defined by the virtual structure sheaf). 
The associativity of the quantum $K$-product
was established in \cite{Lee}, based on a sheaf-theoretic version of
an argument of WDVV-type.

Let us recall the Chevalley-type formula for the small quantum cohomology
ring $QH^*(G/B)$ of $G/B$.  For type $A$, this formula was first proved
in~\cite{FGP}.  In general type, it was proved by D.~Peterson (unpublished) and
by Fulton and Woodward~\cite{FW} (who, in fact, obtained a more general formula
for $G/P$).  
Again, as a $\Z[q]$-module, $QH^*(G/B):=H^*(G/B)\otimes \Z[q]$.
Thus, the quantum cohomology ring has a $\Z[q]$-basis basis 
given by the cohomology classes of $X_{w_\circ w}$, which we denoted by 
$\langle w\rangle$. 

The Chevalley-type formula in $QH^*(G/B)$ can be stated using the 
{\em quantum Bruhat operators} defined in \cite{BFP}. These are
operators on the group algebra $\Z[q][W]$ of the Weyl group $W$ over $\Z[q]$.
For each positive root $\alpha$, 
the quantum Bruhat operator $Q_\alpha$ is defined by 
$$
Q_{\alpha}(w)=\left\{ \begin{array}{ll}  w s_\alpha &\mbox{if}\;\; \ell(w
s_\alpha)=\ell(w)+1, \\ 
q^{d(\alpha)}\,w s_\alpha &\mbox{if}\;\; \ell(w
s_\alpha)=\ell(w)-2\,\mathrm{ht}(\alpha^\vee)+1,\\ 
0 &\mbox{otherwise\,,} \end{array}
\right.  
$$
where $\mathrm{ht}(\alpha^\vee)=(\rho,\alpha^\vee)$ is the height of 
the coroot $\alpha^\vee$,
and $q^{d(\alpha)} = q_1^{d_1} \cdots q_r^{d_r}$,
for $\alpha^\vee = d_1 \alpha_1^\vee + \cdots + d_r \alpha_r^\vee$,
i.e., $d_i = (\omega_i,\alpha^\vee)$.
Also define $Q_{\alpha}:=-Q_{-\alpha}$ if $\alpha$ is a negative root.  It was
proved in \cite{BFP} that the operators $Q_{\alpha}$ satisfy the Yang-Baxter
equation. 

The map $w\mapsto \langle w\rangle$ extends linearly to 
the isomorphism $\Z[q][W]\rightarrow QH^*(G/B)$
of $\Z[q]$-modules, for which we use the same notation $a\to \langle a\rangle$. 
Similarly, we extend the map $w\mapsto  [w]$.
The Chevalley formula in quantum cohomology can now be stated, as follows,
see~\cite{FW, BFP}.
\begin{equation}
\label{eq:quant}
\langle s_i\rangle*\langle w\rangle=\sum_{\alpha\in\Phi^+}(\omega_i,\alpha^\vee)\,\langle Q_{\alpha}(w)\rangle\,,
\end{equation} 
where $s_i$ is a simple reflection  and $*$ denotes
the product in $QH^*(G/B)$. 

Based on Corollary~\ref{cor:K-hyperplane-sections}
and~(\ref{eq:quant}), we formulate the following conjecture. 

\begin{conjecture}
Fix a simple reflection $s_i$.
Let $(\beta_1,\dots,\beta_l)$ be a $(-\omega_i)$-chain of roots.
Then we have
$$
[s_i]\circ [w]=[(1-(1- Q_{\beta_1})\cdots (1- Q_{\beta_l}))(w)]\,,
$$
where $\circ $ denotes the product in the ring $QK(G/B)$.
\end{conjecture}

The conjectured formula in $QK(G/B)$ 
specializes to Corollary~\ref{cor:K-hyperplane-sections},
upon setting $q_1=\cdots=q_r=0$. It also specializes to 
$QH$-Chevalley formula (\ref{eq:quant}), upon
taking the linear terms in the expansion of the operator $1-(1-
Q_{\beta_1})\cdots (1- Q_{\beta_l})$, cf.~Remark~\ref{classchev}.   
One can extend this conjecture to the quantum $T$-equivariant $K$-theory of 
$G/B$, see~\cite{Lee} for the definition of the ring $QK_T(G/B)$. In order to do this, one has to consider the operator $R_q^{[-\omega_i]}$ obtained from
$R^{[-\omega_i]}$ by replacing all Bruhat operators $B_\beta$
with the quantum Bruhat operators $Q_\beta$, cf.~Theorem~\ref{thm:K-Chevalley}.
It is not hard to extend the above conjecture to generalized partial flag
varieties $G/P$, as well.

\section{Appendix: Foldings of Galleries, LS-galleries, and LS-paths}
\label{sec:galleries}

In this appendix, we introduce admissible foldings of galleries,
and use this notion to reformulate our model for the characters
of the irreducible representations (Corollary~\ref{cor:model-for-character}) 
and for the Demazure characters 
(Corollary~\ref{cor:model-for-demazure-character}). 
For regular weights, admissible foldings of galleries are 
similar, but not equivalent, to the LS-galleries of 
Gaussent and Littelmann~\cite{GaLi}.
We clarify this relationship by showing that
it is based on Dyer's theorem~\cite{Dyer} 
about the EL-shellability of the Bruhat order.  
Then we compare the computational 
complexity of our model for characters with that of the model based on 
LS-paths and root operators.

\subsection{Admissible Foldings}

\begin{definition}
A {\it gallery\/} is a sequence 
$\gamma=(F_0,A_0,F_1,A_1, F_2, \dots , F_l, A_l, F_{l+1})$ 
such that $A_0,\dots,A_l$ are alcoves;
$F_j$ is a codimension one common face of the alcoves $A_{j-1}$ and $A_j$,
for $j=1,\dots,l$; $F_0$ is a vertex of the first alcove $A_0$; and 
$F_{l+1}$ is a vertex of the last alcove $A_l$. 
Furthermore, we require that $F_0=\{0\}$ and 
$F_{l+1}=\{\mu\}$ for some weight $\mu\in\Lambda$,
which is called the {\it weight\/} of the gallery.
We say that a gallery is {\it unfolded\/} if $A_{j-1}\ne A_j$, 
for $j=1,\dots,l$.  
\label{def:gallery}
\end{definition}

These galleries are special cases of the generalized galleries in~\cite{GaLi}.

In this subsection, we will consider only galleries such that $A_0=A_\circ$ is
the fundamental alcove.  Unfolded galleries of weight $\mu$ with
$A_0=A_\circ$ are in one-to-one correspondence with alcove paths
$(A_\circ,\dots,A_l)$ such that $\mu\in A_l$.  Indeed, $F_j$ should be the
unique common wall of two adjacent alcoves $A_{j-1}$ and $A_j$, for
$j=1,\dots,l$.

\begin{definition}
Let us say that a gallery $\gamma$ of weight $\mu$ is {\it reduced\/} if 
$A_0=A_\circ$,  and  $\gamma$ has has minimal length among all 
galleries of weight $\mu$ with $A_0=A_\circ$.
Clearly, every reduced gallery is unfolded.
\label{def:reduced-gallery}
\end{definition}

\begin{lemma}
Let $\lambda$ be a dominant weight.
Then the last alcove in a reduced gallery of 
weight $-\lambda$ is $A_l = A_{-\lambda}$.
Hence, reduced galleries with an antidominant
weight $-\lambda$ are in one-to-one correspondence with reduced 
alcove paths from $A_\circ$ to $A_{-\lambda}$, 
which, in turn, correspond to reduced 
decompositions of $v_{-\lambda}\in\Waff$.
\end{lemma}

\begin{proof}
The number of hyperplanes $H_{\alpha,k}$ that separate the point 
$E=\{-\lambda\}$ from the fundamental alcove $A_\circ$
is $m=\sum_{\alpha\in\Phi^+} (\lambda,\alpha^\vee)$.
Thus, the length of any alcove path from $A_\circ$ to an alcove
$A_l$ with vertex $E$ should be at least $m$.  The number $m$ is
precisely the length of a reduced alcove path from $A_\circ$
to $A_{-\lambda}$.
On the other hand, for any other alcove $A' \ne A_{-\lambda}$
such that $E$ is a vertex of $A'$, the number of hyperplanes that separate $A'$
from $A_\circ$ is strictly greater than $m$.
\end{proof}

For a gallery 
$\gamma=(F_0,A_0,F_1,\dots,F_l, A_l, F_{l+1})$, 
let $r_1,\dots,r_l\in \Waff$ denote the affine reflections with respect  
to the affine hyperplanes containing the faces $F_1,\dots,F_l$.
For $j=1,\dots,l$, let the $j$-th
{\it tail-flip operator\/} $f_j$ be the operator that sends the gallery 
$\gamma=(F_0,A_0, F_1, \dots,F_l, A_l, F_{l+1})$ to
the gallery $f_j(\gamma)$ given by 
$$
f_j(\gamma):=(F_0,A_0, F_1, A_1, \dots, A_{j-1},
 F_j'=F_j,  A_{j}', F_{j+1}', A_{j+1}', \dots,  A_l', F_{l+1}'),
$$
where $A_i' := r_j(A_i)$ and $F_i':=r_j(F_i)$, for 
$i=j,\dots,l+1$.
In other words, the operator $f_j$ leaves the initial segment 
of the gallery from $A_0$ to $A_{j-1}$ intact 
and reflects the remaining tail by $r_j$.
Clearly, the operators $f_j$ commute.
Hence, they determine an action of the group $(\Z/2\Z)^l$ on
galleries.
Every gallery is obtained from an unfolded gallery by applying several
tail-flips.
Equivalently, using the operators $f_j$, one can always transform (unfold) 
an arbitrary gallery into a uniquely defined unfolded gallery.

\begin{lemma}
If $\gamma$ is a gallery of weight $\mu$, then
$f_{j_1} \cdots f_{j_s}(\gamma)$ is a gallery of weight
$r_{j_1}\cdots r_{j_s}(\mu)$,
for any $1\leq j_1<\cdots< j_s\leq l$. 
\label{sec:weight-of-tail-flipping-reflection}
\end{lemma}

\begin{proof}
First, let us apply $f_{j_s}$ to $\gamma$.  We obtain a gallery
of weight $r_{j_s}(\mu)$.  Applying the tail-flip
$f_{j_{s-1}}$ to $f_{j_s}(\gamma)$ changes its weight to  
$r_{j_{s-1}} r_{j_s}(\mu)$, etc.
\end{proof}

\begin{definition} 
Let $\gamma$ be an unfolded gallery, and let 
$r_1,\dots,r_l$ be the affine reflections with respect 
to the faces of $\gamma$.  
An {\it admissible folding\/} of $\gamma$
is a gallery of the form $f_{j_1}\cdots f_{j_s} (\gamma)$ 
for some $1\leq j_1<\cdots< j_s\leq l$  such that
$$
1\lessdot \bar r_{j_1} \lessdot \bar r_{j_1} \bar r_{j_2} \lessdot 
\cdots \lessdot \bar r_{j_1} \bar r_{j_2} \cdots \bar r_{j_s}
$$
is a saturated increasing chain in the Bruhat order on the Weyl group $W$.
More generally, for $u\in W$,  a {\it $u$-admissible folding\/} of $\gamma$ 
is a gallery
of the form $f_{j_1}\cdots f_{j_s} (\gamma)$ 
for some $1\leq j_1<\cdots< j_s\leq l$  such that
$$
u \gtrdot u\,\bar r_{j_1} \gtrdot u\,\bar r_{j_1} \bar r_{j_2} \gtrdot 
\cdots \gtrdot u\,\bar r_{j_1} \bar r_{j_2} \cdots \bar r_{j_s} 
$$
is a saturated decreasing chain in the Bruhat order on 
the Weyl group $W$.
We allow $s=0$, so the gallery $\gamma$ itself is an admissible 
($u$-admissible) folding of $\gamma$.
Notice that admissible foldings are precisely $w_\circ$-admissible foldings.
\label{def:u-admissible}
\end{definition}

We can also give the following intrinsic characterization of $u$-admissible
foldings.

\begin{lemma}  Let 
$\gamma'=(A_0', F_1', \dots,  F_l', A_l', E')$ be a gallery, and
$r_1',\dots,r_l'$ be the affine reflections with respect to 
the faces $F_1',\dots,F_l'$.  Let $\{j_1<\cdots<j_s\} :=
\{j\in\{1,\dots,l\} \mid A_{j-1}'= A_j'\}$.
Then the gallery $\gamma'$ is a $u$-admissible folding of some unfolded
gallery $\gamma$ if and only if
$$
u^{-1} \gtrdot \bar r_{j_1}'\,u^{-1} \gtrdot 
\bar r_{j_1}' \bar r_{j_2}'\,u^{-1} \gtrdot 
\cdots \gtrdot \bar r_{j_1}' \bar r_{j_2}' \cdots \bar r_{j_s}'  \,u^{-1}
$$
is a saturated decreasing chain in the Bruhat order on 
the Weyl group $W$.
\label{lem:intrinsic}
\end{lemma}

\begin{proof}
We have $\gamma'= f_{j_1}\cdots f_{j_s} (\gamma)$. 
Let $r_1,\dots,r_l$ be the reflections with respect to the faces of
the unfolded gallery $\gamma$.  Then
$$
r_{j_1}' = r_{j_1}, \
r_{j_2}' = r_{j_1} r_{j_2} r_{j_1}, \ 
r_{j_3}' = r_{j_1} r_{j_2} r_{j_3} r_{j_2} r_{j_1}, \  \dots
$$
This implies $r_{j_1}' r_{j_2}' \cdots r_{j_i}' =
(r_{j_1} r_{j_2} \cdots r_{j_i})^{-1}$, for $i=1,\dots,s$.
Now the lemma follows from Definition~\ref{def:u-admissible}.
\end{proof}

Corollaries~\ref{cor:model-for-demazure-character} 
and~\ref{cor:model-for-character} are equivalent to the following
claim.  Let $\weight(\gamma)$ denote the weight of a gallery $\gamma$.

\begin{corollary}
\label{cor:admissible-foldings}
Let $\lambda$ be a dominant weight, and   
let $\gamma$ be a reduced gallery with $\weight(\gamma)=-\lambda$.

\noindent
{\rm(1)}
The character $ch(V_{\lambda})$ is equal to the sum
$$
ch(V_{\lambda}) = \sum_{\gamma'} e^{-\mathrm{weight}(\gamma')}
$$
over all admissible foldings $\gamma'$ of the gallery $\gamma$.

\noindent
{\rm(2)}
Let $u\in W$.  The Demazure character $ch(V_{\lambda,u})$ 
is equal to the sum
$$
ch(V_{\lambda,u}) = \sum_{\gamma'} e^{-u(\mathrm{weight}(\gamma'))}
$$
over all $u$-admissible foldings $\gamma'$ of the gallery $\gamma$.
\end{corollary}

\subsection{LS-galleries}

In this section, we discuss the relationship between admissible foldings
and LS-galleries of Gaussent and Littelmann in case of a {\it regular\/} 
weight $\lambda$.  We show that LS-galleries can be associated with 
admissible foldings of some {\it special\/} reduced galleries.

We start by recalling some terminology from \cite{GaLi}. 
Let us fix a dominant regular weight $\lambda$.
Let us say that a gallery $\gamma$ of weight $\lambda$ is
{\it minimal\/} if $\gamma$ crosses only the hyperplanes strictly separating 
$0$ and $\lambda$. Note that in such a gallery
we have $A_0=A_\circ$, and the last alcove $A_l$ is
$w_\circ(A_\circ)+\lambda=-A_\circ+\lambda$.

Recall that the facets of the fundamental alcove are
$H_i = H_{\alpha_i,0}$, for $i=1,\dots,r$; and $H_0=H_{\alpha_0,-1}$.
If $F$ is a face of the fundamental alcove $A_\circ$, we define its {\em type}
by 
$$
\mathrm{type}(F)=\{i\mid F\subset
H_i,\,i=0,1,\ldots,r\}\,.
$$ 
For instance, $\mathrm{type}(\{0\})=\{1,\dots,r\}$ 
and $\mathrm{type}(A_\circ)=\emptyset$.  For an arbitrary face $F$, 
its type is defined as ${\rm type}(F')$, where $F'$ is the
unique face of $A_\circ$ such that $F=w(F')$ for some $w$ in $\Waff$. The type
of a gallery $\gamma=(F_0,A_0,F_1,\dots,A_l,F_{l+1})$ is defined as
$\mathrm{type}(\gamma)=(\mathrm{type}(F_0),\mathrm{type}(A_0),\dots,
\mathrm{type}(F_{l+1}))$.

For a gallery $\gamma=(F_0,A_0,F_1,\dots, A_l,F_{l+1})$, let
$\{j_1<\ldots<j_s\}=\{j\mid A_{j-1}=A_j\}$, and let $r_j$  be the reflections
with respect to the hyperplanes containing the faces $F_j$.  The {\it
companion\/} of $\gamma$ is the sequence $(u_0,\dots,u_s)$ of elements in $W$,
where $u_{0}\in W$ is the unique element such that $u(A_\circ) = A_0$; and
$u_{i} = \bar r_{j_i} u_{i-1}$, for $i=1,\dots,s$.

\begin{definition} {\rm\cite{GaLi} }
For a minimal gallery $\gamma$ of a (dominant regular) weight $\lambda$,
the set $\Gamma_{LS}(\gamma)$ of {\it LS-galleries\/} associated 
with $\gamma$ is the set 
of all galleries $\gamma'$ such that
(1) $\mathrm{type}(\gamma')=\mathrm{type}(\gamma)$; and 
(2) the companion $(u_0,\dots,u_s)$ of $\gamma'$ 
is a saturated decreasing chain in the Bruhat order on $W$. 
\end{definition}

The general definition of LS-galleries given is~\cite{GaLi} 
for arbitrary dominant weights $\lambda$ is more complicated.
They are defined as certain collections of faces of alcoves that satisfy
several conditions, including some positivity and dimension conditions.
The companion of such a gallery is a chain in the Bruhat order on the 
quotient $W/W_\lambda$.  For regular weights, the definition of 
LS-galleries from~\cite{GaLi} is equivalent to the 
simplified definition above.

It was shown in~\cite{GaLi} that, for a minimal gallery $\gamma$ of weight
$\lambda$,
$$
ch(V_\lambda) = \sum_{\gamma'\in \Gamma_{LS}(\gamma)} 
e^{\mathrm{weight}(\gamma')}.
$$

Let us now clarify the relationship between 
Corollary~\ref{cor:admissible-foldings}.(1) and this statement.

Let us say that a gallery of $\gamma=(F_0,A_0,F_1,\dots,A_l,F_{l+1})$ 
is {\it special\/} if $l\geq N=|\Phi^+|$ (the number of positive roots)
and all alcoves $A_0, \dots,A_N$ 
and faces $F_1,\dots,F_N$ are adjacent to the origin $0$.
Let us define the transformation 
$$
t\,:\,\{\textrm{special galleries of weight $-\mu$}\}\longrightarrow
\{\textrm{galleries of weight $\mu$}\}.
$$
For a special gallery $\gamma=(F_0,A_0,F_1,\dots,A_l,F_{l+1})$ 
of weight $-\mu$,
the gallery $t(\gamma)$ is defined as follows:
(1) remove the first $N$ alcoves $A_0,\dots,A_{N-1}$
from the gallery $\gamma$ together with the faces $F_1,\dots,F_{N}$;
(2) translate all remaining alcoves and faces by the weight $\mu$; 
(3) reverse the sequence of alcoves and faces in the gallery.
In other words,
$$
t\,:\,(F_0,A_0,\dots,A_l,F_{l+1})\longmapsto (F_{l+1} + \mu, A_l + \mu,\dots,
F_{N+1}+\mu,A_N+\mu, F_0+\mu),
$$

If $\gamma=(F_0,A_0,F_1,\dots,A_l,F_{l+1})$ is a special reduced gallery
of weight $-\lambda$ (Definition~\ref{def:reduced-gallery}),
then $A_N=w_\circ(A_\circ)$ and $F_i\subset H_{\beta_i,0}$, 
for $i=1,\ldots,N$.   All foldings of $\gamma$ are also special.
The image $t(\gamma)$ of $\gamma$ is a minimal gallery of weight $\lambda$.
Moreover, all minimal galleries are of this form.
Notice that, for a regular weight $\lambda$, we can always find a special
reduced gallery of weight $-\lambda$.

\begin{proposition} 
Let $\gamma$ be a special reduced gallery of weight $-\lambda$, where
$\lambda$ is a regular weight.
Then the map $\gamma'\mapsto t(\gamma')$ is a bijection 
between the set of admissible foldings of $\gamma$ and 
the set $\Gamma_{LS}(t(\gamma))$ of LS-galleries associated with $t(\gamma)$. 
Moreover, we have ${\rm weight}(t(\gamma'))=-{\rm weight}(\gamma')$. 
\label{prop:admissible=LS}
\end{proposition}

The proof of this proposition is based on the following fundamental 
(and nontrivial) result, 
which expresses the {\it EL-shellability\/} of the Bruhat 
order on a Weyl group,
and is closely related to the Verma theorem~\cite{Ver}.
This result was proved for an arbitrary Coxeter group 
in~\cite[Proposition 4.3]{Dyer}.  We also refer to
\cite[Theorem 6.4]{BFP} for a new approach and a different generalization.
Recall that {\em reflection orderings\/} \cite{Hum, Dyer} 
are total orders on roots in $\Phi^+$ that are
associated with reduced decompositions $w_\circ=s_{i_1}\ldots s_{i_N}$
for $w_\circ$, as follows:
$$
\alpha_{i_N} < s_{i_N}(\alpha_{i_{N-1}}) < \ldots <
s_{i_N}s_{i_{N-1}}\ldots s_{i_2}(\alpha_{i_{1}})\,.
$$

\begin{proposition}\cite{Dyer, BFP}\label{uniquedec} Fix a reflection ordering
$\beta_1<\cdots <\beta_N$. For any Weyl group element $w$, there is a unique
saturated increasing chain in Bruhat order from $1$ to $w$ of the form
\begin{equation}\label{equniquedec} 1\lessdot s_{\beta_{j_1}}\lessdot
s_{\beta_{j_1}}s_{\beta_{j_2}}\lessdot\ldots\lessdot s_{\beta_{j_1}}\ldots
s_{\beta_{j_p}}=w\,,\end{equation} where $1\le j_1<\ldots<j_p\le N$.
\end{proposition}

\begin{proof}[Proof of Proposition~\ref{prop:admissible=LS}]
Let $\gamma'$ be an arbitrary admissible folding of $\gamma$. Every tail-flip
operator $f_j$ preserves the type of $\gamma'$, that is,
$\mathrm{type}(\gamma')=\mathrm{type}(f_j(\gamma'))$, and changes its weight by
a multiple of a root.  Hence, the transformation $t$ applied to $\gamma'$ can
be viewed as a composition of the translation by $\lambda$ with a translation
by an element of the root lattice. Note that the second translation is an
element of $\Waff$. Recalling that $\gamma$ is mapped to $t(\gamma)$ via the
translation by $\lambda$, we conclude that the gallery $t(\gamma')$ has the
same type as $t(\gamma)$.

Let us now examine the companion of $t(\gamma')$. Let $r_1,\ldots,r_l$ and
$r_1',\ldots,r_l'$ be the affine reflections with respect to the faces of
$\gamma$ and $\gamma'$, respectively. Let $p$ be such that $j_p\le N$ and
$j_{p+1}>N$. Assume that $\gamma' = 
f_{j_1} \cdots f_{j_s}(\gamma)$, where $j_1<\cdots<j_s$,
so 
$$
1\lessdot \bar r_{j_1} \lessdot \bar r_{j_1} \bar r_{j_2} \lessdot 
\cdots \lessdot \bar r_{j_1} \bar r_{j_2} \cdots \bar r_{j_s}
$$
is a saturated decreasing chain in the Bruhat order.
The companion of $t(\gamma')$ is the sequence
$$
(u_0=\bar{r}_{j_1}\ldots\bar{r}_{j_s},\,\bar{r}_{j_s}'u_0,\,\bar{r}_{j_{s-1}}'\bar{r}_{j_s}'u_0,\,\ldots,\,\bar{r}_{j_{p+1}}'\ldots\bar{r}_{j_s}'u_0)\,.
$$
But since $r_{j_1}' r_{j_2}' \cdots r_{j_i}' = (r_{j_1} r_{j_2} \cdots
r_{j_i})^{-1}$, for $i=1,\dots,s$ (see the proof of Lemma \ref{lem:intrinsic}),
the companion of $t(\gamma')$ is the sequence
$$
(\bar{r}_{j_1}\ldots\bar{r}_{j_s},\,\bar{r}_{j_1}\ldots\bar{r}_{j_{s-1}},\,\ldots,\,\bar{r}_{j_1}\ldots\bar{r}_{j_p})\,,
$$
which is a saturated decreasing chain in Bruhat order. 
We have thus shown
that the image of map $t$ is contained in $\Gamma_{LS}(\gamma)$. 

It suffices
to construct the inverse map. Recall that the first $N$ faces $F_i$ of $\Gamma$
satisfy $F_i\subset H_{\beta_i,0}$. This gives a reflection ordering
$\beta_1<\cdots<\beta_N$, according
to Lemma \ref{lem:chains=decompostions}. Given a gallery $\gamma''$ in
$\Gamma_{LS}(t(\gamma))$, assume that its companion ends at some $w$ in $W$.
According to Proposition~\ref{uniquedec}, there is a unique way of writing
$w=s_{\beta_{j_1}}\ldots s_{\beta_{j_p}}$ for $1\le j_1<\ldots<j_p\le N$, such
that (\ref{equniquedec}) holds. 

Let us now relabel the faces of $\gamma''$ as follows:
$(F_{l+1}',A_l',F_l',A_{l-1}',F_{l-1}',\ldots)$. Let
$\{j_{p+1}<\ldots<j_s\}=\{j \mid A_{j-1}'=A_j'\}$. We associate with $\gamma''$
the gallery $f_{j_1}\ldots f_{j_p}f_{j_{p+1}}\ldots f_{j_s}(\gamma)$. 
The facts stated above imply that this construction gives
the inverse map to $t$.  
\end{proof}

\begin{remark}\label{relationgalleries} (i) 
For a nonregular weight $\lambda$, it is not clear how to associate
LS-galleries with our admissible foldings.

\noindent
(ii)  According to~\cite{GaLi}, one can associate a collection of 
continuous piecewise-linear Littelmann paths with the set of LS-galleries 
$\Gamma_{LS}(\gamma)$ by connecting the centers of the lower dimensional faces in the galleries.  
In~\cite{LP}, we show that we do not obtain Littelmann paths by applying the same procedure (or similar ones) to our model.
\end{remark}

\subsection{Comparison of Computational Complexities}

We conclude with a comparison between the computational complexities of our
construction and the construction of LS-paths based on root operators.

Fix a root system of rank $r$ with $N$ positive roots, a dominant weight
$\lambda$, and a Weyl group element $u$ of length $l$. We want to determine the
character of the Demazure module $V_{\lambda,u}$. Let $d$ be its dimension, and
let $L$ be the length of the affine Weyl group element $v_{-\lambda}$ (that is,
the number of affine hyperplanes separating the fundamental alcove $A_\circ$
and $A_\circ-\lambda$). Note that $L=2(\lambda, \rho^\vee)$, where
$\rho^\vee=\frac{1}{2}\sum_{\beta\in\Phi^+}\beta^\vee$. We claim that the
complexity of our character formula is $O(d\, l L)$. Indeed, we start by
determining an alcove path via the method described at the end of
Section~\ref{sec:K-chevalley-formula}, which involves
sorting a sequence of $L$ rational numbers. The complexity is $O(L\,\log L)$,
and note that $\log L$ is, in general, much smaller than $d$ (see below for
some examples). Whenever we examine some subword of the word of length $L$ we
fixed at the beginning, we have to check at most $L-1$ ways to add an extra
reflection at the end. On the other hand, in each case, we have to check
whether, upon multiplying by the corresponding nonaffine reflection, the length
decreases by precisely 1. The complexity of the latter operation is $O(l)$,
based on the Strong Exchange Condition \cite[Theorem~5.8]{Hum}. Then, for each
``good'' subword, we have to do a calculation, namely applying at most $2l$
affine reflections to $-\lambda$. In fact, it is fairly easy to implement this
algorithm.

Now let us examine at the complexity of the algorithm 
based on root operators for constructing the LS-paths associated with $\lambda$.
In other words, we
are looking at the complexity of constructing the corresponding crystal graph.
We have to generate the whole crystal graph first, and then figure out which
paths give weights for the Demazure module. For each path, we can apply $r$
root operators. Each path has at most $N$ linear steps, so applying a root
operator has complexity $O(N)$. But now we have to check whether the result is
a path already determined, so we have to compare the obtained path with the
other paths (that were already determined) of the same rank in the crystal
graph (viewed as a ranked poset). This has complexity $O(N M)$, where $M$ is
the maximum number of elements of the same rank. Since we have at most $N+1$
ranks, $M$ is at least $d/(N+1)$. In conclusion, the complexity is $O(d r N
M)$, which is at least $O(d^2 r)$.

Let us get a better picture of how the two results compare. Assume we are in a
classical type, and let us first take $\lambda$ to be the $i$-th fundamental
weight, with $i$ fixed, plus $u=w_\circ$. Clearly $l$ is $O(r^2)$, $L$ is
$O(r)$, and $d$ is $O(r^i)$, so the complexity of our formula is $O(r^{i+3})$.
For LS-paths, we get at least $O(r^{2i+1})$. So the ratio between the
complexity in the model based on LS-paths and our model is at least
$O(r^{i-2})$.

Let us also take $\lambda=\rho$. In this case $d=2^N$, and a simple calculation
shows that $L$ is $O(r^3)$. Our formula has complexity $O(2^N r^5)$, while the
model based on LS-paths has complexity at least $O(2^{2N} r)$. So the ratio
between the complexities is at least $O(2^N/r^4)$, where $N$ is $r(r+1)/2$,
$r^2$, and $r^2-r$ in types $A$, $B/C$, and $D$, respectively.



\begin{thebibliography}{FuWo}

\bibitem[BFP]{BFP}
F.~Brenti, S.~Fomin, and A.~Postnikov.
\newblock Mixed {B}ruhat operators and {Y}ang-{B}axter equations for {W}eyl
  groups.
\newblock {\em Internat. Math. Res. Notices}, 8:419--441, 1999.

\bibitem[Brion]{Bri}
M.~Brion.
\newblock Positivity in the {G}rothendieck group of complex flag varieties.
\newblock {\em J. Algebra}, 258:137--159, 2002.

\bibitem[BrLa]{BL} 
M.~Brion and V.~Lakshmibai.
\newblock A geometric approach to standard monomial theory.
\newblock {\em Represent. Theory}, 7:651-680, 2003.

\bibitem[Cher]{Cher}
I.~Cherednik.
\newblock Quantum {K}nizhnik-{Z}amolodchikov equations and affine root systems.
\newblock {\em Comm. Math. Phys.}, 150:109--136, 1992.

\bibitem[Chev]{Chev}
C.~Chevalley.
\newblock Sur les d\mbox{\'{e}}compositions cellulaires des espaces
  \mbox{$G/B$}.
\newblock In {\em Algebraic \mbox{G}roups and \mbox{G}eneralizations:
  \mbox{C}lassical \mbox{M}ethods}, volume 56 Part 1 of {\em Proceedings and
  Symposia in Pure Mathematics}, pages 1--23. Amer. Math. Soc., 1994.

\bibitem[Dem]{Dem}
M.~Demazure.
\newblock D\mbox{\'e}singularization des vari\mbox{\'et\'e}s de
  \mbox{S}chubert.
\newblock {\em Annales E.N.S.}, 6:53--88, 1974.

\bibitem[Deo1]{Deo1}
V.~V. Deodhar.
\newblock Some characterizations of {B}ruhat ordering on a {C}oxeter group and
  determination of the relative {M}\"obius function.
\newblock {\em Invent. Math.}, 39:187--198, 1977.

\bibitem[Deo2]{Deo2}
V.~V. Deodhar.
\newblock A splitting criterion for the {B}ruhat orderings on {C}oxeter groups.
\newblock {\em Comm. Algebra}, 15:1889--1894, 1987.

\bibitem[Dyer]{Dyer}
M.~J. Dyer.
\newblock Hecke algebras and shellings of {B}ruhat intervals.
\newblock {\em Compositio Math.}, 89:91--115, 1993.


\bibitem[FGP]{FGP}
S.~Fomin, S.~Gelfand, and A.~Postnikov.
\newblock Quantum \mbox{S}chubert polynomials.
\newblock {\em J. Amer. Math. Soc.}, 10:565--596, 1997.


\bibitem[FuLa]{FL} 
W.~Fulton, A.~Lascoux, 
\newblock A Pieri formula in the Grothendieck ring of a flag bundle.
\newblock {\em Duke Math. J.}, 76:711--729, 1994.

\bibitem[FuWo]{FW}
W.~Fulton and C.~Woodward.
\newblock On the quantum product of {S}chubert classes.
\newblock {\em J. Algebraic Geom.}, 13:641--661, 2004. 

\bibitem[GaLi]{GaLi}
S.~Gaussent and P.~Littelmann.
\newblock {LS}-galleries, the path model and {MV}-cycles.
\newblock {\em Duke Math. J.},  127:35--88, 2005.

\bibitem[GiLe]{GiLe}
A.~Givental and Y.-P. Lee.
\newblock Quantum {$K$}-theory on flag manifolds, finite-difference {T}oda
  lattices and quantum groups.
\newblock {\em Invent. Math.}, 151:193--219, 2003.

\bibitem[GrRa]{GR}
S.~Griffeth and A.~Ram.
\newblock Affine {H}ecke algebras and the {S}chubert calculus.
\newblock  {\em European J. Combin.}, 25:1263--1283, 2004.

\bibitem[Hum]{Hum}
J.~E.~Humphreys.
\newblock {\em Reflection groups and {C}oxeter groups}, volume~29 of {\em
  Cambridge Studies in Advanced Mathematics}.
\newblock Cambridge University Press, Cambridge, 1990.

\bibitem[Kost]{Kost}
B.~Kostant.
\newblock Powers of the {E}uler product and commutative subalgebras of a
  complex simple {L}ie algebra.
\newblock{\em Invent. Math.}, 158:181--226, 2004.

\bibitem[KoKu]{KK}
B.~Kostant and S.~Kumar.
\newblock \mbox{{\it T}}-equivariant \mbox{{\it K}}-theory of generalized flag
  varieties.
\newblock {\em J. Diff. Geom.}, 32:549--603, 1990.

\bibitem[Kum]{Kum}
S.~Kumar.
\newblock Private communication.

\bibitem[LaLi]{LL}
V.~Lakshmibai and P.~Littelmann.
\newblock Richardson varieties and equivariant {$K$}-theory.
\newblock {\em J. Algebra}, 260:230--260, 2003.

\bibitem[LLM]{LLM}
V.~Lakshmibai, P.~Littelmann, and P.~Magyar.
\newblock Standard monomial theory and applications.
\newblock In {\em Representation theories and algebraic geometry (Montreal, PQ,
  1997)}, volume 514 of {\em NATO Adv. Sci. Inst. Ser. C Math. Phys. Sci.},
  pages 319--364. Kluwer Acad. Publ., Dordrecht, 1998.

\bibitem[LaSe]{LS1}
V.~Lakshmibai and C.~S. Seshadri.
\newblock Standard monomial theory.
\newblock In {\em Proceedings of the Hyderabad Conference on Algebraic Groups
  (Hyderabad, 1989)}, pages 279--322, Madras, 1991. Manoj Prakashan.

\bibitem[LaSc]{LS2}
A.~Lascoux and M.-P. Sch\mbox{\"{u}}tzenberger.
\newblock Structure de \mbox{H}opf de l'anneau de cohomologie et de l'anneau de
  \mbox{G}rothendieck d'une vari\mbox{\'et\'e} de drapeaux.
\newblock {\em C. R. Acad. Sci. Paris \mbox{S\'e}r. I Math.}, 295:629--633,
  1982.

\bibitem[Lee]{Lee}
Y.-P. Lee.
\newblock Quantum {$K$}-theory {I}: {F}oundations.
\newblock  {\em Duke Math. J.},  121:389--424, 2004.

\bibitem[Len]{Len}
C.~Lenart.
\newblock A \mbox{$K$}-theory version of {M}onk's formula and some related
  multiplication formulas.
\newblock {\em J. Pure Appl. Algebra}, 179:137--158, 2003.

\bibitem[LePo]{LP}
C.~Lenart and A.~Postnikov.
\newblock A combinatorial model for crystals of {K}ac-{M}oody algebras.
\newblock {\tt arXiv:math.RT/0502147}.

\bibitem[LeSo]{lasptf}
C.~Lenart and F.~Sottile.
\newblock {A Pieri-type formula for the $K$-theory of a flag manifold}.
\newblock {\tt arXiv:math.CO/0407412}, to appear in {\em Trans. Amer. Math. Soc.}

\bibitem[Lit1]{Li1}
P.~Littelmann.
\newblock {A Littlewood-Richardson rule for symmetrizable Kac-Moody algebras}.
\newblock {\em Invent. Math.}, 116:329--346, 1994.

\bibitem[Lit2]{Li2}
P.~Littelmann.
\newblock {Paths and root operators in representation theory}.
\newblock {\em Ann. of Math. (2)}, 142:499--525, 1995.

\bibitem[Lit3]{Li3}
P.~Littelmann.
\newblock {Contracting modules and standard monomial theory for symmetrizable {K}ac-{M}oody algebras}.
\newblock {\em J. Amer. Math. Soc.}, 11:551--567, 1998.

\bibitem[LiSe]{LS3}
P.~Littelmann and C.~S. Seshadri.
\newblock A {P}ieri-{C}hevalley type formula for {$K(G/B)$} and standard
  monomial theory.
\newblock In {\em Studies in memory of Issai Schur (Chevaleret/Rehovot, 2000)},
  volume 210 of {\em Progr. Math.}, pages 155--176. Birkh\"auser Boston,
  Boston, MA, 2003.


\bibitem[Mat]{Mat}
O.~Mathieu.
\newblock Positivity of some intersections in {$K\sb 0(G/B)$}.
\newblock {\em J. Pure Appl. Algebra}, 152:231--243, 2000.

\bibitem[PiRa]{PR}
H.~Pittie and A.~Ram.
\newblock A \mbox{P}ieri-\mbox{C}hevalley formula in the \mbox{$K$}-theory of a
  \mbox{$G/B$}-bundle.
\newblock {\em Electron. Res. Announc. Amer. Math. Soc.}, 5:102--107, 1999.

\bibitem[Ver]{Ver}
D.-N.~Verma. 
\newblock M\"obius inversion for the Bruhat ordering on a Weyl group.
\newblock {\em Ann. Sci. École Norm. Sup. (4)}, 4:393--398, 1971. 
\end{thebibliography}

\end{document}